%% file: main.tex
\documentclass[11pt,twoside,reqno]{amsart}
\usepackage[utf8]{inputenc}
\usepackage[T1]{fontenc}
\usepackage{import}
\usepackage{newclude}
\include{preamble}

\title{Quantative universality for cokernels of matrices with symmetries}
\author{Jiahe Shen}

\date{\today}

\begin{document}

\thanks{I thank Mehtaab Sawhney for suggesting this project and for many helpful discussions, and I also thank Roger Van Peski for numerous valuable conversations and suggestions. This work was supported by NSF grant DMS-2246576 and Simons Investigator grant 929852.}

\maketitle

\begin{abstract} 
We prove universality for cokernels of random integral matrices with symmetries via an approach different from the classical surjection moment method introduced by Wood \cite{wood2017distribution}. In the symmetric case, we reprove Hodges’ universality theorem \cite{hodges2024distribution}, i.e. the version incorporating the canonical pairing from Wood’s setting, and in the alternating case we reprove the local universality theorem of Nguyen-Wood \cite{nguyen2025local}. A key advantage of our method is that it is quantitative: we obtain explicit error bounds, which are exponentially small in most regimes, thereby addressing Wood’s question on effective convergence rates raised in \cite[Section 1.6]{wood2017distribution}. Our argument is inspired by Maples’ exposure-process and coupling viewpoint \cite{maples2013cokernels,Maples_symma_2013} and uses a generalized form of Fourier-analytic estimates in the exponentially sharp style of Ferber-Jain-Sah-Sawhney \cite{ferber2023random}.
\end{abstract}



\tableofcontents

\include*{Introduction}

\include*{Preliminaries}

\include*{The_symmetric_pairing_on_the_cokernel_of_a_nonsingular_symmetric_matrix}

\include*{The_cokernel_of_an_alternating_matrix_and_its_torsion_subgroup}

\include*{The_non_sparse_assumptions}

\include*{Proof_of_the_symmetric_case}

\include*{Proof_of_the_alternating_case}

\include*{Index_of_notations}


\input{main.bbl}
\end{document}

%% file: preamble.tex
\usepackage{amsmath,amssymb,amsthm,mathtools}
\usepackage[toc,page]{appendix}
\usepackage{longtable,booktabs,array}
\usepackage{anysize}
\usepackage{amscd}
\usepackage{stmaryrd}
\usepackage{wasysym}
\usepackage{mathrsfs}
\usepackage[scr=boondox]{mathalfa}
\usepackage{enumitem}
\usepackage{accents}
\usepackage{dsfont}
\usepackage{soul}
\usepackage{color}
\usepackage[all]{xy}
\usepackage{float}
\usepackage{bm}
\usepackage{makecell}
\usepackage{thmtools}      
\usepackage{setspace}
\usepackage{comment}
\usepackage{tikz}
\usepackage{tikz-cd}
\usepackage{mathdots}
\usepackage{mathtools}
\usepackage{graphicx}
\usepackage[myheadings]{fullpage}

\usepackage[hypertexnames=false,colorlinks,linkcolor=black,citecolor=black]{hyperref}
\usepackage[capitalize,nameinlink]{cleveref}

\crefname{appendix}{Appendix}{Appendices}
\Crefname{appendix}{Appendix}{Appendices}

\numberwithin{equation}{section}
\newcommand\mtop{.95in}
\newcommand\mbottom{.95in}
\newcommand\mleft{1in}
\newcommand\mright{1in}
\usepackage[top = \mtop, bottom = \mbottom, left = \mleft, right=\mright]{geometry}
\DeclareMathOperator{\val}{val}
\DeclareMathOperator{\Mat}{Mat}

\newtheorem{thm}{Theorem}[section]

\newtheorem{prop}[thm]{Proposition}

\newtheorem{lemma}[thm]{Lemma}

\newtheorem{cor}[thm]{Corollary}
\theoremstyle{definition}
\newtheorem{defi}[thm]{Definition}

\newtheorem{rmk}[thm]{Remark}

\newcommand\reallywidehat[1]{%
\savestack{\tmpbox}{\stretchto{%
  \scaleto{%
    \scalerel*[\widthof{\ensuremath{#1}}]{\kern-.6pt\bigwedge\kern-.6pt}%
    {\rule[-\textheight/2]{1ex}{\textheight}}
  }{\textheight}%
}{0.5ex}}%
\stackon[1pt]{#1}{\tmpbox}%
}
\DeclareSymbolFont{bbold}{U}{bbold}{m}{n}
\DeclareSymbolFontAlphabet{\mathbbold}{bbold}

\makeatletter
\def\@tocline#1#2#3#4#5#6#7{\relax
  \ifnum #1>\c@tocdepth 
  \else
    \par \addpenalty\@secpenalty\addvspace{#2}%
    \begingroup \hyphenpenalty\@M
    \@ifempty{#4}{%
      \@tempdima\csname r@tocindent\number#1\endcsname\relax
    }{%
      \@tempdima#4\relax
    }%
    \parindent\z@ \leftskip#3\relax \advance\leftskip\@tempdima\relax
    \rightskip\@pnumwidth plus4em \parfillskip-\@pnumwidth
    #5\leavevmode\hskip-\@tempdima
      \ifcase #1
       \or\or \hskip 1em \or \hskip 2em \else \hskip 3em \fi%
      #6\nobreak\relax
    \hfill\hbox to\@pnumwidth{\@tocpagenum{#7}}\par
    \nobreak
    \endgroup
  \fi}
\makeatother

\newcommand{\Z}{\mathbb{Z}}
\newcommand{\Q}{\mathbb{Q}}

\newcommand{\F}{\mathbb{F}}

\newcommand{\E}{\mathbb{E}}

\renewcommand{\l}{\lambda}

\newcommand{\Y}{\mathbb{Y}}

\DeclareMathOperator{\Len}{Len}

\DeclareMathOperator{\Alt}{Alt}

\DeclareMathOperator{\Sym}{Sym}

\DeclareMathOperator{\Cok}{Cok}

\DeclareMathOperator{\Aut}{Aut}

\DeclareMathOperator{\supp}{supp}
\DeclareMathOperator{\wt}{wt}
\DeclareMathOperator{\inv}{inv}
\DeclareMathOperator{\nil}{nil}
\DeclareMathOperator{\res}{res}
\DeclareMathOperator{\tors}{tors}
\DeclareMathOperator{\sym}{sym}
\DeclareMathOperator{\alt}{alt}
\DeclareMathOperator{\rank}{rank}
\DeclareMathOperator{\corank}{corank}

\DeclareMathOperator{\Dep}{Dep}
\DeclareMathOperator{\congsim}{\stackrel{\text{cong}}{\sim}}

\DeclareMathOperator{\Vol}{Vol}

\DeclareMathOperator{\GL}{GL}

\DeclareMathOperator{\Sp}{Sp}

\DeclareMathOperator{\diag}{diag}

\DeclareMathOperator{\Sur}{Sur}

%% file: Introduction.tex
\section{Introduction}\label{sec: intro}

A central theme in random matrix theory is universality: in high dimensions, many asymptotic
phenomena depend mainly on the symmetry class of the ensemble and only weakly on the fine details
of the entry distribution.  Traditionally, the primary objects of interest are eigenvalue
statistics---global laws, local spacing, and related spectral observables. When one passes from
random matrices over $\mathbb R$ or $\mathbb C$ to discrete models over $\mathbb Z$ (or over $p$-adic
and finite rings), eigenvalues become less directly accessible, and it is natural to seek arithmetic
invariants that play a comparable role. One such invariant is the cokernel: it records, in a discrete way, how far a matrix is from being invertible and refines rank information modulo primes, serving as an analogue of eigenvalue (or singular-value) data in the non-archimedean or discrete setting.
This analogy has been explicitly pointed out in several places in the literature; see, for instance, M\'esz\'aros \cite{meszaros2024phase}, Kovaleva \cite{kovaleva2023distribution}, Van Peski \cite{van2021limits,van2024local,van2023p}, and the ICM talk by Wood \cite{wood2022probability}.

For an integer matrix $Q_n\in \Mat_n(\Z)$, the cokernel is
\[
\Cok(Q_n)\coloneqq \Z^n/Q_n\Z^n,
\]
a finitely generated abelian group, finite precisely when $Q_n$ is nonsingular.
The symmetry class of $Q_n$ determines what extra canonical structure $\Cok(Q_n)$ carries, and hence
what the natural limiting object should be. In the absence of any bilinear pairing structure, one simply studies
random finite abelian groups, in line with the heuristics for class groups of number
fields conjectured by Cohen-Lenstra \cite{cohen2006heuristics}. In the symmetric setting, $\Cok(Q_n)$ comes with a natural
perfect symmetric pairing, a viewpoint motivated for instance by sandpile groups of random graphs
(cokernels of reduced Laplacians) and their canonical duality pairing, see Clancy-Kaplan-Leake-Payne-Wood \cite{clancy2015cohen}.
In the alternating (also called skew-symmetric or anti-symmetric) setting, one analogously encounters natural perfect alternating pairings,
matching Cohen-Lenstra type heuristics for elliptic curves in which Selmer and
Tate-Shafarevich groups carry the Cassels-Tate pairing raised by Poonen-Rains \cite{poonen2012random} and Bhargava-Kane-Lenstra-Poonen-Rains \cite{bhargava2015modeling}.

Motivated by these heuristics from different research fields, the study of cokernel distributions
for random integral matrices has attracted substantial
attention in recent years. A major breakthrough of Wood \cite{wood2017distribution} is a universality theorem for cokernels
of random symmetric integral matrices. Her result shows that for a very broad family of random
symmetric integer matrices with independent upper-triangular entries satisfying an anti-concentration
condition (see \Cref{defi: epsilon balanced}), the Sylow $p$-subgroup of the cokernel has a limiting
distribution as $n\to\infty$ which does not depend on the fine details of the entry distribution.
To obtain such convergence, Wood introduced a new method for studying discrete random matrices via
surjection statistics, now often called the surjection moment method. Rather than attempting to
analyze $\Cok(Q_n)$ directly, the method computes, for each fixed finite abelian group $G$, the
expected number of surjections $\mathbb E[\#\Sur(\Cok(Q_n),G)]$. Wood proved that these surjection
moments converge to explicit limits, and that the resulting collection of limits uniquely determines
the limiting distribution of the Sylow $p$-subgroups of the cokernel. Moreover, she also established joint convergence for the $P$-primary part of the cokernel for any
fixed finite set of primes $P$. Here the $P$-primary part means the product of the Sylow
$p$-subgroups over $p\in P$, and we write $\Cok(Q_n)_P$ for the corresponding $P$-primary subgroup of
$\Cok(Q_n)$. This yields a universal limiting law for $\Cok(Q_n)_P$, compatible with the
prime-by-prime limits. As an application, she deduced the limiting distribution of the Sylow $p$-subgroups of sandpile
groups of Erd\H{o}s-R\'enyi random graphs. The resulting limiting formulas confirm conjectures and
predictions from earlier work, notably that of Clancy-Leake-Payne \cite{clancy2015note} and of
Clancy-Kaplan-Leake-Payne-Wood \cite{clancy2015cohen}.

Since then, the surjection moment method has become a standard tool in this area and has led to many
extensions and applications. For random integral matrices with no symmetry constraint, Nguyen-Wood \cite{nguyen2022random} proved universality for random
integral matrices at the level of both surjectivity and the full cokernel distribution, obtaining
precise limiting formulas (in particular, agreeing with Cohen-Lenstra type predictions) under sparse entry distributions. In the symmetric setting, Hodges \cite{hodges2024distribution} strengthens Wood's results by determining the limiting distribution of sandpile
groups together with their canonical pairings. His proof requires upgrading the moment
method to moments of groups equipped with pairing structure. For the alternating case, universality results for cokernels have been proven by Nguyen-Wood \cite{nguyen2025local}. See also Lee \cite{lee2023universality}, Nguyen-Van Peski        \cite{nguyen2024universality}, Huang-Nguyen-Van Peski \cite{huang2025cohen}, Cheong-Yu \cite{cheong2023distribution}, and the current author \cite{shen2025universality} for further related developments.

However, as Wood notes in \cite[Section 1.6]{wood2017distribution}, the moment method
does not provide quantitative information on the rate of convergence, leaving open the problem of
obtaining effective error bounds. The goal of the present paper is to give a new proof of universality for
cokernels of random integral matrices with symmetries, including both the symmetric and alternating cases, via a
completely different approach. It is worth highlighting that our method yields exponentially small error bounds in
most settings, except for a certain special case (the symmetric case at $p=2$), where we still obtain a
stretched-exponential speed of convergence with stretching exponent $1/2$, thereby answering Wood's question by providing explicit convergence
rates.

\begin{defi}\label{defi: epsilon balanced}
We say a random integer $\xi$ is $\epsilon$\emph{-balanced} if for every prime $p$ and every integer $r$, we have
$$\mathbf{P}(\xi\equiv r\mod p)\le 1-\epsilon.$$
\end{defi}

For example, if the integer $\xi$ takes the value $1$ with probability $\epsilon$ and the value $0$ with probability $1-\epsilon$,
then $\xi$ is $\epsilon$-balanced. For all integers $n\ge 1$, let
$$\Sym_n(\Z):=\{M_n\in\Mat_n(\Z):M^T=M\}$$ 
denote the set of symmetric integer matrices of size $n$, and
$$\Alt_n(\Z):=\{A_n\in\Mat_n(\Z):A^T=-A\}$$ 
denote the set of alternating integer matrices of size $n$. Before stating our result in the symmetric case, we briefly recall the canonical pairing on the
cokernel. Let $M_n\in\Sym_n(\Z)$ be nonsingular. Then $M_n$ induces a perfect symmetric pairing $\langle\cdot,\cdot\rangle$ on the finite abelian group $\Cok(M_n)=\Z^n/M_n\Z^n$ as follows: 
\begin{align}
\begin{split}
\langle\cdot,\cdot\rangle: \Cok(M_n)\times \Cok(M_n)&\rightarrow\Q/\Z\\
(x,y)&\mapsto X^TM_n^{-1}Y\quad\mod \Z.
\end{split}
\end{align}
Here, $X,Y\in\Z^n$ are the lifts of $x,y\in\Z^n/M_n\Z^n=\Cok(M_n)$, respectively. For a finite set of primes $P$, we
say that a finite abelian group is a \emph{$P$-group} if its order is a product of powers of primes in
$P$. We also write $\Cok(M_n)_P$ for the $P$-primary
subgroup of $\Cok(M_n)$, and we denote
by $\langle\cdot,\cdot\rangle_P$ the restriction of $\langle\cdot,\cdot\rangle$ to $\Cok(M_n)_P\times
\Cok(M_n)_P$. Thus $(\Cok(M_n)_P,\langle\cdot,\cdot\rangle_P)$ is a finite abelian $P$-primary group
equipped with a perfect symmetric pairing. The readers may refer to \Cref{sec: cok of sym matrix} for more details.

\begin{thm}\label{thm: exponential convergence_sym}
(Symmetric case)
Fix $\epsilon>0$, a finite set of primes $P=\{p_1,\ldots,p_\ell\}$, and a finite abelian group $G$
equipped with a perfect symmetric pairing
\[
\langle\cdot,\cdot\rangle_G: G\times G\to \Q/\Z,
\]
and assume that $P$ contains all primes dividing $\#G$.
Then there exists a constant $K>0$ depending on $G,P,\epsilon$, and a constant $c>0$ depending on $G,P,\epsilon$, such that the following
holds for every $n\ge 1$ and every random matrix $M_n\in \Sym_n(\Z)$ whose entries on and above the
diagonal are independent $\epsilon$-balanced random integers:
\begin{multline}
\left|\mathbf{P}\Bigl(M_n\ \textup{is nonsingular and}\ 
(\Cok(M_n)_P,\langle\cdot,\cdot\rangle_P)\simeq (G,\langle\cdot,\cdot\rangle_G)\Bigr)
-\frac{\prod_{i=1}^\ell\prod_{k\ge 1}(1-p_i^{1-2k})}{|G|\cdot|\Aut(G,\langle\cdot,\cdot\rangle_G)|}\right|\\
\le K\exp(-c n^{1/2}),
\end{multline}
where $\Aut(G,\langle\cdot,\cdot\rangle_G)$ denotes the group of automorphisms of $G$ preserving the
pairing.

Moreover, if $2\notin P$, then there exists a constant $K'>0$ depending on $G,P,\epsilon$, and a constant $c'>0$ depending on $G,P$, such
that under the same assumptions,
\begin{multline}
\left|\mathbf{P}\Bigl(M_n\ \textup{is nonsingular and}\ 
(\Cok(M_n)_P,\langle\cdot,\cdot\rangle_P)\simeq (G,\langle\cdot,\cdot\rangle_G)\Bigr)
-\frac{\prod_{i=1}^\ell\prod_{k\ge 1}(1-p_i^{1-2k})}{|G|\cdot|\Aut(G,\langle\cdot,\cdot\rangle_G)|}\right|\\
\le K'\exp(-\epsilon c'n).
\end{multline}
\end{thm}

Before stating our result for alternating matrices, we introduce some notation. Let $\mathcal S$
denote the class of finite abelian groups of the form $H\oplus H$. For a finite set of primes $P$, we write $\mathcal{S}_P$ for the set of $P$-groups in $\mathcal S$. For a prime $p$, we similarly write
$\mathcal{S}_p$ for the set of $p$-groups in $\mathcal S$.

A finite abelian $P$-group admits a perfect alternating pairing to $\Q/\Z$ if and only if it lies in
$S_P$. Such a group $G$ has a unique perfect alternating pairing up to isomorphism; we write
$\Sp(G)$ for the group of automorphisms of $G$ preserving this pairing. If $A_n\in \Alt_n(\Z)$, then $A_n$
has an even rank, and the torsion subgroup of the cokernel, which we denote by $\Cok_{\tors}(A_n)$, lies in $\mathcal S$. See Bhargava-Kane-Lenstra-Poonen-Rains \cite[Sections 3.4-3.5]{bhargava2015modeling} for further
discussions.

\begin{thm}\label{thm: exponential convergence_alt}
(Alternating case)
Fix $\epsilon>0$, a finite set of primes $P=\{p_1,\ldots,p_\ell\}$, and a finite abelian group $G$.
Assume that $P$ contains all primes dividing $\#G$. 

If $G\in\mathcal S_P$, then there exists a constant $K>0$ depending on $G,P,\epsilon$, and a constant $c>0$ depending on $G,P$, such that
the following holds for every $n\ge 1$ and random alternating matrices
$A_{2n}\in\Alt_{2n}(\Z)$ and $A_{2n+1}\in\Alt_{2n+1}(\Z)$ whose entries above the diagonal are
independent $\epsilon$-balanced random integers:
\[
\left|\mathbf{P}\bigl(\Cok(A_{2n})_P\simeq G\bigr)
-\frac{|G|}{|\Sp(G)|}\prod_{i=1}^\ell\prod_{k\ge 1}(1-p_i^{1-2k})\right|
\le K\exp(-\epsilon c n),
\]
and
\[
\left|\mathbf{P}\bigl(\Cok_{\tors}(A_{2n+1})_P\simeq G,\ \corank(A_{2n+1})=1\bigr)
-\frac{1}{|\Sp(G)|}\prod_{i=1}^\ell\prod_{k\ge 1}(1-p_i^{-1-2k})\right|
\le K\exp(-\epsilon c n).
\]
Otherwise, if $G\notin\mathcal S_P$, then
\[
\mathbf{P}\bigl(\Cok(A_{2n})_P\simeq G\bigr)=
\mathbf{P}\bigl(\Cok_{\tors}(A_{2n+1})_P\simeq G\bigr)=0.
\]
\end{thm}

\Cref{thm: exponential convergence_sym} strengthens Hodges \cite[Corollary 10.5]{hodges2024distribution}, which is a multi-prime generalization of
\cite[Theorem 1.1]{hodges2024distribution}; \Cref{thm: exponential convergence_alt} is a stronger version of Nguyen-Wood \cite[Theorem 1.13]{nguyen2025local}. As far as we know, our results provide the first quantitative convergence rates for cokernel universality in random integral matrix ensembles with symmetry. Moreover, our method also yields quantitative convergence for the joint distribution of the
outermost corners: for any fixed number of successive principal minors near full size, we obtain
effective error bounds for their joint cokernel laws; see \Cref{prop: joint distribution_sym} for
the symmetric case, and \Cref{prop: joint distribution_alt} for the alternating case. These results may be
viewed as non-archimedean analogues of the GOE and aGUE corners processes, respectively, which
have attracted substantial attention in integrable probability. We are not aware of any way to access these corner joint distributions via
the surjection moment method.

\subsection{A brief summary of our method}

We have so far discussed the progress on studying cokernels of random matrices via the surjection moment method.
However, there are also other approaches for proving the relevant universality phenomena. A particularly notable example is the work of Maples \cite{maples2013cokernels},
which proves universality for cokernels of random integral matrices with no symmetry constraint (over $\mathbb Z_p$ and, more generally, over $\mathbb Z/a\mathbb Z$) under the same $\epsilon$-balanced hypotheses on the entries, and moreover obtains exponential convergence to the Cohen–Lenstra distribution.

The main ingredient in Maples' proof is a column-exposure process that can be viewed as a random walk on
the space of finite $R$-modules (with $R=\mathbb Z_p$ or $\mathbb Z/a\mathbb Z$). As a heuristic of the main idea, one reveals the
columns $X_1,\dots,X_n$ sequentially and considers the evolving quotient modules
\[
L_k \;=\; R^n/\langle X_{k+1},\dots,X_n\rangle ,
\]
so that $(L_k)_{k=0}^n$ forms a random process whose increments are induced by adding one random column at a time. The key step is to compare this chain of quotient modules with the corresponding Haar model $(L_k')_{k=0}^n$, in which the columns
$X_1',\dots,X_n'$ are independent and Haar-distributed. Maples proves that, at each step, the transition kernel of
$(L_k)$ is exponentially close (in total variation distance) to the Haar transition kernel, uniformly over the
current state, except on a negligible exceptional set. Consequently, the law of the terminal module $L_0=\Cok(Q_n)$
is exponentially close to the law of the Haar cokernel $L_0'=\Cok(Q_n')$. Combining this comparison with the explicit Haar computation of
Friedman-Washington \cite{friedman1989distribution} yields Cohen-Lenstra type universality and exponential convergence.

Our approach for symmetric and alternating random matrices is motivated by an exposure-process
viewpoint similar to that of the above. Our first step is to pass to quotient rings:
as in the surjection moment method, we work over $\Z/a\Z$ for a sufficiently large integer $a$ in order to
encode the desired cokernel information simultaneously at a finite set of primes. After this
reduction, we exploit an exposure process to compare the resulting dynamics over the quotient ring to the
corresponding uniform model. We now focus on the case of random symmetric matrices, and the
alternating case is analogous. Let $M_n\in \Sym_n(\Z)$ be as in
\Cref{thm: exponential convergence_sym}. For each $1\le t\le n$, let $M_t$ denote the $t\times t$ upper-left corner of $M_n$, and write
$M_t/a$ for its reduction modulo $a$, viewed as an element of $\Sym_t(\Z/a\Z)$.
This yields the exposure process
\[
M_1/a \subset M_2/a \subset \ldots \subset M_n/a,
\]
obtained by iteratively revealing one new random row and column each time. Under reasonable non-sparsity assumptions (discussed in \Cref{sec: nonsparse}), the one-step
transition kernel of this $\epsilon$-balanced exposure process over $\Z/a\Z$ is close to that of the
corresponding uniform model on symmetric matrices over $\Z/a\Z$. By the Chinese remainder theorem, the uniform
model modulo $a=\prod_{i=1}^lp_i^{e_{p_i}}$ decomposes into independent uniform models modulo each
$p_i^{e_{p_i}}$. Furthermore, the reduction of a Haar-distributed matrix over $\Z_{p_i}$ modulo $p_i^{e_{p_i}}$
is uniform on $\Z/p_i^{e_{p_i}}\Z$. Consequently, we may couple the $\epsilon$-balanced exposure process with the Haar
models. The limiting distribution of the Haar case is already explicitly understood by previous works:
in the symmetric case it is resolved in Clancy-Kaplan-Leake-Payne-Wood \cite[Theorem 2]{clancy2015cohen}, and in the alternating case
in Bhargava-Kane-Lenstra-Poonen-Rains \cite[Section 3]{bhargava2015modeling}. Combining this exposure process comparison with those explicit Haar results yields the
desired universality statements and error bounds.

Maples also adopted the exposure-process point of view to study the rank distribution of symmetric random matrices over the finite field $\F_p$ in the unpublished manuscript \cite{Maples_symma_2013}. The key
analytic input there is a collection of quadratic Littlewood-Offord estimates (in the spirit of
the inverse theory initiated earlier by Costello-Tao-Vu \cite{costello2006random}), which are used
to control the probability that a newly revealed row lies in an atypical subspace. However, the
quantitative bounds obtained in Maples' manuscript are only of polynomial order, which is not
strong enough for our purposes and does not by itself yield
\Cref{thm: exponential convergence_sym} and \Cref{thm: exponential convergence_alt}.

Fortunately, a crucial improvement is due to Ferber-Jain-Sah-Sawhney \cite{ferber2023random}:
working in the same finite-field setting $\F_p$ and within the same exposure-process framework, they
strengthen these quadratic Littlewood-Offord estimates so as to upgrade Maples' polynomial error
terms to exponential decay. Building on their strategy, we implement an adaptation in the quotient
ring setting $\Z/a\Z$. This yields typically exponentially small error bounds; in the special
symmetric case at $p=2$, the best bound we obtain is stretched-exponential with stretching exponent
$1/2$. As a byproduct, this approach also extends to quantitative controls of the joint distribution
of corners given in \Cref{prop: joint distribution_sym} and \Cref{prop: joint distribution_alt}.

\begin{rmk}
One might wonder why a stretched-exponential error bound appears in the symmetric case at $p=2$.
This comes from a technical obstruction in our Fourier-analytic estimates: for symmetric matrices, the relevant
character sums naturally involve quadratic terms, and when $p=2$ these quadratic contributions become genuinely
harder to deal with. In fact, the Fourier-estimate input we build on (in the spirit of Ferber-Jain-Sah-Sawhney) is only available for
odd primes. To overcome this difficulty at $p=2$, we impose a revised non-sparsity assumption (see \Cref{defi: revised sym p=2}), under which our method still yields quantitative convergence but only at a stretched-exponential
rate.

By contrast, in the alternating case, the diagonal entries are identically zero, so the corresponding Fourier
analysis does not produce quadratic terms that need special treatment. As a result, our error bounds in the
alternating setting are always exponential.
\end{rmk}

\subsection{Further questions}

Our new approach also raises a number of further questions. We expect that the error bounds in \Cref{thm: exponential convergence_sym} and \Cref{thm: exponential convergence_alt} are not sharp and can be improved. In particular, it is plausible that one could
upgrade the stretched-exponential bound in the symmetric case at $p=2$ to exponential convergence, but we do
not currently see how to achieve this within our framework (in fact, even over the finite field $\F_2$, we do not know how to prove an exponentially small error
term for the distribution of the corank of a $\epsilon$-balanced random symmetric matrix.). Moreover, the constants in \Cref{thm: exponential convergence_sym} and \Cref{thm: exponential convergence_alt} may well be made absolute (that is, depending only on the number of primes in the set $P$), as is known for unconstrained random matrices (i.e., with no symmetry constraint) in the work of Maples.

Another direction is to understand regimes in which the prime $p$ grows with $n$, or the matrix entries are sparse
relative to $n$. The large $p$ case has been treated in the finite field setting of Ferber-Jain-Sah-Sawhney \cite[Theorem 1.2]{ferber2023random}, where they allow
$p$ to grow stretched-exponentially with stretching exponent $1/4$ when the matrix size goes to infinity. Sparse-entry models have also been studied by Nguyen-Wood \cite{nguyen2022random} for random integral matrices without symmetry constraint. It may be possible to adapt aspects of their approach to our symmetric and alternating
integral settings.

The results in this paper concern local universality, in the sense that we always work with a fixed finite set of
primes $P$ and study the $P$-primary part of the cokernel. However, universality questions at the global level are
also of considerable interest. For example, Lorenzini \cite{lorenzini2008smith} asked how often sandpile groups of graphs are cyclic, and this
was resolved by Nguyen and Wood in \cite{nguyen2025local} by proving the existence of a limiting probability as
$n\to\infty$, though without a quantitative rate of convergence. With our new approach, it seems plausible that one
could obtain nontrivial error bounds for such global questions, a direction for which we are not aware of any prior
progress.

A closely related symmetry class not treated in this paper is the Hermitian setting over the ring of
integers $\mathcal O_K$ of a quadratic extension $K=\Q(\sqrt d)$, where one considers matrices
$H_n\in \Mat_n(\mathcal{O}_K)$ satisfying $H^{\ast}_n=H_n$ with respect to the nontrivial Galois conjugation
on $K$. In the local field setting, Lee \cite[Theorem 1.6]{lee2023universality} established
universality for cokernels of random $p$-adic Hermitian matrices over the ring of integers of a
quadratic extension of $\Q_p$, together with an explicit limiting distribution. Lee's argument is
based on Wood's surjection moment method and, accordingly, does not provide quantitative rates of
convergence; moreover, it is local at a single prime $p$ and does not address the
multi-prime setting relevant to our results. We nevertheless believe that the exposure-process viewpoint developed in our paper should also apply in the
Hermitian context (after localizing at primes and working over suitable quotient rings), potentially
yielding exponentially small error term and allowing one to treat multi-prime statistics
in a unified way.

A closely related line of work studies the probability that a discrete random matrix is nonsingular (equivalently,
that its cokernel is finite). In the symmetric
Bernoulli model, a recent breakthrough of Campos-Jenssen-Michelen-Sahasrabudhe \cite{campos2025singularity}
establishes that the singularity probability of an $n\times n$ random symmetric $\{\pm1\}$-matrix $M_n$ is exponentially
small, i.e.\ $\mathbf P(\det M_n=0)\le e^{-cn}$ for some absolute $c>0$. Their argument can also be adapted to the alternating setting, leading to an exponential bound on
the singularity probability of a random $2n\times 2n$ alternating Bernoulli matrix. In our $\epsilon$-balanced setting, Costello-Tao-Vu \cite[Theorem 6.5]{costello2006random}
already states that the singularity probability goes to $0$ as $n\to\infty$. It is natural to expect exponentially fast convergence in this more general
setting as well.

\subsection{Outline of the paper}

In \Cref{sec: Preliminaries}, we collect background on finite abelian groups, together with several
probabilistic and Fourier-analytic tools used throughout the paper. In \Cref{sec: cok of sym matrix} and \Cref{sec: cok of alt matrix}, we clarify the
structure of cokernels of symmetric and alternating matrices, with particular emphasis on the
relationship between the isomorphism class of the cokernel and congruence equivalence of matrices. In \Cref{sec: nonsparse}, we show that, with high probability, the random symmetric and random alternating matrices
considered in this paper satisfy non-sparsity assumptions. These conditions provide the technical foundation
needed to carry out the Fourier-analytic estimates for the transition probabilities in the exposure process. In \Cref{sec: Proof of the symmetric case} and \Cref{sec: Proof of the alternating case}, we give the proof of \Cref{thm: exponential convergence_sym} and \Cref{thm: exponential convergence_alt} respectively. Finally, in \Cref{sec: index of notations}, we include a table that lists the notations appearing throughout the paper, together with brief explanations.

\subsection{Notations}

We write $[n]$ for the set $\{1,2,\ldots,n\}$. For an index set $I\subset [n]$, we write $I^c$ for the complement of $I$ in $[n]$. We denote the order of groups and sets using either absolute value sign $|\cdot|$ or $\#$.

\textbf{Probability.} We write $\mathbf{P}(\cdot)$ for probability, $\E(\cdot)$ for expectations, $\mathcal{L}(\cdot)$ for the law of a random variable. Moreover, we use $\mathbf{P}(\cdot\mid\cdot),\E(\cdot\mid\cdot),\mathcal{L}(\cdot\mid\cdot)$ for conditional probability, conditional expectation and conditional law. 

For two probability measures $\mu$ and $\nu$ on a countable set $E$, the $L^q$-distance between $\mu$ and $\nu$ is defined by 
$$D_{L^q}(\mu,\nu)=\left(\sum_{x\in E}(|\mu(x)-\nu(x)|)^q\right)^{1/q}.$$
We also use the abbreviation $D_{L^q}(\mathcal{X},\mathcal{Y}):=D_{L^q}(\mathcal{L}(\mathcal{X}),\mathcal{L}(\mathcal{Y}))$ for the $L^q$-distance of two random variables $\mathcal{X}$ and $\mathcal{Y}$. In particular, we write $\mathcal{X}\stackrel{d}{=}\mathcal{Y}$ when two random variables $\mathcal{X},\mathcal{Y}$ have the same distribution, and we use the notation $\mathcal{X}_n\stackrel{d}{\rightarrow} \mathcal{X}$ for weak convergence.

\textbf{Analysis.} We use $\exp(\cdot)$ for the exponential function. For a set of parameters $\mathbf{S}$, we write $f(n,\ldots)=O_\mathbf{S}(g(n,\ldots))$ if $|f|\le K|g|$ for some large constant $K=K(\mathbf{S})>0$ depending on $\mathbf{S}$, and we write $f(n,\ldots)=\Omega_\mathbf{S}(g(n,\ldots))$ if $f\ge c|g|$ for some small constant $c=c(\mathbf{S})>0$ depending on $\mathbf{S}$.

\textbf{Linear algebra.} 
For a vector $\bm{\xi}=(\xi_1,\ldots,\xi_n)$, we denote by $\supp(\bm{\xi}):=\{i\in[n]:\xi_i\ne 0\}$ the \emph{support}, and $\wt(\bm{\xi}):=\#\supp(\bm{\xi})$ the \emph{weight} of $\bm{\xi}$. For a given index set $J\subset[n]$, when it appears as a subscript (for example, the symbol $\bm{\xi}_J$), we usually refer to a vector indexed from $J$. In particular, when $\bm{\xi}=(\xi_1,\ldots,\xi_n)$ is a vector with index set $[n]$, the symbol $\bm{\xi}_J$ often refers to the subvector of $\bm{\xi}$ restricted to its components indexed from $J$. 

For $I,J\subset[n]$, the matrix $B_{I\times J}$ is the submatrix of the rows and columns indexed from $I$ and $J$ respectively. 
Unless specifically emphasized, all vectors without the transpose superscript are column vectors, despite being written as rows. 

For a matrix $B$ with entries in $\Z$ and an integer $a\ge 2$, we write $B/a$ for the matrix with coefficients in $\Z/a\Z$ obtained from $B$ by reduction mod $a$. Furthermore, if the entries of $B$ are in $\Z/a\Z$ and $a'\mid a$, we write $B/a'$ for the canonical reduction to $\Z/a'\Z$. We also adapt analogous notation $\bm{\xi}/a$ for a vector $\bm{\xi}$. Given two vectors $\bm{\xi},\bm{\eta}$ of the same length (for example, they are indexed from the same set $J\subset[n]$), we write $\bm{\xi}\cdot\bm{\eta}$ as their dot product.

For two $n\times n$ matrices $B_n,B_n'$ with entries in the $p$-adic integers $\Z_p$, we write $B_n\congsim B_n'$ if they are congruent, i.e., there exists $U\in\GL_n(\Z_p)$ such that $UB_nU^T=B_n'$. Moreover, when we say "congruent", "congruent transformation" or "congruence equivalence" for matrices in $\Q_p$, we always refer to a congruence transform of a matrix in $\GL_n(\Z_p)$ instead of $\GL_n(\Q_p)$.

%% file: Preliminaries.tex
\section{Preliminaries}\label{sec: Preliminaries}

\subsection{Preparations from coding theory}

\begin{defi}
For a prime number $p$, denote by 
$$H_p(x):=x\log_p(p-1)-x\log_p x-(1-x)\log_p(1-x),\quad 0<x<1$$
the $p$-\emph{ary entropy function}. 
\end{defi}

\begin{defi}
For all integers $0\le t\le n$ and prime numbers $p$,  denote by $\Vol_p(n,t)$ the order of the Hamming ball in $\F_p^n$ of radius $t$ (i.e., the elements in $\F_p^n$ that have at most $t$ nonzero entries). Clearly we have
$$\Vol_p(n,t)=\sum_{i=0}^t\binom{n}{i}(p-1)^i.$$
\end{defi}

The following proposition is quoted from \cite[Proposition 3.3.1]{guruswami2012essential}, which shows that the $p$-ary entropy function describes the asymptotics of the Hamming ball.

\begin{prop}\label{prop: asymptotic of Hamming ball}
Let $p$ be a prime, and $0<x<1-1/p$ be a real number. When $n$ tends to infinity, we have
$$p^{(H_p(x)-o(1))n}\le \Vol_p(n,xn)\le p^{H_p(x)n}.$$
\end{prop}




\subsection{Partitions and finite abelian groups}

\begin{defi}
We denote by $\Y$ the set of \emph{partitions} $\lambda=(\lambda_1,\lambda_2,\ldots)$, which are (finite or infinite) sequences of nonnegative integers $\lambda_1\ge\lambda_2\ge\cdots$ that are eventually zero. We do not distinguish between two such sequences that differ only by a string of zeros at the end. The integers $\l_i>0$ are called the \emph{parts} of $\l$. We set $|\l| := \sum_{i\ge 1}\l_i,n(\lambda):=\sum_{i\ge 1} (i-1)\lambda_i$, and $m_k(\l):= \#\{i\mid \l_i = k\}$. We also write $\Len(\lambda)=\#\{i\mid \l_i >0\}$ as the \emph{length} of $\lambda$, which is the number of nonzero parts.
\end{defi}

In the coming sections, we usually write $P$ for a finite set of prime numbers, and $p$ a prime number. We also denote by $\Q_p$ the $p$-adic field, $\val(\cdot):\Q_p^\times\rightarrow\Z$ the valuation function, $\Z_p:=\{\alpha\in\Q_p^\times:\val(\alpha)\ge 0\}\cup\{0\}$ the $p$-adic integers, and $\F_p:=\Z/p\Z$ the finite field of size $p$. For a finite abelian group $G$ and a prime $p$, we write $G_p$ for the Sylow $p$-subgroup of $G$, so that
$$G=\bigoplus_{p\text{ prime}}G_p.$$
We will also write $G_P:=\bigoplus_{p\in P} G_p$ the $P$-primary subgroup of $G$. As a finite abelian $p$-group, $G_p$ is isomorphic to $\bigoplus_{i\ge 1}(\Z/p^{\l_i}\Z)$ for some $\l=(\l_1,\l_2\ldots)\in\Y$. In this case, we say $G$ has type $\l$ at $p$, and we denote by $\Dep_p(G):=\l_1$ the $p$-\emph{depth} of $G$, i.e., the smallest integer such that $p^{\Dep_p(G)}G_p=0$. We also denote by $\Len_p(G):=\Len(\l)=n$ the $p$-\emph{length} of $G$, i.e., the number of cyclic factors in $G_p$. In particular, the $p$-torsion subgroup $G_p$ is trivial if and only if $\Dep_p(G)=\Len_p(G)=0$.

\subsection{Auxiliary tools}

In this subsection, we introduce some useful lemmas that will be useful in later sections. We recall the following basic linear algebra lemma about full-rank principal minors. The symmetric part is essentially contained in \cite[Lemma 3.1]{ferber2023random}, and the alternating part follows analogously.

\begin{lemma}\label{lem: full rank principal minor}
Let $p$ be a prime number, and $n\ge 1$. 
\begin{enumerate}
\item (Symmetric case) Suppose $M_n\in\Sym_n(\F_p)$ has rank $r$. Then there is a $r\times r$ principal minor with determinant in $\F_p^\times$.
\item (Alternating case) Suppose $A_n\in\Alt_n(\F_p)$ has rank $2r$. Then there is a $2r\times 2r$ principal minor with determinant in $\F_p^\times$.
\end{enumerate}
\end{lemma}

\begin{rmk}\label{rem:principal-minor-lifts}
\Cref{lem: full rank principal minor} has natural analogues for matrices with entries in
$\Z_p$ or in the quotient ring $\Z/p^{e_p}\Z$ (for any $e_p\ge 1$). Indeed, given
$M_n\in \Sym_n(\Z_p)$ (or $M_n\in \Sym_n(\Z/p^{e_p}\Z)$), one may reduce modulo $p$ to obtain a matrix in
$\Sym_n(\F_p)$ and apply the lemma there to find a full-rank principal minor. The corresponding
principal minor of the original matrix then has determinant in $\Z_p^\times$ (respectively in
$(\Z/p^{e_p}\Z)^\times$). The same applies in the alternating case.
In later arguments we will routinely use this lifted version of the lemma for matrices over $\Z_p$
or $\Z/p^{e_p}\Z$.
\end{rmk}

\begin{defi}
We say that a random variable $\xi\in\Z/a\Z$ is \emph{not $\epsilon$-concentrated} if $\mathbf{P}(\xi\equiv r\mod p)\le1-\epsilon$ for all $p\mid a$ and $r\in \F_p$. We also say that a random integer $\xi\in\Z$ is \emph{not $\epsilon$-concentrated at} $a$ if $\xi+a\Z\in\Z/a\Z$ is not $\epsilon$-concentrated.
\end{defi}

\begin{lemma}\label{lem: expectation of not concentrated}
Let $z$ be random and not $\epsilon$-concentrated in $\Z/a\Z$. Then we have
$$\left|\E\left(\exp\left(\frac{2\pi\sqrt{-1}}{a}z\right)\right)\right|\le\exp(-\frac{\epsilon}{a^2}).$$
\end{lemma}

\begin{proof}
See \cite[Lemma 4.2]{wood2017distribution}.
\end{proof}

%% file: The_symmetric_pairing_on_the_cokernel_of_a_nonsingular_symmetric_matrix.tex
\section{The symmetric pairing on the cokernel of a nonsingular symmetric matrix}\label{sec: cok of sym matrix}

The goal of this section is to set up the structural input needed
for the exposure-process argument in \Cref{sec: Proof of the symmetric case}. We first formalize the notion of a \emph{paired}
finite abelian group (see \Cref{defi: symmetric abelian group}) and recall that the cokernel of a nonsingular symmetric integral (or $p$-adic)
matrix carries a canonical perfect symmetric pairing. Working prime-by-prime, in \Cref{thm: add invertible to same cokernel with pairing}, we relate the
resulting isomorphism class $(\Cok(M),\langle\cdot,\cdot\rangle)$ to congruence equivalence of symmetric
matrices (up to adjoining invertible blocks), thereby reducing questions about cokernels with pairing
to statements about congruence classes. Finally, in \Cref{defi: quasi pairing of Cok}, we pass to reductions modulo $a$ and introduce the
notion of a \emph{cokernel with quasi-pairing} for symmetric matrices over $\mathbb{Z}/p^{e_p}\mathbb{Z}$.
With a suitable choice of modulus $a$, this provides a finite state space that records exactly the
$P$-primary cokernel together with its pairing, and it is the framework in which we later define and
estimate the transition probabilities of the symmetric exposure process.

\begin{defi}\label{defi: symmetric abelian group}
Let $G$ be a finite abelian group, equipped with a bilinear symmetric pairing $\langle\cdot,\cdot\rangle_G:G\times G\rightarrow\Q/\Z$. We say that $(G,\langle\cdot,\cdot\rangle_G)$ is a \emph{paired abelian group} (or simply a \emph{paired group}) if the pairing $\langle\cdot,\cdot\rangle_G$ is perfect. Additionally, given a prime number $p$, we say $(G,\langle\cdot,\cdot\rangle_G)$ is a \emph{paired abelian $p$-group} (or simply a \emph{paired $p$-group}) if it is a paired abelian group and a $p$-group. We say that two paired abelian groups $(G,\langle\cdot,\cdot\rangle_G)$ and $(G',\langle\cdot,\cdot\rangle_{G'})$ are isomorphic and write $(G,\langle\cdot,\cdot\rangle_G)\simeq(G',\langle\cdot,\cdot\rangle_{G'})$, if there exists a group isomorphism that preserves the pairing.
\end{defi}

The group $G$ together with its pairing $\langle\cdot,\cdot\rangle_G$ in \Cref{thm: exponential convergence_sym} form a finite paired abelian group. From now on, all paired abelian groups are assumed to be finite. As discussed in \Cref{sec: intro}, paired abelian groups naturally emerge when considering the cokernel of a nonsingular symmetric integral matrix $M_n\in\Sym_n(\Z)$. Now, for a given prime $p$, we focus on the Sylow $p$-subgroup of this cokernel, where we regard the matrix entries as elements in the $p$-adic integers $\Z_p$. When $M_n$ is nonsingular, the $p$-torsion part of $\Cok(M_n)$, denoted by $\Cok(M_n)_p:=\Z_p^n/M_n\Z_p^n$, gives a paired abelian $p$-group. Indeed, its symmetric pairing, denoted by $\langle\cdot,\cdot\rangle_p$, is given by
\begin{align}
\begin{split}
\langle\cdot,\cdot\rangle_p: \Cok(M_n)_p\times \Cok(M_n)_p&\rightarrow\Q_p/\Z_p\\
(x,y)&\mapsto X^TM_n^{-1}Y\quad\mod \Z_p.
\end{split}
\end{align}
Here, $X,Y\in\Z_p^n$ are the lifts of $x,y\in\Cok(M_n)_p$, respectively. 

Whenever it is clear from the text that we are in the localized-at-$p$ setting where $M_n$ is regarded as a matrix with entries in $\Z_p$, we will sometimes drop the subscript $p$ and simply write $(\Cok(M_n),\langle\cdot,\cdot\rangle)$ for $(\Cok(M_n)_p,\langle\cdot,\cdot\rangle_p)$. For the global setting, we will be interested in the localizations of $M_n\in\Sym_n(\Z)$ at a finite set of primes, denoted by $P$. In this case, we write $(\Cok(M_n)_P,\langle\cdot,\cdot\rangle_P)$ for the restriction of $(\Cok(M_n),\langle\cdot,\cdot\rangle)$ to the primes in $P$, or equivalently, the collection of abelian groups with pairings $((\Cok(M_n)_p,\langle\cdot,\cdot\rangle_p))_{p\in P}$.

Following the above definitions, it is clear that when $M_{n_1}\in\Sym_{n_1}(\Z_p)$ is nonsingular, and $M_{n-n_1}^c\in\Sym_{n-n_1}(\Z_p)\cap\GL_{n-n_1}(\Z_p)$, we have 
$$(\Cok(M_{n_1}),\langle\cdot,\cdot\rangle)\congsim\left(\Cok\begin{pmatrix}M_{n-n_1}^c & 0 \\ 0 & M_{n_1}\end{pmatrix},\langle\cdot,\cdot\rangle\right).$$

\begin{prop}\label{prop: same pairing matrix of same size}
Let $n\ge 1$, and $M_n,M_n'\in\Sym_n(\Z_p)$ be nonsingular. Then, $(\Cok(M_n),\langle\cdot,\cdot\rangle)\simeq(\Cok(M_n'),\langle\cdot,\cdot\rangle)$ if and only if there exists a nonsingular $M_n''\in\Sym_n(\Z_p)$, such that $M_n''\congsim M_n'$, and $M_n^{-1}-(M_n'')^{-1}\in\Sym_n(\Z_p)$.
\end{prop}

\begin{proof}
On the one hand, suppose that such a matrix $M_n''$ exists. Then we have
$$(\Cok(M_n),\langle\cdot,\cdot\rangle)\simeq(\Cok(M_n''),\langle\cdot,\cdot\rangle)\simeq(\Cok(M_n'),\langle\cdot,\cdot\rangle).$$
On the other hand, suppose $(\Cok(M_n),\langle\cdot,\cdot\rangle)\simeq(\Cok(M_n'),\langle\cdot,\cdot\rangle)$. By the definition of isomorphism between paired abelian $p$-groups given in \Cref{defi: symmetric abelian group}, there exists $U\in\GL_n(\Z_p)$ such that $M_n'':=UM_n'U^T$ and $M_n$ generate the same pairing that takes values in $\Q_p/\Z_p$, i.e., $M_n^{-1}-(M_n'')^{-1}\in\Sym_n(\Z_p)$. This completes the proof.
\end{proof}

\Cref{prop: same pairing matrix of same size} motivates us to ask for the simplest form of a symmetric matrix in $\Z_p$ under the congruence transformations of $\GL_n(\Z_p)$, which is given in the following proposition. 

\begin{thm}\label{prop: simplest matrix_sym}
Let $M_n\in\Sym_n(\Z_p)$ be nonsingular. 
\begin{enumerate}
\item If $p$ is odd, take a fixed non-square $\mathfrak{r}_p\in\Z_p^\times$. Then $M_n$ is congruent to a diagonal matrix of the form
\begin{equation}\label{eq: simplest matrix sym odd}
\diag_{n\times n}(p^{\l_1},\ldots,p^{\l_n})\cdot\diag(\varepsilon_1,\ldots,\varepsilon_n),\end{equation}
where $\l:=(\l_1,\ldots,\l_n)\in\Y$, and $\varepsilon_1,\ldots,\varepsilon_n\in\{1,\mathfrak{r}_p\}$ such that for $1\le i\le n-1$ with $\l_i=\l_{i+1}$, we have $\varepsilon_i=1$. Furthermore, no two distinct matrices in \eqref{eq: simplest matrix sym odd} are congruent under the action of $\GL_n(\Z_p)$.
\item If $p=2$, then $M_n$ is congruent to a block diagonal matrix whose diagonal blocks are one of the following form:
\begin{equation}\label{eq: simplest matrix sym 2}
(2^d),(3\times2^d),(5\times2^d),(7\times2^d),\begin{pmatrix}0 & 2^d\\ 2^d & 0\end{pmatrix},\begin{pmatrix}2^{d+1} & 2^d\\ 2^d & 2^{d+1}\end{pmatrix}.
\end{equation}
Here, $d\ge 0$ is a non-negative integer.
\end{enumerate}
\end{thm}

\begin{proof}
For the $p>2$ part see \cite[Chapter 8.3]{cassels2008rational}, and for the $p=2$ part see \cite[Chapter 8.4]{cassels2008rational}.
\end{proof}

Note that when $p=2$, we do not assert that the simplest form is unique expressed by the components listed in \eqref{eq: simplest matrix sym 2}. 

\begin{prop}\label{prop: locally congruent equivalence at p} 
Let $p$ be a prime number. Let $n\ge 1$, and $M_n\in\Sym_n(\Z_p)$ be nonsingular. When $p$ is odd, let $e_p=\Dep_p(\Cok(M_n))+1$, and when $p=2$, let $e_2=\Dep_2(\Cok(M_n))+3$. Then for all $\Delta M_n\in\Sym_n(p^{e_p}\Z_p)$, we have $M_n\congsim M_n+\Delta M_n$.
\end{prop}

\begin{proof}
Applying a congruence transformation by $\GL_n(\Z_p)$, there is no loss of generality to assume that $M_n$ is already the simplest form given in \Cref{prop: simplest matrix_sym}. We will only give the proof of the case $p=2$, and the odd prime case follows similarly but more straightforwardly. 

We have already assumed that $M_n$ is a diagonal block matrix, whose diagonal blocks are one of the forms in \eqref{eq: simplest matrix sym 2}. Each block has an exponent $d$, and $\Dep_2(\Cok(M_n))$ is the maximum of all these exponents. Turning to $M_n+\Delta M_n$, we can use the diagonal blocks in the corresponding positions as pivots to clear the off-block-diagonal entries, and under such simultaneous row-column eliminations, each entry in the diagonal blocks varies only $2^{e_2}\Z_2$. 
Thus, we have $M_n+\Delta M_n\congsim M_n+\Delta M_n'$ for some $\Delta M_n'\in\Sym_n(2^{e_2}\Z_2)$, which has the same block structure as $M_n$ (i.e., the diagonal blocks occur in the same positions and have the same sizes).

Now, we only need to show that, for any block listed in \eqref{eq: simplest matrix sym 2}, modifying its entries by elements of $2^{e_2}\Z_2$ while maintaining symmetry does not affect its congruent equivalence class. Since every entry contains a factor $2^d$ such that $e_2\ge d+3$, we may assume $d=0,e_2=3$ without loss of generality by factoring it out, thereby reducing the situation to the following cases.
\begin{enumerate}
\item For $b\in\Z_2$, we have congruent equivalence
$$(1)\congsim (1+8b),(3)\congsim (3+8b),(5)\congsim (5+8b),(7)\congsim (7+8b)$$ 
due to the fact that $(\Z_2^\times)^2=1+8\Z_2$.
\item For $\xi,b_1,b_2\in\Z_2$, we claim congruent equivalence $\begin{pmatrix}0 & 1\\ 1 & 0\end{pmatrix}\congsim\begin{pmatrix}8b_1 & 1+8\xi\\ 1+8\xi & 8b_2\end{pmatrix}$. In fact, applying Hensel's lemma, we can find $\eta\in\Z_2$ that satisfies 
$$8b_2\eta^2-2(1+8\xi)\eta+8b_1=0$$
due to the fact that the quadratic term and constant term have valuation $\ge 3$, and the linear term has valuation $1$. Therefore, we have
\begin{align}
\begin{split}
\begin{pmatrix}8b_1 & 1+8\xi\\ 1+8\xi & 8b_2\end{pmatrix}&\congsim\begin{pmatrix}8b_2\eta^2-2(1+8\xi)\eta+8b_1 & 1+8\xi-8b_2\eta\\ 1+8\xi-8b_2\eta & 8b_2\end{pmatrix}\\
&=\begin{pmatrix}0 & 1+8\xi-8b_2\eta\\ 1+8\xi-8b_2\eta & 8b_2\end{pmatrix}\\
&\congsim\begin{pmatrix}0 & 1+8\xi-8b_2\eta\\ 1+8\xi-8b_2\eta & 0\end{pmatrix}\\
&\congsim\begin{pmatrix}0 & 1\\ 1 & 0\end{pmatrix}.
\end{split}
\end{align}
Here, in the congruence transformations used in the first, third, and fourth line, the matrices multiplied on the left are 
$$\begin{pmatrix}1 & -\eta\\ 0 & 1\end{pmatrix},\begin{pmatrix}1 & 0\\ -4b_2(1+8\xi-8b_2\eta)^{-1} & 1\end{pmatrix},\begin{pmatrix}1 & 0\\\ 0 & (1+8\xi-8b_2\eta)^{-1} \end{pmatrix},$$
respectively.
\item For $\xi,b_1,b_2\in\Z_2$, we claim congruent equivalence $\begin{pmatrix}2 & 1\\ 1 & 2\end{pmatrix}\congsim\begin{pmatrix}2+8b_1 & 1+8\xi\\ 1+8\xi & 2+8b_2\end{pmatrix}$. In fact, we can find $\eta_1\in\Z_2,\eta_2\in\Z_2^\times$ such that $$(2+8b_2)\eta_1^2-2(1+8\xi)\eta_1+8b_1=0,\quad\det\begin{pmatrix}2+8b_1 & 1+8\xi\\ 1+8\xi & 2+8b_2\end{pmatrix}\cdot\eta_2^2=\det\begin{pmatrix}2 & 1\\ 1 & 2\end{pmatrix}=3.$$
Therefore, we have
\begin{align}
\begin{split}
\begin{pmatrix}2+8b_1 & 1+8\xi\\ 1+8\xi & 2+8b_2\end{pmatrix}&\congsim\begin{pmatrix}(2+8b_2)\eta_1^2-2(1+8\xi)\eta_1+2+8b_1 & 1+8\xi-8b_2\eta_1\\ 1+8\xi-8b_2\eta_1 & 2+8b_2\end{pmatrix}\\
&=\begin{pmatrix}2 & 1+8\xi-8b_2\eta_1\\ 1+8\xi-8b_2\eta_1 & 2+8b_2\end{pmatrix}\\
&\congsim\begin{pmatrix}2 & (1+8\xi-8b_2\eta_1)\eta_2\\ (1+8\xi-8b_2\eta_1)\eta_2 & (2+8b_2)\eta_2^2\end{pmatrix}\\
&\congsim\begin{pmatrix}2 & 1\\ 1 & 2\end{pmatrix}.
\end{split}
\end{align}
Here, in the congruence transformations used in the first, third, and fourth line, the matrices multiplied on the left are 
$$\begin{pmatrix}1 & -\eta_1\\ 0 & 1\end{pmatrix},\begin{pmatrix}1 & 0\\ 0 & \eta_2\end{pmatrix},\begin{pmatrix}1 & 0\\\ (1-(1+8\xi-8b_2\eta_1)\eta_2)/2 & 1\end{pmatrix},$$
respectively.
\end{enumerate}
Now, we complete the proof.
\end{proof}

It is clear that when two nonsingular matrices in $\Sym_n(\Z_p)$ are congruent, they induce the same paired $p$-group. However, one has to be careful for the converse assertion, because two nonsingular matrices in $\Sym_n(\Z_p)$ might induce the same paired $p$-group even if they do not lie on the same orbit under the congruence transformation of $\GL_n(\Z_p)$. Let us consider the following four sets of block matrices in \eqref{eq: simplest matrix sym 2}, where $p=2$:
\begin{equation}\label{eq: same pairing block}
\mathcal{T}_1=\{(2),(6),(10),(14)\},\mathcal{T}_2=\left\{\begin{pmatrix}0 & 2\\ 2 & 0\end{pmatrix},\begin{pmatrix}4 & 2\\ 2 & 4\end{pmatrix}\right\},\mathcal{T}_3=\{(4),(20)\},\mathcal{T}_4=\{(12,28)\}.
\end{equation}
The readers can verify that for $1\le i\le 4$, the matrices in $\mathcal{T}_i$ are pairwise non-congruent but still induce the same symmetric pairing.


\begin{prop}\label{prop: same pairing for block congruence}
\begin{enumerate}
\item ($p$ odd) Let $\alpha,\alpha'\in p\Z_p$ be nonzero, such that $\alpha^{-1}-(\alpha')^{-1}\in\Z_p$. Then, there exists $u\in\Z_p^\times$ such that $\alpha'=\alpha u^2$.
\item ($p=2$) Let $M$ be one of the matrix in \eqref{eq: simplest matrix sym 2} with $d\ge 1$. Furthermore, suppose $M'$ is symmetric, nonsingular, has the same size as $M$, and $M^{-1}-(M')^{-1}$ has entries in $\Z_p$. In this case, if $M$ is contained in $\mathcal{T}_i$ for some $1\le i\le 4$, then $M'$ is congruent to a matrix in $\mathcal{T}_i$; otherwise, we have $M\congsim M'$.
\end{enumerate}
\end{prop}

\begin{proof}
When $p$ is odd, we have
$$\val(\frac{\alpha'}{\alpha}-1)=\val(\alpha^{-1}-(\alpha')^{-1})+\val(\alpha')\ge 1,$$
and therefore $\frac{\alpha'}{\alpha}\in1+p\Z_p$ is a square in $\Z_p^\times$.

When $p=2$, we proceed by a case-by-case discussion according to the congruence class of $M$. If $M\in\mathcal{T}_1$, the only entry of $(M')^{-1}$ lies in $1/2+\Z_2$. The cases in which this entry lies in 
$$1/2+4\Z_2,3/2+4\Z_2,5/2+4\Z_2,7/2+4\Z_2$$
correspond to the congruence classes $(2),(6),(10),(14)$, respectively.

If $M\in\mathcal{T}_2$, then we have $M'\in\begin{pmatrix}0 & 1/2 \\ 1/2 & 0\end{pmatrix}+\Sym_2(\Z_2)$. When $(M')^{-1}$ lies in 
$$\left\{\begin{pmatrix}1+2b_1 & 1/2+b_2 \\ 1/2+b_2 & 1+2b_3\end{pmatrix}:b_1,b_2,b_3\in\Z_2\right\},$$
we have $(M')^{-1}\congsim\begin{pmatrix}1 & 1/2 \\ 1/2 & 1\end{pmatrix}$, and therefore $M'\congsim\begin{pmatrix}4 & 2 \\ 2 & 4\end{pmatrix}$; otherwise, we have $M'\congsim\begin{pmatrix}0 & 2 \\ 2 & 0\end{pmatrix}$.

If $M\in\mathcal{T}_3$, the only entry of $(M')^{-1}$ lies in $1/4+\Z_2$. The cases in which this entry lies in $1/4+2\Z_2,5/4+2\Z_2$ correspond to the congruence classes $(4),(20)$, respectively.

If $M\in\mathcal{T}_4$, the only entry of $(M')^{-1}$ lies in $3/4+\Z_2$. The cases in which this entry lies in $3/4+2\Z_2,7/4+2\Z_2$ correspond to the congruence classes $(12),(28)$, respectively.

If $M=\begin{pmatrix}0 & 4\\ 4 & 0\end{pmatrix}$, then $(M')^{-1}\in\begin{pmatrix}0 & 1/4\\ 1/4 & 0\end{pmatrix}+\Sym_2(\Z_2)$. Therefore, we have $(M')^{-1}\congsim\begin{pmatrix}0 & 1/4\\ 1/4 & 0\end{pmatrix}$, and $M'\congsim M$.

If $M=\begin{pmatrix}8 & 4\\ 4 & 8\end{pmatrix}$, then 
$$(M')^{-1}\in\begin{pmatrix}1/6 & -1/12\\ -1/12 & 1/6\end{pmatrix}+\Sym_2(\Z_2)=\begin{pmatrix}1/2 & 1/4\\ 1/4 & 1/2\end{pmatrix}+\Sym_2(\Z_2).$$ 
Therefore, we have $(M')^{-1}\congsim\begin{pmatrix}1/2 & 1/4\\ 1/4 & 1/2\end{pmatrix}\congsim M^{-1}$, and $M'\congsim M$.

If $M$ is one of the matrix in \eqref{eq: simplest matrix sym 2} such that $d\ge 3$, we deduce that $2^d M^{-1}$ is invertible with entries in $\Z_p$, $\Dep_2(\Cok(2^d M^{-1}))=0$, and that $2^d M^{-1}-2^d (M')^{-1}$ has entries in $2^d\Z_p$. By \Cref{prop: locally congruent equivalence at p}, we have $2^d M^{-1}\congsim 2^d(M')^{-1}$, and therefore $M\congsim M'$. This completes the discussion and ends the proof.
\end{proof}

\begin{prop}\label{thm: add invertible to same cokernel with pairing}
Let $n_1,n_2\ge 1$, and $M_{n_1}\in\Sym_{n_1}(\Z_p),M_{n_2}'\in\Sym_{n_2}(\Z_p)$ be nonsingular. Then, the following are equivalent:
\begin{enumerate}
\item We have $(\Cok(M_{n_1}),\langle\cdot,\cdot\rangle)\simeq(\Cok(M_{n_2}'),\langle\cdot,\cdot\rangle)$. \label{item: same cokernel with pairing} 
\item There exists $m\le \min\{n_1,n_2\}$, and
$$M_{n_1}^{\inv}\in\Sym_{n_1-m}(\Z_p)\cap\GL_{n_1-m}(\Z_p),M_{n_1}^{\nil}\in\Sym_m(p\Z_p),$$
$$M_{n_2}'^{\inv}\in\Sym_{n_2-m}(\Z_p)\cap\GL_{n_2-m}(\Z_p),M_{n_2}'^{\nil}\in\Sym_m(p\Z_p),$$
such that $(M_{n_1}^{\nil})^{-1}-(M_{n_2}'^{\nil})^{-1}\in\Sym_m(\Z_p)$, and
$$M_{n_1}\congsim\begin{pmatrix}M_{n_1}^{\inv} & 0 \\ 0 & M_{n_1}^{\nil}\end{pmatrix},\quad M_{n_2}'\congsim\begin{pmatrix}M_{n_2}'^{\inv} & 0 \\ 0 & M_{n_2}'^{\nil}\end{pmatrix}.$$ \label{item: difference of inverse}
\item There exists $n\ge\max\{n_1,n_2\}$, and
$$M_{n-n_1}^c\in\Sym_{n-n_1}(\Z_p)\cap\GL_{n-n_1}(\Z_p),M_{n-n_2}'^c\in\Sym_{n-n_2}(\Z_p)\cap\GL_{n-n_2}(\Z_p),$$ 
such that $\begin{pmatrix}M_{n-n_1}^c & 0 \\ 0 & M_{n_1}\end{pmatrix}\congsim\begin{pmatrix}M_{n-n_2}'^c & 0 \\ 0 & M_{n_2}'\end{pmatrix}$. \label{item: add matrix congruent}
\end{enumerate}
\end{prop}

\begin{proof}[Proof of \Cref{thm: add invertible to same cokernel with pairing}]
First, notice that adding a new invertible block at the upper-left corner and congruent transformations do not change the cokernel with pairing. Therefore, we have \eqref{item: difference of inverse} implies \eqref{item: same cokernel with pairing}, and \eqref{item: add matrix congruent} implies \eqref{item: same cokernel with pairing}.

Next, we prove that \eqref{item: same cokernel with pairing} implies \eqref{item: difference of inverse}. Suppose that we already have $(\Cok(M_{n_1}),\langle\cdot,\cdot\rangle)\simeq(\Cok(M_{n_2}'),\langle\cdot,\cdot\rangle)$. In this case, the two matrices $M_{n_1}/p\in\Sym_{n_1}(\F_p),M_{n_2}'/p\in\Sym_{n_2}(\F_p)$ have the same corank, which is exactly the $m$ we take. Applying \Cref{lem: full rank principal minor}, there exists a $(n_1-m)\times (n_1-m)$ principal minor of $M_{n_1}$ with determinant in $\Z_p^\times$. There is no loss of generality to assume that the $(n_1-m)\times (n_1-m)$ upper-left block is invertible, so that
$$M_{n_1}=\begin{pmatrix}B & C \\ C^T & D\end{pmatrix}, B\in\Sym_{n_1-m}(\Z_p)\cap\GL_{n_1-m}(\Z_p),C\in\Mat_{(n_1-m)\times m}(\Z_p),D\in\Sym_{m}(\Z_p).$$
In this case, we have
$$M_{n_1}\congsim\begin{pmatrix}B & 0 \\ 0 & D-C^TB^{-1}C\end{pmatrix},$$
where $D-C^TB^{-1}C\in\Sym_{m}(p\Z_p)$. Let $M_{n_1}^{\inv}=B$, and $M_{n_1}^{\nil}=D-C^TB^{-1}C\in\Sym_{m}(p\Z_p)$. Moreover, $M_{n_1}^{\nil}$ is nonsingular, and 
$(\Cok(M_{n_1}),\langle\cdot,\cdot\rangle)\simeq(\Cok(M_{n_1}^{\nil}),\langle\cdot,\cdot\rangle)$. We construct $M_{n_2}'^{\inv},M_{n_2}'^{\nil}$ in the analogous way, so that 
$$(\Cok(M_{n_1}^{\nil}),\langle\cdot,\cdot\rangle)\simeq(\Cok(M_{n_1}),\langle\cdot,\cdot\rangle)\simeq(\Cok(M_{n_2}'),\langle\cdot,\cdot\rangle)\simeq(\Cok(M_{n_2}'^{\nil}),\langle\cdot,\cdot\rangle).$$
Applying \Cref{prop: same pairing matrix of same size}, we can adjust $M_{n_2}'^{\nil}$ by a congruent transformation (with a slight abuse of notation, we denote the adjusted matrix again by $M_{n_2}'^{\nil}$), so that $(M_{n_1}^{\nil})^{-1}-(M_{n_2}'^{\nil})^{-1}\in\Sym_n(\Z_p)$.

In the end, we prove that \eqref{item: difference of inverse} implies \eqref{item: add matrix congruent}. Suppose that the property in \eqref{item: difference of inverse} holds. It suffices to show that we can adjoin invertible matrices of the same size to the upper-left corners of $M_{n_1}^{\nil}$ and $M_{n_2}'^{\nil}$, respectively, so that the resulting enlarged matrices are congruent.

Applying the same congruence transformation to $M_{n_1}^{\nil}$ and $M_{n_2}'^{\nil}$, there is no loss of generality to assume that $M_{n_1}^{\nil}$ is already a diagonal block matrix given in \Cref{prop: simplest matrix_sym}. In this case, $(M_{n_1}^{\nil})^{-1}$ is block-diagonal with the same block structure (i.e., the diagonal blocks occur in the same positions and have the same sizes). Turning to $(M_{n_2}'^{\nil})^{-1}$, we can use the diagonal blocks in the corresponding positions as pivots to clear the off-block-diagonal entries, and under such simultaneous left-right elementary matrix operations, each entry in the diagonal blocks changes only by an element of $\Z_p$. 

Now, we have adjusted $M_{n_2}'^{\nil}$, so that it not only satisfies $(M_{n_1}^{\nil})^{-1}-(M_{n_2}'^{\nil})^{-1}\in\Sym_m(\Z_p)$, but also has the same block structure as $M_{n_1}^{\nil}$. Therefore, it suffices to treat the diagonal blocks in each corresponding position separately. By \Cref{prop: same pairing for block congruence}, if $p$ is odd, we already have $M_{n_1}^{\nil}\congsim M_{n_2}'^{\nil}$; and if $p=2$, we only need to deal with matrices in the sets $\mathcal{T}_1,\mathcal{T}_2,\mathcal{T}_3,\mathcal{T}_4$. For the rest of the proof, we will treat these four sets one by one. Specifically, in each case we enlarge the matrix by adding an invertible block in the upper-left corner, and then derive the required congruence relation.

For the matrices in $\mathcal{T}_1$, we have
$$\begin{pmatrix}1 & 0 \\ 0 & 2\end{pmatrix}\congsim\begin{pmatrix}3 & 2 \\ 2 & 2\end{pmatrix}\congsim\begin{pmatrix}3 & 0 \\ 0 & 2/3\end{pmatrix}\congsim \begin{pmatrix}3 & 0 \\ 0 & 6\end{pmatrix},$$
$$\begin{pmatrix}9 & 0 \\ 0 & 6\end{pmatrix}\congsim\begin{pmatrix}15 & 6 \\ 6 & 6\end{pmatrix}\congsim\begin{pmatrix}15 & 0 \\ 0 & 18/5\end{pmatrix}\congsim \begin{pmatrix}15 & 0 \\ 0 & 10\end{pmatrix},$$
$$\begin{pmatrix}5 & 0 \\ 0 & 10\end{pmatrix}\congsim\begin{pmatrix}15 & 10 \\ 10 & 10\end{pmatrix}\congsim\begin{pmatrix}15 & 0 \\ 0 & 10/3\end{pmatrix}\congsim \begin{pmatrix}15 & 0 \\ 0 & 14\end{pmatrix},$$
$$\begin{pmatrix}1 & 0 \\ 0 & 14\end{pmatrix}\congsim\begin{pmatrix}15 & 14 \\ 14 & 14\end{pmatrix}\congsim\begin{pmatrix}15 & 0 \\ 0 & 14/15\end{pmatrix}\congsim \begin{pmatrix}15 & 0 \\ 0 & 2\end{pmatrix}.$$
Therefore, we deduce that
$$\begin{pmatrix}1 & 0 & 0\\ 0 & 9 & 0\\ 0 & 0 & 2\end{pmatrix}\congsim\begin{pmatrix}3 & 0 & 0\\ 0 & 9 & 0\\ 0 & 0 & 6\end{pmatrix}\congsim\begin{pmatrix}3 & 0 & 0\\ 0 & 15 & 0\\ 0 & 0 & 10\end{pmatrix},$$
$$\begin{pmatrix}9 & 0 & 0\\ 0 & 5 & 0\\ 0 & 0 & 6\end{pmatrix}\congsim\begin{pmatrix}15 & 0 & 0\\ 0 & 5 & 0\\ 0 & 0 & 10\end{pmatrix}\congsim\begin{pmatrix}15 & 0 & 0\\ 0 & 15 & 0\\ 0 & 0 & 14\end{pmatrix}.$$

For the matrices in $\mathcal{T}_2$, we have 
$$\begin{pmatrix}1 & 0 & 0\\ 0 & 0 & 2\\ 0 & 2 & 0\end{pmatrix}\congsim\begin{pmatrix}1 & 0 & 0\\ 0 & 0 & 2\\ 0 & 2 & 4\end{pmatrix}\congsim\begin{pmatrix}1 & 2 & 0\\ 2 & 4 & 2\\ 0 & 2 & 4\end{pmatrix}\congsim\begin{pmatrix}-1/3 & 0 & 0\\ 0 & 4 & 2\\ 0 & 2 & 4\end{pmatrix}\congsim\begin{pmatrix}5 & 0 & 0\\ 0 & 4 & 2\\ 0 & 2 & 4\end{pmatrix}.$$

For the matrices in $\mathcal{T}_3$, we have
$$\begin{pmatrix}1 & 0 \\ 0 & 4\end{pmatrix}\congsim\begin{pmatrix}5 & 4 \\ 4 & 4\end{pmatrix}\congsim\begin{pmatrix}5 & 0 \\ 0 & 4/5\end{pmatrix}\congsim \begin{pmatrix}5 & 0 \\ 0 & 20\end{pmatrix}.$$

For the matrices in $\mathcal{T}_4$, we have
$$\begin{pmatrix}1 & 0 \\ 0 & 12\end{pmatrix}\congsim\begin{pmatrix}13 & 12 \\ 12 & 12\end{pmatrix}\congsim\begin{pmatrix}13 & 0 \\ 0 & 12/13\end{pmatrix}\congsim \begin{pmatrix}13 & 0 \\ 0 & 28\end{pmatrix}.$$

This completes the discussion and ends the proof.
\end{proof} 

\begin{cor}\label{cor: new added row and column cokernel with pairing_sym}
Let $n_1,n_2\ge 1$, and $M_{n_1}\in\Sym_{n_1}(\Z_p),M_{n_2}'\in\Sym_{n_2}(\Z_p)$ be nonsingular, such that $(\Cok(M_{n_1}),\langle\cdot,\cdot\rangle)\simeq(\Cok(M_{n_2}'),\langle\cdot,\cdot\rangle)$. Let $z,\xi_1,\xi_2,\ldots$ be i.i.d. Haar-distributed in $\Z_p$, and let $\bm{\xi}_1:=(\xi_1,\ldots,\xi_{n_1})\in\Z_p^{n_1},\bm{\xi}_2:=(\xi_1,\ldots,\xi_{n_2})\in\Z_p^{n_2}$. Then, we have 
$$\left(\Cok\begin{pmatrix}M_{n_1} & \bm{\xi}_1 \\ \bm{\xi_1}^T &z \end{pmatrix},\langle\cdot,\cdot\rangle\right)\stackrel{d}{=}\left(\Cok\begin{pmatrix}M_{n_2}' & \bm{\xi}_2 \\ \bm{\xi_2}^T &z \end{pmatrix},\langle\cdot,\cdot\rangle\right).$$
\end{cor}

\begin{proof}
By \Cref{thm: add invertible to same cokernel with pairing}, we can find $n\ge\max\{n_1,n_2\}$, and
$$M_{n-n_1}^c\in\Sym_{n-n_1}(\Z_p)\cap\GL_{n-n_1}(\Z_p),M_{n-n_2}'^c\in\Sym_{n-n_2}(\Z_p)\cap\GL_{n-n_2}(\Z_p),$$ 
such that $\begin{pmatrix}M_{n-n_1}^c & 0 \\ 0 & M_{n_1}\end{pmatrix}\congsim\begin{pmatrix}M_{n-n_2}'^c & 0 \\ 0 & M_{n_2}'\end{pmatrix}$. Let $$\bm{\xi}_1^c=(\xi_{n_1+1},\ldots,\xi_n)\in\Z_p^{n-n_1},\bm{\xi}_2^c(\xi_{n_2+1},\ldots,\xi_n)\in\Z_p^{n-n_2},$$ 
so that $(\bm{\xi}_1^c,\bm{\xi}_1),(\bm{\xi}_2^c,\bm{\xi}_2)\in\Z_p^n$ are Haar-distributed. Notice that the Haar probability measure over $\Z_p^n$ is invariant under the action of $\GL_n(\Z_p)$. Therefore, we have
\begin{align}
\begin{split}
\left(\Cok\begin{pmatrix}M_{n_1} & \bm{\xi}_1 \\ \bm{\xi_1}^T &z \end{pmatrix},\langle\cdot,\cdot\rangle\right)&\stackrel{d}{=}\left(\Cok\begin{pmatrix}M_{n-n_1}^c & 0 & \bm{\xi_1}^c\\
0 & M_{n_1} & \bm{\xi}_1 \\ (\bm{\xi_1}^c)^T & \bm{\xi_1}^T & z \end{pmatrix},\langle\cdot,\cdot\rangle\right)\\
&\stackrel{d}{=}\left(\Cok\begin{pmatrix}M_{n-n_2}'^c & 0 & \bm{\xi_2}^c\\
0 & M_{n_2}' & \bm{\xi}_2 \\ (\bm{\xi_2}^c)^T & \bm{\xi_2}^T & z \end{pmatrix},\langle\cdot,\cdot\rangle\right)\\
&\stackrel{d}{=}\left(\Cok\begin{pmatrix}M_{n_2}' & \bm{\xi}_2 \\ \bm{\xi_2}^T &z \end{pmatrix},\langle\cdot,\cdot\rangle\right).
\end{split}
\end{align}
Here, the first line holds because we can use the block matrix in the upper-left corner to eliminate $\bm{\xi}_1^c$ and $(\bm{\xi}_1^c)^T$ on the new added column and row, and similarly, the third line also holds.
\end{proof}

Based on \Cref{cor: new added row and column cokernel with pairing_sym}, we are able to give the following definition.

\begin{defi}
Let $p$ be a prime number, and $(G^{(1)},\langle\cdot,\cdot\rangle_{G^{(1)}}),(G^{(2)},\langle\cdot,\cdot\rangle_{G^{(2)}})$ be two paired $p$-groups. Take a nonsingular matrix $M_n\in\Sym_n(\Z_p)$ for some $n\ge 1$, such that $(\Cok(M_n),\langle\cdot,\cdot\rangle)\simeq(G^{(1)},\langle\cdot,\cdot\rangle_{G^{(1)}})$. Let $z,\xi_1,\xi_2,\ldots$ be i.i.d. Haar-distributed in $\Z_p$, and let $\bm{\xi}:=(\xi_1,\xi_2,\ldots.\xi_n)\in\Z_p^n$. Then, the \emph{transition probability} from $(G^{(1)},\langle\cdot,\cdot\rangle_{G^{(1)}})$ to $(G^{(2)},\langle\cdot,\cdot\rangle_{G^{(2)}})$ is defined as
\begin{multline}
\mathbf{P}\left(\left(G^{(2)},\langle\cdot,\cdot\rangle_{G^{(2)}}\right)\bigg|\left(G^{(1)},\langle\cdot,\cdot\rangle_{G^{(1)}}\right)\right):=\\
\mathbf{P}\left(\begin{pmatrix}M_n & \bm{\xi} \\ \bm{\xi}^T &z \end{pmatrix}\text{ nonsingular, }\left(\Cok\begin{pmatrix}M_n & \bm{\xi} \\ \bm{\xi}^T &z \end{pmatrix},\langle\cdot,\cdot\rangle\right)\simeq(G^{(2)},\langle\cdot,\cdot\rangle_{G^{(2)}})\right).
\end{multline}
\end{defi}

Indeed, by \Cref{cor: new added row and column cokernel with pairing_sym}, the transition  probability $\mathbf{P}\left((G^{(2)},\langle\cdot,\cdot\rangle_{G^{(2)}})\mid(G^{(1)},\langle\cdot,\cdot\rangle_{G^{(1)}})\right)$ does not rely on the matrix $M_n$ we choose, and therefore is well defined.

Although our approach does not involve the surjection moment method, we still transfer matrices in $\Sym_n(\Z)$ to $\Sym_n(\Z/a\Z)$ by taking the reduction mod $a$. Suppose that we have chosen $(G,\langle\cdot,\cdot\rangle_G),P$ in the same way as in \Cref{thm: exponential convergence_sym}. Then, the integer $a$ we pick has to be large enough to distinguish the matrices in $\Sym_n(\Z)$ that generate the pairing $(G,\langle\cdot,\cdot\rangle_G)$ under the prime set $P$. When we do not care about the pairing structure, it suffices to take $a=\prod_{p\in P}p^{\Dep_p(G)+1}$, see the proof of \cite[Corollary 9.2]{wood2017distribution} for instance. Nevertheless, one has to be more careful when taking the pairing into account, which is the main concern of the following proposition. 

\begin{prop}\label{cor: globally same pairing}
Let $(G,\langle\cdot,\cdot\rangle_G)$ be a paired abelian group, and $P=\{p_1,\ldots,p_l\}$ be a finite set of primes that includes all those that divide $\#G$. Let $a:=p_1^{e_{p_1}}\cdots p_l^{e_{p_l}}$, where for all $1\le i\le l$, we have
$$e_{p_i}=\begin{cases}
\Dep_{p_i}(G)+3 & p_i=2\\
\Dep_{p_i}(G)+1 & p_i>2.
\end{cases}.$$
Let $M_n\in\Sym_n(\Z)$ be nonsingular, such that $(\Cok(M_n)_P,\langle\cdot,\cdot\rangle_P)\simeq (G,\langle\cdot,\cdot\rangle_G)$. Then for all $\Delta M_n\in\Sym_n(a\Z)$, we also have $M_n+\Delta M_n$ is nonsingular, and
$$(\Cok(M_n+\Delta M_n)_P,\langle\cdot,\cdot\rangle_P)\simeq (G,\langle\cdot,\cdot\rangle_G).$$
\end{prop}

\begin{proof}
This immediately follows from \Cref{prop: locally congruent equivalence at p} since everything naturally factors over $p_1,\ldots,p_l$.
\end{proof}

We now turn to symmetric matrices with entries in $\Z/a\Z$. Let $H_n\in\Sym_n(\Z/a\Z)$, then we naturally get a $\Z/a\Z$-module $\Cok(H_n)=(\Z/a\Z)^n/H_n(\Z/a\Z)^n$. Nevertheless, it can be impractical to talk about the pairing structure of $\Cok(H_n)$. Let us consider this simple example: take $n=1$, $a=3$ is a prime number, then the matrix $H_1=(0)$ gives the cokernel $\Z/3\Z=\F_3$, but we cannot naturally equip this cokernel with a perfect pairing structure. Thus, we need to introduce the following definition, which is inspired by \Cref{thm: add invertible to same cokernel with pairing}.

\begin{defi}\label{defi: quasi pairing of Cok}
Let $p$ be a prime number, and $e_p\ge 1$. We define an equivalence relation over the set of symmetric matrices with entries in $\Z/p^{e_p}\Z$ and use the notation $\Cok_*^{(p^{e_p})}$ for an equivalence class, which we call a cokernel with \emph{quasi-pairing} as follows. Given two matrices $H_{n_1}\in\Sym_{n_1}(\Z/p^{e_p}\Z),H_{n_2}'\in\Sym_{n_2}(\Z/p^{e_p}\Z)$, we say that the cokernels with quasi-pairing induced by $H_{n_1},H_{n_2}'$ are isomorphic, and write
$$\Cok_*^{(p^{e_p})}(H_{n_1})\simeq\Cok_*^{(p^{e_p})}(H_{n_2}'),$$ if there exist $n\ge\max\{n_1,n_2\}$, and
$$H_{n-n_1}^c\in\Sym_{n-n_1}(\Z/p^{e_p}\Z)\cap\GL_{n-n_1}(\Z/p^{e_p}\Z),H_{n-n_2}'^c\in\Sym_{n-n_2}(\Z/p^{e_p}\Z)\cap\GL_{n-n_2}(\Z/p^{e_p}\Z),$$ 
such that $\begin{pmatrix} H_{n-n_1}^c & 0\\ 0 & H_{n_1}\end{pmatrix}$ is congruent to $\begin{pmatrix} H_{n-n_2}'^c & 0\\ 0 & H_{n_2}'\end{pmatrix}$ under the action of $\GL_n(\Z/p^{e_p}\Z)$. Furthermore, for $a\ge 2$ with the prime factorization $a=p_1^{e_{p_1}}\cdots p_l^{e_{p_l}}$, we also define the cokernel with quasi-pairing $\Cok_*^{(a)}$ as follows. Given two matrices $H_{n_1}\in\Sym_{n_1}(\Z/a
\Z),H_{n_2}'\in\Sym_{n_2}(\Z/a\Z)$, we say that the cokernels with quasi-pairing induced by $H_{n_1},H_{n_2}'$ are isomorphic, and write
$$\Cok_*^{(a)}(H_{n_1})\simeq\Cok_*^{(a)}(H_{n_2}'),$$ 
if $\Cok_*^{(p_i^{k_i})}(H_{n_1}/p_i^{k_i})\simeq\Cok_*^{(p_i^{k_i})}(H_{n_2}'/p_i^{k_i})$ for all $1\le i\le l$.
\end{defi}

$\Cok_*^{(a)}$ is not a true pairing structure, since we do not equip the cokernel with a bilinear form. In particular, when $H_n,H_n'\in\Sym_n(\Z/a\Z)$ are congruent under the action of $\GL_n(\Z/a\Z)$, we have $\Cok_*^{(a)}(H_n)\simeq\Cok_*^{(a)}(H_n')$. Also, we must have $\Cok_*^{(a)}(H)\simeq\Cok_*^{(a)}\begin{pmatrix}H^c & 0\\ 0 & H\end{pmatrix}$ when the symmetric matrix $H^c$ is invertible. Furthermore, for $H_{n_1}\in\Sym_{n_1}(\Z/a\Z),H_{n_2}'\in\Sym_{n_2}(\Z/a\Z)$ such that $\Cok_*^{(a)}(H_{n_1})\simeq\Cok_*^{(a)}(H_{n_2}')$, we have an isomorphism of finite abelian groups $\Cok(H_{n_1})\simeq\Cok(H_{n_2}')$. 

\begin{prop}\label{prop: same quasi pairing implies same pairing}
Let $(G,\langle\cdot,\cdot\rangle_G)$ be a paired abelian group, and $P=\{p_1,\ldots,p_l\}$ be a finite set of primes that includes all those that divide $\#G$. Let $a\ge 2$ be an integer with prime factorization $a=p_1^{e_{p_1}}\cdots p_l^{e_{p_l}}$, where for all $1\le i\le l$, we have $e_{p_i}\ge\Dep_{p_i}(G)+3$ if $p_i=2$, and $e_{p_i}\ge\Dep_{p_i}(G)+1$ if $p_i>2$. Suppose $M_{n_1}\in\Sym_{n_1}(\Z)$ such that $(\Cok(M_{n_1})_P,\langle\cdot,\cdot\rangle_P)\simeq (G,\langle\cdot,\cdot\rangle_G)$. Then, for all $M_{n_2}'\in\Sym_{n_2}(\Z)$ such that $\Cok_*^{(a)}(M_{n_1}/a)\simeq\Cok_*^{(a)}(M_{n_2}'/a)$, the matrix $M_{n_2}'$ is also nonsingular, and $(\Cok(M_{n_2}')_P,\langle\cdot,\cdot\rangle_P)\simeq (G,\langle\cdot,\cdot\rangle_G)$.
\end{prop}

\begin{proof}
Notice that the primes $p_1,\ldots,p_l$ naturally factors. Thus, there is no loss of generality to assume $a=p_1^{e_{p_1}}$. Since $\Cok_*^{(p_1^{e_{p_1}})}(M_{n_1}/p_1^{e_{p_1}})\simeq\Cok_*^{(p_1^{e_{p_1}})}(M_{n_2}'/p_1^{e_{p_1}})$, there exists $n\ge\max\{n_1,n_2\}$, $\Delta M\in\Sym_n(p_1^{e_{p_1}}\Z_{p_1})$, and
$$M_{n-n_1}^c\in\Sym_{n-n_1}(\Z_{p_1})\cap\GL_{n-n_1}(\Z_{p_1}),M_{n-n_2}'^c\in\Sym_{n-n_2}(\Z_p)\cap\GL_{n-n_2}(\Z_{p_1}),$$
such that $\begin{pmatrix}M_{n-n_1}^c & 0 \\ 0 & M_{n_1}\end{pmatrix}\congsim\begin{pmatrix}M_{n-n_2}'^c & 0 \\ 0 & M_{n_2}'\end{pmatrix}+\Delta M$. Thus, the matrix $\begin{pmatrix}M_{n-n_2}'^c & 0 \\ 0 & M_{n_2}'\end{pmatrix}+\Delta M$ is nonsingular, and induces the same paired $p_1$-group $(G,\langle\cdot,\cdot\rangle_G)$. By \Cref{cor: globally same pairing}, the matrix $\begin{pmatrix}M_{n-n_2}'^c & 0 \\ 0 & M_{n_2}'\end{pmatrix}$ is nonsingular, and therefore $M_{n_2}'$ is nonsingular. Moreover, 
$$(\Cok(M_{n_2}')_{p_1},\langle\cdot,\cdot\rangle_{p_1})\simeq\left(\Cok\begin{pmatrix}M_{n-n_2}'^c & 0 \\ 0 & M_{n_2}'\end{pmatrix}_{p_1},\langle\cdot,\cdot\rangle_{p_1}\right)\simeq(G,\langle\cdot,\cdot\rangle_G).$$
This completes the proof.
\end{proof}

Based on \Cref{prop: same quasi pairing implies same pairing}, we have the following definition.

\begin{defi}\label{defi: quotient version of paired group}
Let $(G,\langle\cdot,\cdot\rangle_G)$ be a paired abelian group, and $P=\{p_1,\ldots,p_l\}$ be a finite set of primes that includes all those that divide $\#G$. Let $a\ge 2$ be an integer with prime factorization $a=p_1^{e_{p_1}}\cdots p_l^{e_{p_l}}$, where for all $1\le i\le l$, we have $e_{p_i}\ge\Dep_{p_i}(G)+3$ if $p_i=2$, and $e_{p_i}\ge\Dep_{p_i}(G)+1$ if $p_i>2$. We denote by $(G,\langle\cdot,\cdot\rangle_*^{(a)})$ the equivalence class under $\Cok_*^{(a)}$, given by the set of matrices
$$\{M/a:M\in\Sym_n(\Z)\text{ for some }n\ge 1, M\text{ nonsingular, } (\Cok(M)_P,\langle\cdot,\cdot\rangle_P)\simeq (G,\langle\cdot,\cdot\rangle_G)\}.$$
\end{defi}



The motivation for us to define the quasi-pairing $\Cok_*^{(a)}$ is the following proposition, which claims that when we add a new uniformly random row and column to a fixed $n\times n$ symmetric matrix over $\Z/a\Z$, the distribution of the cokernel with quasi-pairing of the new $(n+1)\times(n+1)$ matrix only relies on the cokernel with quasi-pairing of the original $n\times n$ corner.

\begin{prop}\label{prop: transition probability relies on cokernel_sym}
Let $a\ge 2$, and $n_1,n_2\ge 1$. Suppose that we have fixed matrices $M_{n_1}\in\Sym_{n_1}(\Z),M_{n_2}'\in\Sym_{n_2}(\Z)$ such that $\Cok_*^{(a)}(M_{n_1}/a)\simeq\Cok_*^{(a)}(M_{n_2}'/a)$.  Let $z,\xi_1,\xi_2,\ldots$ be independent and uniformly distributed in $\{0,1,\ldots,a-1\}$. Then we have
$$\Cok_*^{(a)}
\left(
\begin{array}{c|c}
M_{n_1}/a & 
\begin{matrix}
\xi_1/a \\
\vdots \\
\xi_{n_1}/a
\end{matrix}
\\ \hline
\begin{matrix}
\xi_1/a & \cdots & \xi_{n_1}/a
\end{matrix}
& z/a
\end{array}
\right)\stackrel{d}{=}\Cok_*^{(a)}
\left(
\begin{array}{c|c}
M_{n_2}'/a & 
\begin{matrix}
\xi_1/a \\
\vdots \\
\xi_{n_2}/a
\end{matrix}
\\ \hline
\begin{matrix}
\xi_1/a & \cdots & \xi_{n_2}/a
\end{matrix}
& z/a
\end{array}
\right).$$
\end{prop}

\begin{proof}
Notice that the primes that divide $a$ naturally factor. Thus, there is no loss to assume $a=p^{e_p}$. In this case, the proof follows the same as \Cref{cor: new added row and column cokernel with pairing_sym}.
\end{proof}



\begin{defi}
Let $a\ge 2$, and $n\ge 1$. Let $H_n\in\Sym_n(\Z/a\Z)$ be uniformly distributed. Denote by 
$$\mu_{n}^{\sym,(a)}:=\mathcal{L}(\Cok_*^{(a)}(H_n))$$
as the law of $\Cok_*^{(a)}(H_n)$.
\end{defi}

The following proposition shows that for the uniform random model, when the size $n$ of the random matrix goes to infinity, the cokernel with quasi-pairing gives exponential convergence.

\begin{prop}\label{prop: uniform model exponential convergence_sym} 
Let $a\ge 2$. When $n$ goes to infinity, $\mu_n^{\sym,(a)}$ weakly converges to a limiting probability measure, denoted by $\mu_\infty^{\sym,(a)}$. Furthermore, we have
$$D_{L_1}\left(\mu_n^{\sym,(a)},\mu_\infty^{\sym,(a)}\right)=O_a(\exp(-\Omega_a( n))).$$
\end{prop}

\begin{proof}
Suppose $a$ has prime factorization $a=p_1^{e_{p_1}}\cdots p_l^{e_{p_l}}$. In this case, we have
$$\mu_n^{\sym,(a)}=\mu_n^{\sym,(p_1^{e_{p_1}})}\times\cdots\times\mu_n^{\sym,(p_l^{e_{p_l}})}.$$
Notice that the assertions
$$\mu_n^{\sym,(p_i^{e_{p_i}})}\stackrel{d}{\rightarrow}\mu_\infty^{\sym,(p_i^{e_{p_i}})},D_{L_1}\left(\mu_n^{\sym,(p_i^{e_{p_i}})},\mu_\infty^{\sym,(p_i^{e_{p_i}})}\right)=O_{p_i^{e_{p_i}}}(\exp(-\Omega_{p_i^{e_{p_i}}}( n))),\quad 1\le i\le l,$$
directly implies the conclusion
$$\mu_n^{\sym,(a)}\stackrel{d}{\rightarrow}\mu_\infty^{\sym,(a)}:=\mu_\infty^{\sym,(p_1^{e_{p_1}})}\times\cdots\times\mu_\infty^{\sym,(p_l^{e_{p_l}})},D_{L_1}\left(\mu_n^{\sym,(a)},\mu_\infty^{\sym,(a)}\right)=O_a(\exp(-\Omega_a( n))).$$
Thus, there is no loss of generality to assume that $a$ only has one prime factor, i.e., $a=p^{e_p}$.

Now, we claim that for all integers $0\le k\le n$,
and fixed matrix $B_n\in\Sym_n(\F_p)$ such that $\corank(B_n)=k$, we have 
\begin{equation}\label{eq: distribution conditioned on residue}
\mathcal{L}\left(\Cok_*^{(p^{e_p})}(H_n)\bigg| H_n/p=B_n\right)=\mathcal{L}\left(\Cok_*^{(p^{e_p})}(H_k)\bigg| H_k\in\Sym_k(p\Z/p^{e_p}\Z)\right).
\end{equation}
Here, on the right hand side, $H_k$ is also uniformly distributed. To see this, first recall from \Cref{lem: full rank principal minor} that there exists a $(n-k)\times(n-k)$ principal minor of $B_n$ that is invertible. There is no loss of generality to assume that the $(n-k)\times(n-k)$ upper-left block of $B_n$ is invertible. In this case, when we uniformly randomly sample a matrix $H_n\in\Sym_n(\Z/p^{e_p}\Z)$ such that $H_n/p=B_n$, we can always use its $(n-k)\times(n-k)$ upper-left block to eliminate adjacent rows and columns, and the remaining $k\times k$ lower-right block is uniformly distributed in $\Sym_k(p\Z/p^{e_p}\Z)$. Thus, we have confirmed \eqref{eq: distribution conditioned on residue}. As a consequence, when $B_n$ runs through all the matrices in $\Sym_n(\F_p)$ of corank $k$, we deduce that 
$$\mathcal{L}\left(\Cok_*^{(p^{e_p})}(H_n)\bigg|\corank(H_n/p)=k\right)=\mathcal{L}\left(\Cok_*^{(p^{e_p})}(H_k)\bigg| H_k\in\Sym_k(p\Z/p^{e_p}\Z)\right).$$
Denote $\nu_n^{\sym,p}:=\mathcal{L}(\corank(H_n/p))$. Then, for all $0\le k\le n$ and fixed matrices $H_k'\in\Sym_k(p\Z/p^{e_p}\Z)$, we have
\begin{equation}\label{eq: distribution of quasi of size n}
\mu_n^{\sym,(p^{e_p})}\left(\Cok_*^{(p^{e_p})}(H_k')\right)=\nu_n^{\sym,p}(k)\cdot\mathbf{P}\left(\Cok_*^{(p^{e_p})}(H_k)\simeq\Cok_*^{(p^{e_p})}(H_k')\bigg| H_k\in\Sym_k(p\Z/p^{e_p}\Z)\right).
\end{equation}
Furthermore, by \cite[Theorem 4.1]{fulman2015stein}\footnote{In the notation of \cite[Theorem 4.1]{fulman2015stein}, we are taking $q=p$, $\mathcal{Q}_q=\nu_\infty^{\sym,p}$, and $\mathcal{Q}_{q,n}=\nu_n^{\sym,p}$.}, we have 
$$D_{L_1}(\nu_n^{\sym,p},\nu_\infty^{\sym,p})=O_p(\exp(-\Omega_p(n))),$$
where limit distribution $\nu_\infty^{\sym,p}$ is given by
$$\nu_\infty^{\sym,p}(k)=
\lim_{n\rightarrow\infty}\nu_n^{\sym,p}(k)=\frac{\prod_{i=0}^{\infty}(1-p^{-(2i+1)})}{\prod_{i=1}^k(p^i-1)},\quad\forall k\ge 0.$$
Therefore, when $n$ goes to infinity, $\mu_n^{\sym,(p^{e_p})}$ weakly converges to the limit distribution $\mu_\infty^{\sym,(p^{e_p})}$ given by
\begin{equation}\label{eq: distribution of quasi of size infinity}
\mu_\infty^{\sym,(p^{e_p})}\left(\Cok_*^{(p^{e_p})}(H_k')\right)=\nu_\infty^{\sym,p}(k)\cdot\mathbf{P}\left(\Cok_*^{(p^{e_p})}(H_k)\simeq\Cok_*^{(p^{e_p})}(H_k')\bigg| H_k\in\Sym_k(p\Z/p^{e_p}\Z)\right)
\end{equation}
for all $k\ge 0$ and fixed matrices $H_k'\in\Sym_k(p\Z/p^{e_p}\Z)$. Furthermore, comparing \eqref{eq: distribution of quasi of size n} and \eqref{eq: distribution of quasi of size infinity}, we have
\begin{equation}
D_{L_1}\left(\mu_n^{\sym,(p^{e_p})},\mu_\infty^{\sym,(p^{e_p})}\right)=D_{L_1}\left(\nu_n^{\sym,p},\nu_\infty^{\sym,p}\right)=O_p(\exp(-\Omega_p(n))).
\end{equation}
This completes the proof.
\end{proof}

%% file: The_cokernel_of_an_alternating_matrix_and_its_torsion_subgroup.tex
\section{The cokernel of an alternating matrix and its torsion subgroup}\label{sec: cok of alt matrix}

This section is the alternating analogue of \Cref{sec: cok of sym matrix}. Our goal is to clarify the structure of
cokernels arising from alternating matrices. While many ideas parallel the symmetric case, the
alternating setting differs in two main respects:

\begin{enumerate}
\item An alternating matrix $A_{2n}\in\Alt_{2n}(\Z)$ of even size can have nonzero determinant, and in this case $A_{2n}$ is nonsingular. However, an alternating integral matrix $A_{2n+1}\in\Alt_{2n+1}(\Z)$ of odd size must have determinant zero, and therefore $\Cok(A_{2n+1}):=\Z^{2n+1}/A_{2n+1}\Z^{2n+1}$ is an infinite abelian group. Hence, it is a more reasonable strategy to characterize $\Cok_{\tors}(A_{2n+1})$, which is the torsion subgroup of the cokernel.
\item For the symmetric case, two finite paired abelian groups $(G^{(1)},\langle\cdot,\cdot\rangle_{G^{(1)}}),(G^{(2)},\langle\cdot,\cdot\rangle_{G^{(2)}})$ might not be isomorphic even if $G^{(1)}\simeq G^{(2)}$ as abelian groups. However, for the alternating case, up to isomorphism, the pairing structure is totally determined by the abelian group itself. Although the term $|\Sp(G)|$ shows up in the expressions in \Cref{thm: exponential convergence_alt}, it can be regarded as a function over $S_P$. Hence, in the alternating case, no pairing structure is actually involved, which simplifies our discussion considerably.
\end{enumerate}

Let $n\ge 1$, and $A_n\in\Alt_n(\Z)$. Then $\Cok_{\tors}(A_n)$ is a finite abelian group. For a prime number $p$, the $p$-torsion subgroup of $\Cok(A_n)$ is the torsion subgroup of $\Z_p^n/A_n\Z_p^n$, where we regard $A_n$ as a matrix with entries in $\Z_p$. In particular, when $n$ is even and $A_n$ is nonsingular, the abelian group $\Cok(A_n)$ is finite, and therefore $\Cok(A_n)=\Cok_{\tors}(A_n)$. 

The following proposition gives the simplest form of an alternating matrix in $\Z_p$ under the congruence transformation of $\GL_n(\Z_p)$.

\begin{prop}\label{prop: simplest matrix_alt}
Let $n\ge 1$ be an integer, and $p$ be a prime number.
\begin{enumerate}
\item (Even size) Let $A_{2n}\in\Alt_{2n}(\Z_p)$ be nonsingular. Then $A_{2n}$ is congruent to a block diagonal matrix of the form
\begin{equation}\label{eq: simplest matrix alt even size}
\diag_{2n\times 2n}\left(\begin{pmatrix}0 & p^{\l_1} \\ -p^{\l_1} & 0\end{pmatrix},\ldots,\begin{pmatrix} 0 & p^{\l_n} \\ -p^{\l_n} & 0\end{pmatrix}\right),
\end{equation}
where $\l:=(\l_1,\ldots,\l_n)\in\Y$. In this case, we have 
$$\Cok(A_{2n})\simeq\bigoplus_{i=1}^n((\Z/p^{\l_i}\Z)\oplus(\Z/p^{\l_i}\Z))\in\mathcal{S}_p.$$ 
Furthermore, no two distinct matrices in \eqref{eq: simplest matrix alt even size} are congruent.
\item (Odd size) Let $A_{2n+1}\in\Alt_{2n+1}(\Z_p)$, such that $\corank(A)=1$. Then $A_{2n+1}$ is congruent to a block diagonal matrix of the form
\begin{equation}\label{eq: simplest matrix alt odd size}
\diag_{(2n+1)\times (2n+1)}\left(\begin{pmatrix}0 & p^{\l_1} \\ -p^{\l_1} & 0\end{pmatrix},\ldots,\begin{pmatrix} 0 & p^{\l_n} \\ -p^{\l_n} & 0\end{pmatrix},0\right),
\end{equation}
where $\l:=(\l_1,\ldots,\l_n)\in\Y$. In this case, we have 
$$\Cok_{\tors}(A_{2n+1})\simeq\bigoplus_{i=1}^n((\Z/p^{\l_i}\Z)\oplus(\Z/p^{\l_i}\Z))\in\mathcal{S}_p.$$ 
Furthermore, no two distinct matrices in \eqref{eq: simplest matrix alt odd size} are congruent.
\end{enumerate}
\end{prop}

\begin{proof}
See \cite[Proposition 2.1]{shen2024non}.
\end{proof}

\begin{prop}\label{prop: same cok and add new row column_alt}
Let $1\le n_1\le n_2$ be two integers, and $p$ be a prime number.
\begin{enumerate}
\item (Even size) Let $A_{2n_1}\in\Alt_{2n_1}(\Z_p),A_{2n_2}'\in\Alt_{2n_2}(\Z_p)$ be nonsingular, such that $\Cok(A_{2n_1})\simeq\Cok(A_{2n_2}')$. Then there exists $\l=(\l_1,\ldots,\l_{n_1})\in\Y$, such that
$$A_{2n_1}\congsim\diag_{2n_1\times 2n_1}\left(\begin{pmatrix}0 & p^{\l_1} \\ -p^{\l_1} & 0\end{pmatrix},\ldots,\begin{pmatrix} 0 & p^{\l_n} \\ -p^{\l_n} & 0\end{pmatrix}\right),$$
$$A_{2n_2}'\congsim\diag_{2n_2\times 2n_2}\left(\begin{pmatrix}0 & 1\\ -1 & 0\end{pmatrix},\ldots,\begin{pmatrix}0 & 1\\ -1 & 0\end{pmatrix},\begin{pmatrix}0 & p^{\l_1} \\ -p^{\l_1} & 0\end{pmatrix},\ldots,\begin{pmatrix} 0 & p^{\l_n} \\ -p^{\l_n} & 0\end{pmatrix}\right).$$
Furthermore, let $\bm{\xi}_1:=(\xi_1,\ldots,\xi_{2n_1})\in\Z_p^{2n_1}$, and $\bm{\xi}_2:=(\xi_1,\ldots,\xi_{2n_2})\in\Z_p^{2n_2}$. Then, we have
$$\Cok_{\tors}\begin{pmatrix}A_{2n_1} & \bm{\xi}_1 \\ -\bm{\xi_1}^T & 0 \end{pmatrix}\stackrel{d}{=}\Cok_{\tors}\begin{pmatrix}A_{2n_2}' & \bm{\xi}_2 \\ -\bm{\xi_2}^T & 0\end{pmatrix}.$$
\item (Odd size) Let $A_{2n_1+1}\in\Alt_{2n_1+1}(\Z_p),A_{2n_2+1}'\in\Alt_{2n_2+1}(\Z_p)$, such that $\corank(A_{2n_1+1})=\corank(A_{2n_2+1}')=1$, and $\Cok_{\tors}(A_{2n_1+1})\simeq\Cok_{\tors}(A_{2n_2+1}')$. Then there exists $\l=(\l_1,\ldots,\l_{n_1})\in\Y$, such that
$$A_{2n_1+1}\congsim\diag_{(2n_1+1)\times (2n_1+1)}\left(\begin{pmatrix}0 & p^{\l_1} \\ -p^{\l_1} & 0\end{pmatrix},\ldots,\begin{pmatrix} 0 & p^{\l_n} \\ -p^{\l_n} & 0\end{pmatrix},0\right),$$
\setlength{\arraycolsep}{3pt}
$$A_{2n_2+1}'\congsim\diag_{(2n_2+1)\times (2n_2+1)}\left(\begin{pmatrix}0 & 1\\ -1 & 0\end{pmatrix},\ldots,\begin{pmatrix}0 & 1\\ -1 & 0\end{pmatrix},\begin{pmatrix}0 & p^{\l_1} \\ -p^{\l_1} & 0\end{pmatrix},\ldots,\begin{pmatrix} 0 & p^{\l_n} \\ -p^{\l_n} & 0\end{pmatrix},0\right).$$
Furthermore, let $\bm{\xi}_1:=(\xi_1,\ldots,\xi_{2n_1+1})\in\Z_p^{2n_1+1}$, and $\bm{\xi}_2:=(\xi_1,\ldots,\xi_{2n_2+1})\in\Z_p^{2n_2+1}$. Then, we have
\setlength{\arraycolsep}{5pt}
$$\Cok\begin{pmatrix}A_{2n_1+1} & \bm{\xi}_1 \\ -\bm{\xi_1}^T & 0 \end{pmatrix}\stackrel{d}{=}\Cok\begin{pmatrix}A_{2n_2+1}' & \bm{\xi}_2 \\ -\bm{\xi_2}^T & 0\end{pmatrix}.$$
\end{enumerate}
\end{prop}

\begin{proof}
We will only prove the even size case since the case can be deduced routinely. Since $\Cok(A_{2n_1})\simeq\Cok(A_{2n_2}')\in\mathcal{S}_p$, they must be isomorphic to an abelian group of the form 
$$\bigoplus_{i=1}^{n_1}((\Z/p^{\l_i}\Z)\oplus(\Z/p^{\l_i}\Z)),\quad \l_1\ge\ldots\ge\l_{n_1}.$$
In this case, we can deduce their simplest forms under congruence transformations by \Cref{prop: simplest matrix_alt}. Consequently, 
$$A_{2n_2'}\congsim\begin{pmatrix}A_{2n_2-2n_1}^c & 0 \\ 0 & A_{2n_1}\end{pmatrix},A_{2n_2-2n_1}^c:=\diag_{(2n_2-2n_1)\times (2n_2-2n_1)}\left(\begin{pmatrix}0 & 1\\ -1 & 0\end{pmatrix},\ldots,\begin{pmatrix}0 & 1\\ -1 & 0\end{pmatrix}\right).$$
Denote by $\bm{\xi}_1^c:=(\xi_{2n_1+1},\ldots,\xi_{2n_2})\in\Z_p^{2n_2-2n_1}$, so that $(\bm{\xi}_1^c,\bm{\xi}_1)\in\Z_p^{2n_2}$ is Haar-distributed. We have
\begin{align}
\begin{split}
\Cok_{\tors}\begin{pmatrix}A_{2n_1} & \bm{\xi}_1 \\ -\bm{\xi_1}^T & 0 \end{pmatrix}&\stackrel{d}{=}\Cok_{\tors}\begin{pmatrix}
A_{2n_2-2n_1}^c  & 0 & \bm{\xi_1}^c\\
0 & A_{2n_1} & \bm{\xi}_1 \\ 
-(\bm{\xi_1}^c)^T & -\bm{\xi_1}^T & 0 \end{pmatrix}\\
&\stackrel{d}{=}\Cok_{\tors}\begin{pmatrix}A_{2n_2}' & \bm{\xi}_2 \\ -\bm{\xi_2}^T & 0\end{pmatrix}.
\end{split}
\end{align}
Here, the second line holds because the Haar measure over $\Z_p^{2n_2}$ is invariant under the action of $\GL_{2n_2}(\Z_p)$. This completes the proof.
\end{proof}

Based on \Cref{prop: same cok and add new row column_alt}, we are able to give the following definition.

\begin{defi}
Let $p$ be a prime number, and $G^{(1)},G^{(2)}\in\mathcal{S}_p$. Let $\xi_1,\xi_2,\ldots$ be i.i.d. Haar-distributed in $\Z_p$.
\begin{enumerate}
\item (From even to odd) Take a nonsingular matrix $A_{2n}\in\Alt_{2n}(\Z_p)$ for some $n\ge 1$, such that $\Cok(A_{2n})\simeq G^{(1)}$, and let $\bm{\xi}:=(\xi_1,\ldots,\xi_{2n})\in\Z_p^{2n}$. Then, the \emph{even-odd transition probability} from $G^{(1)}$ to $G^{(2)}$ is defined as
$$\mathbf{P}\left(G^{(2)},\text{ odd }\bigg| G^{(1)},\text{ even}\right):=\mathbf{P}\left(\Cok_{\tors}\begin{pmatrix}A_{2n} & \bm{\xi} \\ -\bm{\xi}^T &0 \end{pmatrix}\simeq G^{(2)}\right).$$
\item (From odd to even) Take a matrix $A_{2n+1}\in\Alt_{2n+1}(\Z_p)$ for some $n\ge 1$, such that $\corank(A_{2n+1})=1$, and $\Cok_{\tors}(A_{2n+1})\simeq G^{(1)}$. Let $\bm{\xi}:=(\xi_1,\ldots,\xi_{2n+1})\in\Z_p^{2n+1}$. Then, the \emph{odd-even transition probability} from $G^{(1)}$ to $G^{(2)}$ is defined as
$$\mathbf{P}\left(G^{(2)},\text{ even }\bigg| G^{(1)},\text{ odd}\right):=\mathbf{P}\left(\Cok\begin{pmatrix}A_{2n+1} & \bm{\xi} \\ -\bm{\xi}^T &0 \end{pmatrix}\simeq G^{(2)}\right).$$
\end{enumerate}
\end{defi}

Indeed, by \Cref{prop: same cok and add new row column_alt}, the transition probabilities $\mathbf{P}(G^{(2)},\text{ odd}\mid G^{(1)},\text{ even}),\mathbf{P}(G^{(2)},\text{ even}\mid G^{(1)},\text{ odd})$ do not rely on the matrices we choose, and therefore are well defined.

\begin{prop}\label{cor: globally same pairing_alt}
$P=\{p_1,\ldots,p_l\}$ be a finite set of primes, and $G\in\mathcal{S}_P$. Let $a:=p_1^{e_{p_1}}\cdots p_l^{e_{p_l}}$, where $e_{p_i}=\Dep_{p_i}(G)+1$ for all $1\le i\le l$. 
\begin{enumerate}
\item (Even size) Let $A_{2n}\in\Alt_{2n}(\Z)$, such that $\Cok(A_{2n})_P\simeq G$. Then for all $\Delta A_{2n}\in\Alt_{2n}(a\Z)$, we have $\Cok(A_{2n}+\Delta A_{2n})_P\simeq G$.
\item (Odd size) Let $A_{2n+1}\in\Alt_{2n+1}(\Z)$, such that $\corank(A_{2n+1})=1$, and $\Cok_{\tors}(A_{2n+1})_P\simeq G$. Then for all $\Delta A_{2n+1}\in\Alt_{2n+1}(a\Z)$, we have $\corank(A_{2n+1}+\Delta A_{2n+1})=1$, and $\Cok_{\tors}(A_{2n+1}+\Delta A_{2n+1})_P\simeq G$.
\end{enumerate}
\end{prop}

\begin{proof}
This immediately follows from \Cref{prop: same cok and add new row column_alt} since everything naturally factors over $p_1,\ldots,p_l$.
\end{proof}

\begin{defi}
Let $a\ge 2$, and $n\ge 1$. Let $H_n\in\Alt_n(\Z/a\Z)$ be uniformly distributed. Denote by 
$$\mu_{n}^{\alt,(a)}:=\mathcal{L}(\Cok(H_n))$$
as the law of $\Cok(H_n)$.
\end{defi}

\begin{prop}\label{prop: uniform model exponential convergence_alt}
Let $a\ge 2$ be an integer.
\begin{enumerate}
\item (Even size) When $n$ goes to infinity, $\mu_{2n}^{\alt,(a)}$ weakly converges to a limiting probability measure, denoted by  $\mu_{2\infty}^{\alt,(a)}$. Furthermore, we have
$$D_{L_1}\left(\mu_{2n}^{\alt,(a)},\mu_{2\infty}^{\alt,(a)}\right)=O_a(\exp(-\Omega_a( n))).$$
\item (Odd size) When $n$ goes to infinity, $\mu_{2n+1}^{\alt,(a)}$ weakly converges to a limiting probability measure, denoted by  $\mu_{2\infty+1}^{\alt,(a)}$. Furthermore, we have
$$D_{L_1}\left(\mu_{2n+1}^{\alt,(a)},\mu_{2\infty+1}^{\alt,(a)}\right)=O_a(\exp(-\Omega_a( n))).$$
\end{enumerate}
\end{prop}

\begin{proof}
Suppose $a$ has prime factorization $a=p_1^{e_{p_1}}\cdots p_l^{e_{p_l}}$. In this case, we have
$$\mu_n^{\alt,(a)}=\mu_n^{\alt,(p_1^{e_{p_1}})}\times\cdots\times\mu_n^{\alt,(p_l^{e_{p_l}})}.$$
Notice that for the even size case, the assertions 
$$\mu_{2n}^{\alt,(p_i^{e_{p_i}})}\stackrel{d}{\rightarrow}\mu_{2\infty}^{\alt,(p_i^{e_{p_i}})},D_{L_1}\left(\mu_{2n}^{\alt,(p_i^{e_{p_i}})},\mu_{2\infty}^{\alt,(p_i^{e_{p_i}})}\right)=O_{p_i^{e_{p_i}}}(\exp(-\Omega_{p_i^{e_{p_i}}}( n))),\quad 1\le i\le l,$$
directly implies the conclusion
$$\mu_{2n}^{\alt,(a)}\stackrel{d}{\rightarrow}\mu_{2\infty}^{\alt,(a)}:=\mu_{2\infty}^{\alt,(p_1^{e_{p_1}})}\times\cdots\times\mu_{2\infty}^{\alt,(p_l^{e_{p_l}})},D_{L_1}\left(\mu_{2n}^{\alt,(a)},\mu_{2\infty}^{\alt,(a)}\right)=O_a(\exp(-\Omega_a( n))).$$
Moreover, for the odd size case, the analogous implication also holds. Thus, there is no loss of generality to assume that $a$ only has one prime factor, i.e., $a=p^{e_p}$.

We will first deal with the even size case. We claim that for all integers $0\le k\le n$,
and fixed matrix $B_{2n}\in\Alt_{2n}(\F_p)$ such that $\corank(B_{2n})=2k$, we have 
\begin{equation}\label{eq: distribution conditioned on residue_alt}
\mathcal{L}\left(\Cok(H_{2n})\bigg| H_{2n}/p=B_{2n}\right)=\mathcal{L}\left(\Cok(H_{2k})\bigg| H_{2k}\in\Alt_{2k}(p\Z/p^{e_p}\Z)\right).
\end{equation}
Here, on the right hand side, $H_{2k}$ is also uniformly distributed. To see this, first recall from \Cref{lem: full rank principal minor} that there exists a $(2n-2k)\times(2n-2k)$ principal minor of $B_{2n}$ that is invertible. There is no loss of generality to assume that the $(2n-2k)\times(2n-2k)$ upper-left block of $B_{2n}$ is invertible. In this case, when we uniformly randomly sample a matrix $H_{2n}\in\Alt_{2n}(\Z/p^{e_p}\Z)$ such that $H_{2n}/p=B_{2n}$, we can always use its $(2n-2k)\times(2n-2k)$ upper-left block to eliminate adjacent rows and columns, and the remaining $2k\times 2k$ lower-right block is uniformly distributed in $\Alt_{2k}(p\Z/p^{e_p}\Z)$. Thus, we have confirmed \eqref{eq: distribution conditioned on residue_alt}. As a consequence, when $B_{2n}$ runs through all the matrices in $\Alt_{2n}(\F_p)$ of corank $2k$, we deduce that 
$$\mathcal{L}\left(\Cok(H_{2n})\bigg|\corank(H_{2n}/p)=2k\right)=\mathcal{L}\left(\Cok(H_{2k})\bigg| H_{2k}\in\Alt_{2k}(p\Z/p^{e_p}\Z)\right).$$
Denote $\nu_n^{\alt,p}:=\mathcal{L}(\corank(H_n/p))$. Then, for all $0\le k\le n$ and fixed matrices $H_{2k}'\in\Alt_{2k}(p\Z/p^{e_p}\Z)$, we have
\begin{equation}\label{eq: distribution of quasi of size n_alt}
\mu_{2n}^{\alt,(p^{e_p})}(\Cok(H_{2k}'))=\nu_{2n}^{\alt,p}(2k)\cdot\mathbf{P}\left(\Cok(H_{2k})\simeq\Cok(H_{2k}')\bigg| H_{2k}\in\Alt_{2k}(p\Z/p^{e_p}\Z)\right).
\end{equation}

By \cite[Theorem 5.1]{fulman2015stein}\footnote{In the notation of \cite[Theorem 5.1]{fulman2015stein}, we are taking $q=p$, $\mathcal{Q}_q$ be the probability measure over $\Z_{\ge 0}$ such that $\mathcal{Q}_q(k)=\nu_{2\infty}^{\alt,p}(2k)$, and $\mathcal{Q}_{q,n}$ ($n$ even) be the probability measure over $\{0,1,\ldots,n/2\}$ such that $\mathcal{Q}_{q,n}(k)=\nu_n^{\alt,p}(2k)$. One should also keep in mind the remark in the first paragraph of \cite[Section 5]{fulman2015stein}, which clarifies that although the title refers to symmetric matrices over finite fields with zero diagonal, the alternating case leads to the same formulas.}, we have 
$$D_{L_1}\left(\nu_{2n}^{\alt,p},\nu_{2\infty}^{\alt,p}\right)=O_p(\exp(-\Omega_p(n))),$$
where limit distribution $\nu_{2\infty}^{\alt,p}$ is given by
$$\nu_{2\infty}^{\alt,p}(2k)=
\lim_{n\rightarrow\infty}\nu_{2n}^{\sym,p}(2k)=\prod_{i\ge 0}(1-p^{-2i-1})\frac{p^{2k}}{\prod_{i=1}^{2k}(p^i-1)},\quad\forall k\ge 0.$$
Therefore, when $n$ goes to infinity, $\mu_{2n}^{\alt,(p^{e_p})}$ weakly converges to the limit distribution $\mu_{2\infty}^{\alt,(p^{e_p})}$ given by
\begin{equation}\label{eq: distribution of quasi of size infinity_alt}
\mu_{2\infty}^{\alt,(p^{e_p})}(\Cok(H_{2k}'))=\nu_{2\infty}^{\alt,p}(2k)\cdot\mathbf{P}\left(\Cok(H_{2k})\simeq\Cok(H_{2k}')\bigg| H_{2k}\in\Alt_{2k}(p\Z/p^{e_p}\Z)\right)
\end{equation}
for all $k\ge 0$ and fixed matrices $H_{2k}'\in\Alt_{2k}(p\Z/p^{e_p}\Z)$. Furthermore, comparing \eqref{eq: distribution of quasi of size n_alt} and \eqref{eq: distribution of quasi of size infinity_alt}, we have
\begin{equation}
D_{L_1}\left(\mu_{2n}^{\alt,(p^{e_p})},\mu_{2\infty}^{\alt,(p^{e_p})}\right)=D_{L_1}\left(\nu_{2n}^{\alt,p},\nu_{2\infty}^{\alt,p}\right)=O_p(\exp(-\Omega_p(n))).
\end{equation}

For the odd size case, by \cite[Theorem 5.5]{fulman2015stein}\footnote{In the notation of \cite[Theorem 5.5]{fulman2015stein}, we are taking $q=p$, $\mathcal{Q}_q$ be the probability measure over $\Z_{\ge 0}$ such that $\mathcal{Q}_q(k)=\nu_{2\infty+1}^{\alt,p}(2k+1)$, and $\mathcal{Q}_{q,n}$ ($n$ odd) be the probability measure over $\{0,1,\ldots,(n-1)/2\}$ such that $\mathcal{Q}_{q,n}(k)=\nu_n^{\alt,p}(2k+1)$.}, we have 
$$D_{L_1}\left(\nu_{2n+1}^{\alt,p},\nu_{2\infty+1}^{\alt,p}\right)=O_p(\exp(-\Omega_p(n))),$$
where limit distribution $\nu_{2\infty+1}^{\alt,p}$ is given by
$$\nu_{2\infty+1}^{\alt,p}(2k+1)=
\lim_{n\rightarrow\infty}\nu_{2n+1}^{\alt,p}(2k+1)=\prod_{i\ge 0}(1-p^{-2i-1})\frac{p^{2k+1}}{\prod_{i=1}^{2k+1}(p^i-1)},\quad k\ge 0.$$
Following the same strategy as in the even size case, we have
\begin{equation}
D_{L_1}\left(\mu_{2n+1}^{\alt,(p^{e_p})},\mu_{2\infty+1}^{\alt,(p^{e_p})}\right)=D_{L_1}\left(\nu_{2n+1}^{\alt,p},\nu_{2\infty+1}^{\alt,p}\right)=O_p(\exp(-\Omega_p(n))).
\end{equation}
This completes the proof.
\end{proof}

%% file: The_non_sparse_assumptions.tex
\section{The non-sparse assumptions}\label{sec: nonsparse}

In this section, we introduce the non-sparsity hypotheses that will be used throughout the
Fourier-analytic estimates for the exposure process. As explained in the outline of the paper,
these conditions serve as the technical input that allows us to compare the one-step transition
kernel of the $\epsilon$-balanced process with the corresponding uniform model. 

\begin{defi}\label{defi: non sparse_sym}
(Symmetric case) Let $p$ be a prime number, and $n\ge 1$. We say that a matrix $M_n\in\Sym_n(\Z)$ satisfies
$\mathcal{E}_{n,p}^{\sym}$, if its quotient $M_n/p\in\Sym_n(\F_p)$ has the following properties:
\begin{enumerate}[label=(S\arabic*)]
\item $\corank(M_n/p)\le n^{2/3}$. \label{item: large rank_sym}
\item Suppose $\bm{\xi}\in\F_p^n$ is nonzero, such that $\#I_{\bm{\xi},1}\le n^{3/4}$. Here
$$I_{\bm{\xi},1}:=\{j\in[n]:\frac{n}{1000}<j\le n,\text{ the $j$th row of } M_n/p \text{ is not orthogonal to }\bm{\xi}\}.$$
then $\bm{\xi}$ has at least $\frac{n}{100}$ nonzero entries. \label{item: ortho nonsparse_sym}
\end{enumerate}
\end{defi}

\begin{defi}\label{defi: non sparse_alt}
(Alternating case) Let $p$ be a prime number, and $n\ge 1$. We say that a matrix $A_n\in\Alt_n(\Z)$ satisfies
$\mathcal{E}_{n,p}^{\alt}$, if its quotient $A_n/p\in\Alt_n(\F_p)$ has the following properties:
\begin{enumerate}[label=(A\arabic*)]
\item $\corank(A_n/p)\le n^{2/3}$. \label{item: large rank_alt}
\item Suppose $\bm{\xi}\in\F_p^n$ is nonzero, such that $\#I_{\bm{\xi},1}\le n^{3/4}$. Here
$$I_{\bm{\xi},1}:=\{j\in[n]:1\le j\le n,\text{ the $j$th row of } A_n/p \text{ is not orthogonal to }\bm{\xi}\}.$$
Then $\bm{\xi}$ has at least $\frac{n}{100}$ nonzero entries.
\label{item: ortho nonsparse_alt}
\end{enumerate}
\end{defi}

\begin{defi}\label{defi: revised sym p=2}
(Revised symmetric $p=2$ case) 
Let $n\ge 1$. We say that a matrix $M_n\in\Sym_n(\Z)$ satisfies $\mathcal{E}_n^{\text{\textdagger}}$, if its quotient $M_n/2\in\Sym_n(\F_2)$ has the following properties:
\begin{enumerate}[label=(S\arabic*')]
\item $\corank(M_n/2)\le n^{1/4}$. \label{item: large rank_r2}
\item For all $I_2\subset[n]$ of size $\le n^{1/4}$, one of the following holds:
\begin{enumerate}
\item $(M_n/2)_{I_2^c\times I_2^c}$ has determinant zero. 
\item $(M_n/2)_{I_2^c\times I_2^c}$ is invertible. Denote $\bm{\zeta}_{I_2^c}\in\F_2^{n-\#I_2}$ as the diagonal of $(M_n/2)_{I_2^c\times I_2^c}^{-1}$. Any nonzero linear combination of the columns of 
$(M_n/2)_{I_2^c\times I_2^c}^{-1}(M_n/2)_{I_2^c\times I_2}$ and $\bm{\zeta}_{I_2^c}$ has at least $\frac{1}{2}n^{1/2}$ nonzero entries.
\end{enumerate} 
\label{item: nonsparse combination of inverse diagonal and column}
\end{enumerate}
\end{defi}

Moreover, when $M_n\in\Sym_n(\Z)$ is randomly distributed as in \Cref{thm: exponential convergence_sym}, we can naturally regard $\mathcal{E}_{n,p}^{\sym}$ (here $p$ is a given prime) and $\mathcal{E}_n^{\text{\textdagger}}$ as random events. Similarly, when $A_n\in\Alt_n(\Z)$ is randomly distributed as in \Cref{thm: exponential convergence_alt}, we can naturally regard $\mathcal{E}_{n,p}^{\alt}$ as a random event.

\begin{rmk}\label{rem:S2-vs-A2}
The non-sparsity events in \Cref{defi: non sparse_sym} and \Cref{defi: non sparse_alt} are stated in slightly different forms.
In the symmetric case, our Fourier-analytic estimates involve quadratic terms (coming
from the diagonal and the symmetry constraint), and we control them via a decoupling step that
reduces quadratic expressions to bilinear ones, see the first case of the proof of \Cref{prop: universal transition_sym}. For this technical reason, in \ref{item: ortho nonsparse_sym} we only impose the relevant non-sparsity requirement on
a tail window of indices (e.g.\ $j\in(n/1000,n]$).

By contrast, in the alternating case the diagonal vanishes identically and the corresponding Fourier
analysis does not produce quadratic terms. No decoupling is needed, and the argument works with the
more straightforward non-sparsity condition \ref{item: ortho nonsparse_alt} imposed across all indices.
\end{rmk}

The following theorems are the main results of this section. Roughly speaking, except with exponentially small probability, we can assume the property $\mathcal{E}_{n,p}^{\sym}$ (resp. $\mathcal{E}_{n,p}^{\alt}$) for $\epsilon$-balanced symmetric (resp. alternating) matrices.

\begin{thm}\label{thm: not sparse}
(Symmetric case and alternating case) Let $M_n\in\Sym_n(\Z)$ be randomly distributed as in \Cref{thm: exponential convergence_sym}, and let $A_n\in\Alt_n(\Z)$ be randomly distributed as in \Cref{thm: exponential convergence_alt}. Then we have 
$$\mathbf{P}(M_n\text{ satisfies }\mathcal{E}_{n,p}^{\sym})=1-O_{p,\epsilon}(\exp(-\epsilon\Omega_{p}(n))),$$
and
$$\mathbf{P}(A_n\text{ satisfies }\mathcal{E}_{n,p}^{\alt})=1-O_{p,\epsilon}(\exp(-\epsilon\Omega_{p}(n))).$$
\end{thm}

\begin{thm}\label{thm: not sparse revised p=2}
(Revised symmetric $p=2$ case) Let $M_n\in\Sym_n(\Z)$ be randomly distributed as in \Cref{thm: exponential convergence_sym}. Then we have 
$$\mathbf{P}(M_n\text{ satisfies }\mathcal{E}_n^{\text{\textdagger}})=1-O_{\epsilon}(\exp(-\Omega_\epsilon(n^{1/2}))).$$
\end{thm}

\begin{lemma}\label{lem: small corank}
Let $M_n\in\Sym_n(\Z)$ be randomly distributed as in \Cref{thm: exponential convergence_sym}, and let $A_n\in\Alt_n(\Z)$ be randomly distributed as in \Cref{thm: exponential convergence_alt}. Let $p$ be a prime number, and $c>0$ be an absolute constant. Then we have
$$\mathbf{P}(\corank(M_n/p)\ge n^c)=O_{\epsilon,c}(\exp(-\Omega_{\epsilon,c}(n^{2c}))),$$
and
$$\mathbf{P}(\corank(A_n/p)\ge n^c)=O_{\epsilon,c}(\exp(-\Omega_{\epsilon,c}(n^{2c}))).$$
\end{lemma}

\begin{proof}
We will only prove the symmetric case, and the alternating case is similar. Suppose $\mathbf{P}(\corank(M_n/p)\ge n^c)$. By \Cref{lem: full rank principal minor}, there exist $n^c\le j\le n$ and a subset $I_p\subset [n]$ of size $j$, such that $(M_n/p)_{J_p^c}$ is invertible, and
$$(M_n/p)_{I_p\times I_p}=(M_n/p)_{I_p\times I_p^c}(M_n/p)_{I_p^c\times I_p^c}(M_n/p)_{I_p^c\times I_p}.$$
For each $I_p$, this gives probability $\le (1-\epsilon)^{\frac{1}{2}\#I_p(I_p+1)}\le \exp(-\frac{1}{2}\epsilon(\#I_p)^2)$. Therefore, when $n$ is sufficiently large (depending on $\epsilon,c$), we have
\begin{align}
\begin{split}
\mathbf{P}\left(\corank(M_n/p)\ge n^c\right)&\le\sum_{j=n^c}^{n}\sum_{\#I_p=j}\exp\left(-\frac{1}{2}\epsilon(\#I_p)^2\right)\\
&\le\sum_{j=n^c}^{n} n^j\exp\left(-\frac{1}{2}\epsilon j^2\right)\\
&\le n\cdot \exp\left(-\frac{1}{2}\epsilon n^{2c}-n^c\log n\right)\\
&=O_{\epsilon,c}(\exp(-\Omega_{\epsilon,c}(n^{2c}))).
\end{split}
\end{align}
Here, for the second line, we use the fact that the number of subsets of $[n]$ of size $j$ is $\binom{n}{j}\le n^j$, and for the third line, we replace $j$ in all summation terms with $n^c$ due to monotonicity.
\end{proof}

\subsection{Proof of the symmetric case and alternating case}

\begin{lemma}
Let $\epsilon>0$ be a real number, and  $M_n\in\Sym_n(\Z),A_n\in\Alt_n(\Z)$ be the same as in \Cref{thm: exponential convergence_sym} and \Cref{thm: exponential convergence_alt}. Then, the probability that the requirement \ref{item: ortho nonsparse_sym} (resp. \ref{item: ortho nonsparse_alt}) is not satisfied is $O_{p,\epsilon}(\exp(-\epsilon\Omega_{p}(n)))$.
\end{lemma}

\begin{proof}
We will only prove the symmetric case, and the alternating case is similar. Suppose the requirement \ref{item: ortho nonsparse_sym} does not hold. We will discuss two possible cases:
\begin{enumerate}
\item We have a $\bm{\xi}\in\F_p^n$ with $1\le \wt(\bm{\xi})\le n^{3/4}$, such that $\#I_{\bm{\xi},1}\le n^{3/4}$. For each choice of such $\bm{\xi}$ and $I_{\bm{\xi},1}$, we can find a set $I_{\bm{\xi},0}\ge\frac{n}{2}$ that does not intersect with $[\frac{n}{1000}]\bigcup I_{\bm{\xi},1}\bigcup\supp(\bm{\xi})$. Notice that the entries with row index in $I_{\bm{\xi},0}$ and column index in $\supp(\bm{\xi})$ are independent. The rows with index in $I_{\bm{\xi},0}$ have to be orthogonal to $\bm{\xi}$, each of them having probability $\le 1-\epsilon$. Thus, each choice of $\bm{\xi}$ and $I_{\bm{\xi},1}$ contributes probability $\le (1-\epsilon)^{n/2}\le \exp(-\epsilon n/2)$. On the other hand, there are at most $ n^{3/4}\cdot\binom{n}{n^{3/4}}\le n^{n^{3/4}+1}$ choices for such $\bm{\xi}$, and at most $\le n^{n^{3/4}+1}$ choices for such $I_{\bm{\xi},1}$. In conclusion, this case contributes probability $\le n^{2n^{3/4}+2}\exp(-\epsilon n/2)=O_{\epsilon}\exp(-\epsilon\Omega(n)).$
\item We have a $\bm{\xi}\in\F_p^n$ with $n^{3/4}<\wt(\bm{\xi})\le \frac{n}{100}$, such that $\#I_{\bm{\xi},1}\le n^{3/4}$. Similar to the above case, for each choice of such $\bm{\xi}$ and $I_{\bm{\xi},1}$, we can find a set $I_{\bm{\xi},0}\ge\frac{n}{2}$ that does not intersect with $[\frac{n}{1000}]\bigcup I_{\bm{\xi},1}\bigcup\supp(\bm{\xi})$. In this case, we can take $n$ sufficiently large (depending on $\epsilon$), such that for all $r\in\F_p$ and independent $\epsilon$-stable random integers $\xi_1,\ldots,\xi_{n^{3/4}}$,
$$\mathbf{P}(\xi_1+\cdots+\xi_{n^{3/4}}\equiv r \mod p)\le \frac{1}{p-0.001}.$$
In this case, each choice of $\bm{\xi}$ and $I_{\bm{\xi},1}$ contributes probability $\le\frac{1}{(p-0.1)^{n/2}}$. On the other hand, by \Cref{prop: asymptotic of Hamming ball}, there are at most 
$$\Vol_p(n,n/100)\le p^{H_p(0.01)n}\le p^{n/10}$$
choices for such $\bm{\xi}$, and at most $n^{n^{3/4}+1}$ choices for such $I_{\bm{\xi},1}$. In conclusion, when $n$ is sufficiently large (depending on $\epsilon$), this case contributes probability $\le\frac{p^{n/10}n^{n^{3/4}+1}}{(p-0.001)^{n/2}}=O_{p,\epsilon}\exp(-\Omega_p(n))$.
\end{enumerate}
The above two cases together complete the proof.
\end{proof}

\begin{proof}[Proof of \Cref{thm: not sparse}]
For the symmetric case, each of the requirements \ref{item: large rank_sym}, \ref{item: ortho nonsparse_sym} holds with probability at least $1-O_{p,\epsilon}\exp(-\epsilon\Omega_p(n))$. Therefore, the probability that both requirements hold is also at least $1-O_{p,\epsilon}\exp(-\epsilon\Omega_p(n))$. The alternating case is analogous. 
\end{proof}

\subsection{Proof of the revised symmetric $p=2$ case}

\begin{lemma}\label{lem: inverse nonzero diagonal}
Suppose $M_n/2\in\Sym_n(\F_2)$ is invertible. Then the following are equivalent:
\begin{enumerate}
\item The entries on the diagonal of $M_n/2$ are all zero.
\item The entries on the diagonal of $(M_n/2)^{-1}$ are all zero.
\end{enumerate}
\end{lemma}

\begin{proof}
We only need to prove that the first assertion implies the second. Suppose that the entries on the diagonal of $M_n/2$ are all zero. In this case, $M_n/2$ can also be regarded as an alternating matrix. Since $M_n/2$ is invertible, $n$ must be even. The entries on the diagonal of $(M_n/2)^{-1}$ are exactly the $(n-1)\times(n-1)$ principal minors of $M_n/2$, which must be zero because they are of odd size.
\end{proof}

The following lemma provides an asymptotic for the transition probability of corank when adding a new random row and column. Although we state the lemma for a general prime number $p$, in this subsection we only use the case $p=2$. The proof can be deduced in the same way as in Step 1 of \cite[Lemma 3.2]{ferber2023random} or the unpublished work of Maples \cite[Proposition 3.1]{Maples_symma_2013}, so we omit it here.

\begin{lemma}\label{lem: lower bound of corank decrease probability}

Let $\epsilon>0$ be a real number, and $p$ be a prime number. Suppose $M_j\in\Sym_j(\Z)$ is a fixed matrix that satisfies $\mathcal{E}_{j,p}^{\sym}$. Let $M_{j+1}\in\Sym_{j+1}(\Z)$ be the random matrix adding a new independent $\epsilon$-balanced row and column to $M_j$. Then we have 
$$\corank(M_{j+1}/p)-\corank(M_j/p)\in\{-1,0,1\}.$$ 
Moreover, we have
$$\mathbf{P}(\corank(M_{j+1}/p)-\corank(M_j/p)=-1)\ge 1-\frac{1}{p^2}-O_{p,\epsilon}(\exp(-\Omega_{p,\epsilon}(n)))$$
when $\corank(M_j/p)\ge 2$, and
$$\mathbf{P}(\corank(M_{j+1}/p)-\corank(M_j/p)=-1)\ge 1-\frac{1}{p}-O_{p,\epsilon}(\exp(-\Omega_{p,\epsilon}(n)))$$
when $\corank(M_j/p)=1$.
\end{lemma}

\begin{rmk}
In fact, when $p$ is odd, \cite[Lemma 3.2]{ferber2023random} also provides asymptotics of the transition probabilities that $\corank(M_{j+1}/p)-\corank(M_j/p)=0$ and $\corank(M_{j+1}/p)-\corank(M_j/p)=1$, which all match the uniform model. However, these asymptotics rely on Fourier analysis estimate over the quadratic form, which requires $p\ne 2$. Nevertheless, the transition probability that $\corank(M_{j+1}/2)-\corank(M_j/2)=-1$ only depends on linear form, so the case $p=2$ still applies.
\end{rmk}

\begin{lemma}\label{lem: det zero or invertible at the end}
Let $M_n$ be the same as in \Cref{thm: exponential convergence_sym}. Let $J\subset[n]$ such that $\# J\le n^{1/2}$. With probability no less than $1-O_{\epsilon}(\exp(-\Omega_\epsilon(n^{1/2})))$, at least one of the following holds:
\begin{enumerate}
\item $M_n/2$ has determinant zero. \label{item: determinant zero global}
\item There exists $j\in J$, such that the $j$th principal minor (i.e., the $(n-1)\times(n-1)$ principal minor obtained by removing the $j$th row and column) of $M_n/2$ is invertible. \label{item: principal minor at the end}
\end{enumerate}
\end{lemma}

\begin{proof}
It suffices to prove the case $J=\{j\in[n]:n-n^{1/2}\le j\le n\}$. By \Cref{lem: small corank}, with probability at least $1-O_{\epsilon}(\exp(-\Omega_{\epsilon}(n^{1/2})))$, $(M_n/2)_{n-n^{1/2}}$ has rank $\le n^{1/3}$. By \Cref{thm: not sparse}, with probability at least $1-O_{\epsilon}(\exp(-\epsilon\Omega(n))$, the events $$\mathcal{E}_{j,2}^{\sym},n-n^{1/2}\le j<n-\frac{1}{2}n^{1/2}$$
all happen. Therefore, by a straightforward comparison argument using \Cref{lem: lower bound of corank decrease probability}, we can apply a standard Hoeffding
large-deviation bound to deduce that with probability at least $1-O_{\epsilon}(\exp(-\Omega_\epsilon(n^{1/2})))$, there exists $n-n^{1/2}\le j_0<n-\frac{1}{2}n^{1/2}$ such that $M_{j_0}/2$ is invertible. When this happens, we use the $j_0\times j_0$ upper-left block to eliminate adjacent rows and columns via elementary transform, and the lower-right block will be transformed to
$$M':=(M_n/2)_{[j_0]^c\times[j_0]^c}-(M_n/2)_{[j_0]^c\times [j_0]}\times(M_n/2)_{[j_0]\times [j_0]}^{-1}\times(M_n/2)_{[j_0]\times [j_0]^c}\in\Sym_{n-j_0}(\F_2).$$
With probability at least $1-(1-\epsilon)^{n-j_0}=1-O_{\epsilon}(\exp(-\Omega_\epsilon(n^{1/2})))$, the matrix $M'$ has a $1$ on the diagonal. 

In summary, with probability at least $1-O_{\epsilon}(\exp(-\Omega_\epsilon(n^{1/2})))$, there exists $n-n^{1/2}\le j_0<n-\frac{1}{2}n^{1/2}$ such that $j_0\times j_0$ upper-left block of $M_n/2$ is invertible, and after we use this block to eliminate adjacent rows and columns via elementary transform, with probability at least $1-(1-\epsilon)^{n-j_0}=1-O_{\epsilon}(\exp(-\Omega_\epsilon(n^{1/2})))$, the remaining lower-right block $M'\in\Sym_{n-j_0}(\F_2)$ has a $1$ on the diagonal. Now, if $M'$ has determinant zero, then $M_n/2$ also has determinant zero. Otherwise, if $M'$ is invertible, by \Cref{lem: inverse nonzero diagonal} there exists $j\in[j_0+1,n]$, such that the $(j-j_0)$-th  principal minor of $M'$ is invertible. As a consequence, the $j$-th  principal minor of $M_n/2$ is invertible. This completes the proof.
\end{proof}

\begin{cor}\label{cor: non-sparse diagonal of inverse}
With probability no less than $1-O_\epsilon(\exp(-\Omega_\epsilon(n^{1/2})))$, at least one of the following holds:
\begin{enumerate}
\item $M_n/2$ has determinant zero. 
\item Among all the $n$ principal minors of size $n-1$, at least $n^{1/2}$ of them are invertible.
\end{enumerate}
\end{cor}

\begin{proof}
We partition the set $[n]$ into $n^{1/2}$ intervals, each of them having size $n^{1/2}$. Applying \Cref{lem: det zero or invertible at the end} to each interval, we deduce that with probability $1-n^{1/2}O_\epsilon(\exp(-\Omega_\epsilon(n^{1/2})))=1-O_\epsilon(\exp(-\Omega_\epsilon(n^{1/2})))$, at least one of the two properties appearing in the statement of the Corollary is satisfied.
\end{proof}

\begin{lemma}\label{lem: non-sparse diagonal of inverse and new columns}
Let $M_n/2\in\Sym_n(\F_2)$ be fixed and invertible, with at least $n^{1/2}$ invertible $(n-1)\times (n-1)$ principal minors. Denote by $\bm{\zeta}\in\F_2^n$ the diagonal of $(M_n/2)^{-1}$, so that $\wt(\bm{\zeta})\ge n^{1/2}$. Let
$$\bm{\eta}_1=(\eta_{1,1},\ldots,\eta_{1,n}),\ldots,\bm{\eta}_{n^{2/3}}=(\eta_{n^{2/3},1},\ldots,\eta_{n^{2/3},n})\in\F_2^n$$
be random vectors. Here, the elements 
$\{\eta_{i,j}:1\le i\le n^{2/3},1\le j\le n\}$
are independent elements in $\F_2$ that are not $\epsilon$-concentrated. Then, with probability at least $1-O_\epsilon(\exp(-\Omega_\epsilon(n)))$, every nonzero linear combination of $(M_n/2)^{-1}\bm{\eta}_1,\ldots,(M_n/2)^{-1}\bm{\eta}_{n^{2/3}},\bm{\zeta}$ has weight $\ge n^{1/2}$.
\end{lemma}

\begin{proof}
On the one hand, denote
$$\mathcal{S}:=\{\bm{\eta}\in\F_2^n:\wt((M_n/2)^{-1}\bm{\eta})< n^{1/2},\wt((M_n/2)^{-1}\bm{\eta}+\bm{\zeta})<n^{1/2}\}.$$
we have
$$\#\mathcal{S}\le 2n^{1/2}\binom{n}{n^{1/2}}\le 2n^{n^{1/2}}.$$
On the other hand, there are only $2^{n^{2/3}}-1$ nonzero linear combinations of $\bm{\eta}_1,\ldots,\bm{\eta}_{n^{2/3}}$ in total. Each of them gives a random vector whose entries are independent and not $\epsilon$-concentrated. Therefore, the probability that there exists a nonzero linear combination of $$(M_n/2)^{-1}\bm{\eta}_1,\ldots,(M_n/2)^{-1}\bm{\eta}_{n^{2/3}},\bm{\zeta}$$ 
with weight $\ge n^{1/2}$ is no greater than
$$(2^{n^{2/3}}-1)\cdot(1-\epsilon)^n\#\mathcal{S}\le(2^{n^{2/3}}-1)\cdot2n^{n^{1/2}}\cdot(1-\epsilon)^n=O_\epsilon(\exp(-\Omega_\epsilon(n))).$$
\end{proof}

\begin{proof}[Proof of \Cref{thm: not sparse revised p=2}]
On the one hand, by \Cref{lem: small corank}, with probability at least $1-O_\epsilon(\exp(-\Omega_\epsilon(n^{1/2})))$, the property \ref{item: large rank_r2} holds. 

On the other hand, let us take an arbitrary subset $I_2\subset[n]$ with $\#I_2\le n^{1/4}$. Combining \Cref{cor: non-sparse diagonal of inverse} and \Cref{lem: non-sparse diagonal of inverse and new columns}, we deduce that when $n$ is sufficiently large (depending on $\epsilon$), with probability probability at least $1-O_\epsilon(\exp(-\Omega_\epsilon(n^{1/2})))$, one of the following holds:
\begin{enumerate}
\item $(M_n/2)_{I_2^c\times I_2^c}$ has determinant zero. 
\item $(M_n/2)_{I_2^c\times I_2^c}$ is invertible. Denote $\bm{\zeta}_{I_2^c}\in\F_2^{n-\#I_2}$ as the diagonal of $(M_n/2)_{I_2^c\times I_2^c}^{-1}$. Any nonzero linear combination of the columns of 
$(M_n/2)_{I_2^c\times I_2^c}^{-1}(M_n/2)_{I_2^c\times I_2}$ and $\bm{\zeta}_{I_2^c}$ has at least 
$$(n-\#I_2^c)^{1/2}\ge(n-n^{1/4})^{1/2}\ge\frac{1}{2}n^{1/2}$$ 
nonzero entries.
\end{enumerate} 
Furthermore, there are only at most $n^{1/4}\cdot\binom{n}{n^{1/4}}\le n^{n^{1/4}}$ subsets $I_2\subset[n]$ of size $\le n^{1/4}$. Therefore, with probability at least $1-n^{n^{1/4}}\cdot O_\epsilon(\exp(-\Omega_\epsilon(n^{1/2})))=1-O_\epsilon(\exp(-\Omega_\epsilon(n^{1/2})))$, the property \eqref{item: nonsparse combination of inverse diagonal and column} holds. This completes the proof.
\end{proof}

%% file: Proof_of_the_symmetric_case.tex
\section{Proof of the symmetric case}\label{sec: Proof of the symmetric case}

\begin{thm}\label{thm: quotient simulate_sym}
Let $\epsilon>0$ and $M_n\in\Sym_n(\Z)$ be the same as in \Cref{thm: exponential convergence_sym}, and $a\ge 2$ be an integer. Then we have
$$D_{L^1}\left(\mathcal{L}(\Cok_*^{(a)}(M_n/a)),\mu_\infty^{\sym,(a)}\right)=O_{a,\epsilon}(\exp(-\Omega_{a,\epsilon}(n^{1/2}))).$$
Here, $\mu_\infty^{\sym,(a)}$ is the same as in \Cref{prop: uniform model exponential convergence_sym}. Furthermore, if $a$ is odd, the right hand side of the above can be reduced to $O_{a,\epsilon}(\exp(-\epsilon\Omega_a(n)
)).$ 
\end{thm}

Before we return to the proof of \Cref{thm: quotient simulate_sym}, we see how it implies \Cref{thm: exponential convergence_sym}.

\begin{proof}[Proof of \Cref{thm: exponential convergence_sym}, assuming \Cref{thm: quotient simulate_sym}]
Let $a:=p_1^{e_{p_1}}\cdots p_l^{e_{p_l}}$ be the same as in \Cref{cor: globally same pairing}. In this case, for all $1\le i\le l$ and $n\ge 1$,
$$\mu_n^{\sym,(p_i^{e_{p_i}})}\left(\left(G_{p_i},\langle\cdot,\cdot\rangle^{(p_i^{e_{p_i}})}_*\right)\right)=\mathbf{P}\left(\left(\Cok(M_n^{(p_i')})_{p_i},\langle\cdot,\cdot\rangle_{p_i}\right)\simeq\left(G_{p_i},\langle\cdot,\cdot\rangle_{G_{p_i}}\right)\right).$$
Here, $M_n^{(p_i')}\in\Sym_n(\Z_{p_i})$ is Haar-distributed, and the definition of $(G_{p_i},\langle\cdot,\cdot\rangle_*^{(p_i^{e_{p_i}})})$ follows from \Cref{defi: quotient version of paired group}. Applying \cite[Theorem 2]{clancy2015cohen}, we have
\begin{align}\label{eq: uniform model in S_sym}
\begin{split}
\mu_\infty^{\sym,(a)}\left(\left(G,\langle\cdot,\cdot\rangle^{(a)}_*\right)\right)&=\prod_{i=1}^l\mu_\infty^{\sym,(p_i^{e_{p_i}})}\left(\left(G_{p_i},\langle\cdot,\cdot\rangle^{(p_i^{e_{p_i}})}_*\right)\right)\\
&=\prod_{i=1}^l\lim_{n\rightarrow\infty}\mathbf{P}\left(\left(\Cok(M_n^{(p_i')})_{p_i},\langle\cdot,\cdot\rangle_{p_i}\right)\simeq\left(G_{p_i},\langle\cdot,\cdot\rangle_{G_{p_i}}\right)\right)\\
&=\prod_{i=1}^l\frac{\prod_{k\ge 1}(1-p_i^{1-2k})}{|G_{p_i}||\Aut(G_{p_i},\langle\cdot,\cdot\rangle_{G_{p_i}})|}\\
&=\frac{\prod_{i=1}^l\prod_{k\ge 1}(1-p_i^{1-2k})}{|G||\Aut(G,\langle\cdot,\cdot\rangle_G)|}.
\end{split}
\end{align}
Therefore, due to our assumption in \Cref{thm: quotient simulate_sym}, we have
$$\left|\mathbf{P}\left((\Cok(M_n),\langle\cdot,\cdot\rangle_P)\simeq(G,\langle\cdot,\cdot\rangle_G)\right)-\frac{\prod_{i=1}^l\prod_{k\ge 1}(1-p_i^{1-2k})}{|G||\Aut(G,\langle\cdot,\cdot\rangle_G)|}\right|=O_{a,\epsilon}(\exp(-\Omega_{a,\epsilon}(n^{1/2}))).$$
Moreover, the error term on the right hand side of the above can be reduced to $O_{a,\epsilon}(\exp(-\epsilon\Omega_a(n)))$ when $a$ is odd. Notice that the integer $a$ defined in \Cref{cor: globally same pairing} depends on $G,P$, and in particular, $a$ is odd if and only if $2\notin P$. Therefore, we can rewrite $O_{a,\epsilon}(\exp(-\Omega_{a,\epsilon}(n^{1/2})))$ as $O_{G,P,\epsilon}(\exp(-\Omega_{G,P,\epsilon}(n^{1/2})))$ when $2\in P$, and rewrite $O_{a,\epsilon}(\exp(-\epsilon\Omega_a(n)))$ as $O_{G,P,\epsilon}(\exp(-\epsilon\Omega_{G,P}(n)
))$ when $2\notin P$. This completes the proof.
\end{proof}


Now we return to our proof of \Cref{thm: quotient simulate_sym}. The following proposition provides a more refined characterization based on \cite[Lemma 3.2]{ferber2023random}.

\begin{prop}\label{prop: universal transition_sym}
Let $\epsilon>0$ be a real number, and $a\ge 2$ be an integer with prime factorization $a=p_1^{e_{p_1}}\cdots p_l^{e_{p_l}}$. Suppose $M_n\in\Sym_n(\Z)$ is a fixed matrix that satisfies $\mathcal{E}_{n,p_i}^{\sym}$ for $1\le i\le l$. Furthermore, when $a$ is even, we assume that $M_n$ satisfies $\mathcal{E}_n^{\text{\textdagger}}$. Let $z,\xi_1,\ldots,\xi_n$ be independent $\epsilon$-balanced random integers, and $z',\xi_1',\ldots,\xi_n'$ be independent and uniformly distributed in $\{0,1,\ldots,a-1\}$. Let $\bm{\xi}=(\xi_1,\ldots,\xi_n),\bm{\xi}'=(\xi_1',\ldots,\xi_n')\in\Z^n$. Then we have
$$
D_{L^2}\left(\Cok_*^{(a)}
\begin{pmatrix}M_n/a & \bm{\xi}/a\\
(\bm{\xi}/a)^T & z/a
\end{pmatrix}
,\Cok_*^{(a)}\begin{pmatrix}M_n/a & \bm{\xi}'/a\\
(\bm{\xi}'/a)^T & z'/a
\end{pmatrix}\right)=O_{a,\epsilon}(\exp(-\Omega_{a,\epsilon}(n^{1/2}))).\\
$$
Furthermore, if $a$ is odd, the right hand side of the above can be reduced to $O_{a,\epsilon}(\exp(-\epsilon\Omega_a(n))
).$
\end{prop}

\begin{proof}
For all $1\le i\le l$, by \Cref{lem: full rank principal minor}, there exists $I_{p_i}\subset[n]$ with $\# I_{p_i}=n-\rank(M_n/{p_i})$, such that $(M_n/p_i)_{I_{p_i}^c\times I_{p_i}^c}$ has full rank. By \ref{item: large rank_sym}, we have $\#I_{p_i}\le n^{2/3}$. 

Now, we treat the $p
_i$-part of the cokernel with its quasi-pairing, where we reduce all matrix entries modulo $p_i^{e_{p_i}}$. By paired row-column operations, we can use the invertible matrix $(M_n/p_i^{e_{p_i}})_{I_{p_i}^c\times I_{p_i}^c}$ to eliminate the components of $\bm{\xi}/p_i^{e_{p_i}}$ and $(\bm{\xi}/p_i^{e_{p_i}})^T$ indexed from $I_{p_i}^c$, while the symmetric structure is preserved. In this way, the new added column $(\bm{\xi}/p_i^{e_{p_i}},z/p_i^{e_{p_i}})$ will be transformed into the form
\begin{multline}\bm{w}_i:=
\Bigg((\bm{\xi}/p_i^{e_{p_i}})_{I_{p_i}}-(M_n/p_i^{e_{p_i}})_{I_{p_i}\times I_{p_i}^c}(M_n/p_i^{e_{p_i}})_{I_{p_i}^c\times I_{p_i}^c}^{-1}(\bm{\xi}/p_i^{e_{p_i}})_{I_{p_i}^c},\\
z-(\bm{\xi}/p_i^{e_{p_i}})_{I_{p_i}^c}^T(M_n/p_i^{e_{p_i}})_{I_{p_i}^c\times I_{p_i}^c}^{-1}(\bm{\xi}/p_i^{e_{p_i}})_{I_{p_i}^c}\Bigg)\in(\Z/p_i^{e_{p_i}}\Z)^{\#I_{p_i}+1}.
\end{multline}
Here, the first line is a linear term whose index is in $I_{p_i}$, while the second line is the quadratic term on the lower-right corner. Moreover, we ignore and skip the components with index $I_{p_i}^c$ because they have been eliminated to zero. When $i$ runs over all integers from $1$ to $l$, we obtain a combination of vectors $(\bm{w}_i)_{1\le i\le l}$,
and this combination is a random variable in $\left((\Z/p_i^{e_{p_i}}\Z)^{\#I_{p_i}+1}\right)_{1\le i\le l}$. If we replace $z,\xi_1,\ldots,\xi_n$ by the input $z',\xi_1',\ldots,\xi_n'$, we can obtain a combination of vectors $(\bm{w}_i')_{1\le i\le l}\in\left((\Z/p_i^{e_{p_i}}\Z)^{\#I_{p_i}+1}\right)_{1\le i\le l}$ in the similar way. It is clear that $(\bm{w}_i')_{1\le i\le l}$ is uniformly distributed. We only need to prove that $$D_{L^2}\left((\bm{w}_i)_{1\le i\le l},(\bm{w}_i')_{1\le i\le l}\right)=\begin{cases}
O_{a,\epsilon}(\exp(-\Omega_{a,\epsilon}(n^{1/2}))) & a\text{ even}\\
O_{a,\epsilon}(\exp(-\epsilon\Omega_a(n))) & a\text{ odd}
\end{cases}.$$
To achieve this, applying Parseval's identity, we have
$$D_{L^2}^2\left((\bm{w}_i)_{1\le i\le l},(\bm{w}_i')_{1\le i\le l}\right)=\prod_{i=1}^lp_i^{-e_{p_i}(\#I_{p_i}+1)}\sum\left|\E_{\bm{\xi},z}\exp\left(\sum_{i=1}^l\frac{2\pi \sqrt{-1}}{p_i^{e_{p_i}}}(\bm{\alpha}_i,b_i)\cdot\bm{w}_i\right)\right|^2.$$
Here, the sum on the right hand side ranges through the nonzero elements in the set
$$
\mathcal{V}:=\{(\bm{\alpha}_i,b_i)_{1\le i\le l}: \bm{\alpha}_i:=(\alpha_{i,j_i})_{j_i\in I_{p_i}}\in(\Z/p_i^{e_{p_i}}\Z)^{\#I_{p_i}},b_i\in\Z/p_i^{e_{p_i}}\Z\}.
$$
Therefore, it suffices to show that for all nonzero $(\bm{\alpha}_i,b_i)_{1\le i\le l}\in\mathcal{V}$, we have
\begin{equation}\label{eq: small at nonzero vector_sym}
\left|\E_{\bm{\xi},z}\exp\left(\sum_{i=1}^l\frac{2\pi \sqrt{-1}}{p_i^{e_{p_i}}}(\bm{\alpha}_i,b_i)\cdot\bm{w}_i\right)\right|=\begin{cases}
O_{a,\epsilon}(\exp(-\Omega_{a,\epsilon}(n^{1/2}))) & a\text{ even}\\
O_{a,\epsilon}(\exp(-\epsilon\Omega_a(n))) & a\text{ odd}
\end{cases}.
\end{equation}

The proof of \eqref{eq: small at nonzero vector_sym} proceeds by a case-by-case discussion as follows.

\noindent{\bf Case 1. When $a$ is odd, and $b_i\ne 0$ for some $1\le i\le l$. } Without loss of generality, we can let $b_1\ne 0$. We will use a generalized version of the decoupling trick introduced in the proof of \cite[Lemma 3.2]{ferber2023random}. Since $\#I_{p_i}\le n^{2/3}$ for all $1\le i\le l$, we have $\#\bigcup_{i=1}^l I_{p_i}\le ln^{2/3}$. We denote
$$J=[\frac{n}{1000}]\backslash\bigcup_{i=1}^lI_{p_i}\subset[n],$$
so that $\#J\in[\frac{n}{1000}-ln^{2/3},\frac{n}{1000}]$. Let $I:=J^c$, so that $I_{p_i}\subset I$ for all $1\le i\le l$. Let $\bm{\xi}_I,\bm{\xi}_J$ be the obvious restrictions. Let $\bm{\xi}^{\res}_J$ be an independent resample of $\bm{\xi}_J$. Let $\bm{\xi}^{\res}:=\bm{\xi}_I+\bm{\xi}_J^{\res}$. Let $\bm{\phi}_i:=\bm{\xi}_{I_{p_i}^c},\bm{\phi}^{\res}_i:=\bm{\xi}^{\res}_{I_{p_i}^c}$. We have 
\begin{align}\label{eq: decoupling_sym}
\begin{split}
&\left|\E_{\bm{\xi},z}\exp\left(\sum_{i=1}^l\frac{2\pi \sqrt{-1}}{p_i^{e_{p_i}}}(\bm{\alpha}_i,b_i)\cdot\bm{w}_i\right)\right|^2\\
&=\E_{\bm{\xi}_I,\bm{\xi}_J,\bm{\xi}^{\res}_J}\exp\Bigg(\sum_{i=1}^l\frac{2\pi \sqrt{-1}}{p_i^{e_{p_i}}}\Bigg(b_i(\bm{\phi}_i^T(M_n/p_i^{e_{p_i}})_{I_{p_i}^c\times I_{p_i}^c}^{-1}\bm{\phi}_i-(\bm{\phi}^{\res}_i)^T(M_n/p_i^{e_{p_i}})_{I_{p_i}^c\times I_{p_i}^c}^{-1}\bm{\phi}^{\res}_i)\\
&+\bm{\alpha}_i\cdot((M_n/p_i^{e_{p_i}})_{I_{p_i}\times I_{p_i}^c}(M_n/p_i^{e_{p_i}})_{I_{p_i}^c\times I_{p_i}^c}^{-1}(\bm{\xi}_J-\bm{\xi}^{\res}_J))\Bigg)\Bigg)\\
&\le\E_{\bm{\xi}_J,\bm{\xi}_J^{\res}}\left|\E_{\bm{\xi}_I}\left[\exp\left(\sum_{i=1}^l\frac{2\pi \sqrt{-1}}{p_i^{e_{p_i}}}b_i(\bm{\phi}_i-\bm{\phi}^{\res}_i)^T(M_n/p_i^{e_{p_i}})_{I_{p_i}^c\times I_{p_i}^c}^{-1}(\bm{\phi}_i+\bm{\phi}^{\res}_i)\right)\Bigg|\bm{\xi}_J,\bm{\xi}^{\res}_J\right]\right|\\
&=\E_{\bm{\xi}_J,\bm{\xi}^{\res}_J}\left|\E_{\bm{\xi}_I}\left[\exp\left(\sum_{i=1}^l\frac{4\pi \sqrt{-1}}{p_i^{e_{p_i}}}b_i(\bm{\xi}_J-\bm{\xi}^{\res}_J)^T(M_n/p_i^{e_{p_i}})_{I_{p_i}^c\times I_{p_i}^c}^{-1}\bm{\phi}_i''\right)\Bigg|\bm{\xi}_J,\bm{\xi}^{\res}_J\right]\right|.\\
\end{split}
\end{align}
Here, for each $1\le i\le l$, $\bm{\xi}_J-\bm{\xi}^{\res}_J$ and $\bm{\phi}_i''$ are vectors with entry index in $I_{p_i}^c$. We have abused the notation by using $\bm{\xi}_J-\bm{\xi}_J^{\res}$ to denote the extension of this vector to $I_{p_i}^c$ with coordinates in $I_{p_i}^c\backslash J$ equal to zero and where $\bm{\phi}_i''$ denotes the vector which coincides with $\bm{\xi}$ (and hence $\bm{\xi}^{\res}$) on $I\bigcap I_{p_i}^c$ and has remaining coordinates zero. Moreover, for notational convenience and to save space, we omit the ``$/p_i^{e_{p_i}}$'' following the bold vectors, although their entries are in fact taken in the quotient ring $\Z/p_i^{e_{p_i}}\Z$.

By the estimate in \eqref{eq: decoupling_sym}, we decouple the index sets $I\cap I_{p_i}^c$ and $I\backslash I_{p_i}^c$ for all $1\le i\le l$,  removing the quadratic terms and reducing the problem to a linear one. We now consider two circumstances. To handle the event $\bm{\xi}_J\equiv\bm{\xi}_J^{\res}\pmod {p_1}$, we have 
$$\mathbf{P}(\bm{\xi}_J\equiv\bm{\xi}_J^{\res}\pmod {p_1})\le(1-\epsilon)^{\#J}\le (1-\epsilon)^{n/1000-ln^{2/3}}=O_{a,\epsilon}(\exp(-\epsilon\Omega_{a}(n))).$$ 
Otherwise, let us assume $\bm{\xi}_J\not\equiv\bm{\xi}_J^{\res}\pmod {p_1}$. Under this assumption, when we extend the vector 
$$(M_n/p_1)_{I_{p_1}^c\times I_{p_1}^c}^{-1}(\bm{\xi}_J-\bm{\xi}^{\res}_J)/{p_1}$$ to a $n$-dimensional vector in $\F_{p_1}$ by padding it with $0$s, the resulting vector is orthogonal to the rows of $M_n/p_1$ with index in $I\backslash\bigcup_{i=1}^lI_{p_i}$. Notice that the index set $I\backslash\bigcup_{i=1}^lI_{p_i}$ does not intersect with $[\frac{n}{1000}]$, and
$$\#(I\backslash \bigcup_{i=1}^lI_{p_i})=n-\#J-\sum_{i=1}^l\#I_{p_i}\ge n-\frac{n}{1000}-ln^{2/3}.$$
Therefore, when $n$ is sufficiently large, by \ref{item: ortho nonsparse_sym}, we deduce that the vector
$$(M_n/p_1^{e_{p_1}})_{I_{p_1}^c\times I_{p_1}^c}^{-1}(\bm{\xi}_J-\bm{\xi}^{\res}_J)/p_1^{e_{p_1}}$$ has at least $\frac{n}{100}$ invertible entries. Furthermore, this vector has at least 
$$\frac{n}{100}-\#J-\#I_{p_1}^c\ge\frac{n}{100}-\frac{n}{1000}-n^{2/3}\ge\frac{n}{200}-n^{2/3}$$ invertible entries on $I\cap I_{p_1}^c$, i.e., the support of the random entries of $\bm{\phi}_1''$. Applying \Cref{lem: expectation of not concentrated}, we have 
\begin{align}\label{eq: nonsparse after decoupling}
\begin{split}
\left|\E_{\bm{\xi}_I}\left[\exp\left(\sum_{i=1}^l\frac{4\pi \sqrt{-1}}{p_i^{e_{p_i}}}b_i(\bm{\xi}_J-\bm{\xi}^{\res}_J)^T(M_n/p_i^{e_{p_i}})_{I_{p_i}^c\times I_{p_i}^c}^{-1}\bm{\phi}_i''\right)\right]\right|&\le\exp\left(-\frac{\epsilon}{a^2}\left(\frac{n}{200}-n^{2/3}\right)\right)\\
&=O_{a,\epsilon}(\exp(-\epsilon\Omega_a(n))).
\end{split}
\end{align}

\noindent{\bf Case 2. When $a$ is odd, and $b_i=0$ for all $1\le i\le l$. } In this case, we have
\begin{multline}
\text{LHS}\eqref{eq: small at nonzero vector_sym}=\Bigg|\E_{\bm{\xi},z}\exp\Bigg(\sum_{i=1}^l\frac{2\pi \sqrt{-1}}{p_i^{e_{p_i}}}\bm{\alpha}_i\cdot\Bigg((\bm{\xi}/p_i^{e_{p_i}})_{I_{p_i}}\\
-(M_n/p_i^{e_{p_i}})_{I_{p_i}\times I_{p_i}^c}(M_n/p_i^{e_{p_i}})_{I_{p_i}^c\times I_{p_i}^c}^{-1}(\bm{\xi}/p_i^{e_{p_i}})_{I_{p_i}^c}\Bigg)\Bigg)\Bigg|.
\end{multline}
Here, the random variable inside the exponent does not contain quadratic terms. Without loss of generality, we can let $\bm{\alpha}_1\ne 0$. In this case, let $k$ be the smallest valuation of all the entries of $\bm{\alpha}_1$. Then $0\le k<e_{p_1}$, and there exists 
$\bm{\beta}_1\in(\Z/p_1^{e_{p_1}}\Z)^{\#I_{p_1}}$, such that
$$\bm{\alpha}_1=p_1^k\bm{\beta_1},\bm{\beta}_1/p_1\ne 0.$$
Now, we consider the $n$-dimensional vector 
$$\left(-(\bm{\beta}_1/p_1),(M_n/p_1)_{I_{p_1}^c\times I_{p_1}^c}^{-1}(M_n/p_1)_{I_{p_1}^c\times I_{p_1}}(\bm{\beta}_1/p_1)\right)\in\F_p^n.$$ 
Here, we rearranged the entries: $-(\bm{\beta}_1/p_1)$ refers to the index set $I_{p_1}$, and $$(M_n/p_1)_{I_{p_1}^c\times I_{p_1}^c}^{-1}(M_n/p_1)_{I_{p_1}^c\times I_{p_1}}(\bm{\beta}_1/p_1)$$ 
refers to the index set $I_{p_1}^c$.
This vector is orthogonal to the rows of $M_n/p_1$ with index in $I_{p_1}^c$, which has cardinality $\ge n-n^{2/3}$. By \ref{item: ortho nonsparse_sym}, we deduce that the vector $$(M_n/p_1^{e_{p_1}})_{I_{p_1}^c\times I_{p_1}^c}^{-1}(M_n/p_1^{e_{p_1}})_{I_{p_1}^c\times I_{p_1}}\bm{\alpha}_1=p^k(M_n/p_1^{e_{p_1}})_{I_{p_1}^c\times I_{p_1}^c}^{-1}(M_n/p_1^{e_{p_1}})_{I_{p_1}^c\times I_{p_1}}\bm{\beta}_1$$ 
has at least $\frac{n}{100}-n^{2/3}$ nonzero entries. Applying \Cref{lem: expectation of not concentrated}, we have
$$\text{LHS}\eqref{eq: small at nonzero vector_sym}\le\exp(-\frac{\epsilon (n/100-n^{2/3})}{a^2})=O_{a,\epsilon}(\exp(-\epsilon\Omega_a(n))).$$

\noindent{\bf Case 3. When $a$ is even.} Without loss of generality, we can let $p_1=2$. Since $M_n$ satisfies $\mathcal{E}_n^{\text{\textdagger}}$, we have $\#I_{2}\le n^{1/4}$.
If $2b_1\ne 0$, then $\frac{4\pi\sqrt{-1}}{2^{e_2}}b_1$ is not a multiple of $2\pi\sqrt{-1}$, and we can still use the decoupling method in the first case to obtain the inequalities \eqref{eq: decoupling_sym} and \eqref{eq: nonsparse after decoupling}.
If $b_i\ne 0$ for some $2\le i\le l$, the proof also follows the same as in the first case. If $b_i=0$ for all $1\le i\le l$, the proof follows the same as in the second case. Therefore, the only remaining case is when $b_1\ne 0$ but $2b_1=b_2=\ldots=b_l=0$. In this case, we denote by $\bm{\zeta}_{I_2^c}\in(\Z/2^{e_2}\Z)^{n-\#I_2}$ the diagonal of $(M_n/2^{e_2})_{I_2^c\times I_2^c}^{-1}$. Notice that in the field $\F_2$, we have $\xi_i^2=\xi_i,2\xi_i\xi_j=0$. Therefore, the quadratic term degenerates to a linear form, i.e.,
$$b_1\left(z-(\bm{\xi}/2^{e_2})_{I_2^c}^T(M_n/2^{e_2})_{I_2^c\times I_2^c}^{-1}(\bm{\xi}/2^{e_2})_{I_2^c}\right)=b_1z-b_1\bm{\zeta}_{I_2^c}\cdot(\bm{\xi}/2^{e_2})_{I_2^c}.$$
Hence, we have
\begin{multline}
\text{LHS}\eqref{eq: small at nonzero vector_sym}=\Bigg|\E_{\bm{\xi},z}\exp\Bigg(\frac{2\pi \sqrt{-1}}{2^{e_2}}(b_1z-b_1\bm{\zeta}_{I_2^c}\cdot(\bm{\xi}/2^{e_2})_{I_2^c})\\+\sum_{i=1}^l\frac{2\pi \sqrt{-1}}{p_i^{e_{p_i}}}\bm{\alpha}_i\cdot((\bm{\xi}/p_i^{e_{p_i}})_{I_{p_i}}-(M_n/p_i^{e_{p_i}})_{I_{p_i}\times I_{p_i}^c}(M_n/p_i^{e_{p_i}})_{I_{p_i}^c\times I_{p_i}^c}^{-1}(\bm{\xi}/p_i^{e_{p_i}})_{I_{p_i}^c})\Bigg)\Bigg|.\\
\end{multline}
Similarly as in the second step, there are no quadratic terms to handle. Now, guaranteed by $\mathcal{E}_n^{\text{\textdagger}}$, the vector $-b_1\zeta_{I_2^c}+(M_n/p_i^{e_{p_i}})_{I_{p_i}^c\times I_{p_i}^c}^{-1}(M_n/p_i^{e_{p_i}})_{I_{p_i}^c\times I_{p_i}}\bm{\alpha_1}\in(\Z/2^{e_2}\Z)^{n-\#I_2}$ has at least $\frac{1}{2}n^{1/2}$ nonzero entries.
Applying \Cref{lem: expectation of not concentrated}, we have
$$\text{LHS}\eqref{eq: small at nonzero vector_sym}\le\exp\left(-\frac{\epsilon n^{1/2}}{2a^2}\right)=O_{a,\epsilon}(\exp(-\Omega_{a,\epsilon}(n^{1/2}))).$$

The above steps complete the proof.
\end{proof}

\begin{cor}\label{cor: integral universal transition_sym}
Let $\epsilon>0$ be a real number, and $\left(G^{(1)},\langle\cdot,\cdot\rangle_{G^{(1)}}\right),\left(G^{(2)},\langle\cdot,\cdot\rangle_{G^{(2)}}\right)$ be paired abelian groups of finite size. Let $P=\{p_1,\ldots,p_l\}$ be a finite set of primes that include all those that divide $\#G^{(1)}\#G^{(2)}$. Let $n\ge 1$, and $M_n\in\Sym_n(\Z)$ be a fixed matrix such that the following holds:
\begin{enumerate}
\item $(\Cok(M_n)_P,\langle\cdot,\cdot\rangle_P)\simeq\left(G^{(1)},\langle\cdot,\cdot\rangle_{G^{(1)}}\right)$.
\item $M_n$ satisfies the properties $\mathcal{E}_{n,p_i}^{\sym}$ for $1\le i\le l$. 
\item When $2\in P$, $M_n$ satisfies $\mathcal{E}_n^{\text{\textdagger}}$.
\end{enumerate}
Let $z,\xi_1,\ldots,\xi_n$ be independent $\epsilon$-balanced random integers. Then, we have
\begin{multline}\label{eq: probability isomorphic to G2}
\mathbf{P}\left(\left(\Cok
\left(
\begin{array}{c|c}
M_n & 
\begin{matrix}
\xi_1 \\
\vdots \\
\xi_n
\end{matrix}
\\ \hline
\begin{matrix}
\xi_1 & \cdots & \xi_n
\end{matrix}
& z
\end{array}
\right)_P,\langle\cdot,\cdot\rangle_P\right)\simeq(G^{(2)},\langle\cdot,\cdot\rangle_{G^{(2)}})\right)\\
=\prod_{i=1}^l\mathbf{P}\left(\left(G_{p_i}^{(2)},\langle\cdot,\cdot\rangle_{G_{p_i}^{(2)}}\right)\bigg|\left(G_{p_i}^{(1)},\langle\cdot,\cdot\rangle_{G_{p_i}^{(1)}}\right)\right)+O_{G^{(1)},G^{(2)},P,\epsilon}\left(\exp\left(-\Omega_{G^{(1)},G^{(2)},P,\epsilon}\left(n^{1/2}\right)\right)\right).
\end{multline}
Moreover, when $P$ does not contain the prime $2$, the error term on the last line can be further reduced to $O_{G^{(1)},G^{(2)},P,\epsilon}(\exp(-\epsilon\Omega_{G^{(1)},G^{(2)},P}(n)
))$.
\end{cor}

\begin{proof}
We will only prove the case $2\in P$, since the case $2\notin P$ can be deduced routinely. Let $a:=p_1^{e_{p_1}}\ldots p_l^{e_{p_l}}$, where for all $1\le i\le l$, we have
$$e_{p_i}=\begin{cases}
\max\{\Dep_{p_i}(G^{(1)}),\Dep_{p_i}(G^{(2)})\}+3 & p_i=2\\
\max\{\Dep_{p_i}(G^{(1)}),\Dep_{p_i}(G^{(2)})\}+1 & p_i>2.
\end{cases}.$$
Then, we have $\Cok_*^{(a)}(M_n/a)\simeq\left(G^{(1)},\langle\cdot,\cdot\rangle_*^{(a)}\right)$. Let $z',\xi_1',\ldots,\xi_n'$ be independent and uniformly distributed in $\{0,1,\ldots,a-1\}$. Denote by $\bm{\xi}:=(\xi_1,\ldots,\xi_n)$, and $\bm{\xi}':=(\xi_1',\ldots,\xi_n')$. Applying \Cref{prop: universal transition_sym}, we have 
\begin{align}
\begin{split}
\text{LHS}\eqref{eq: probability isomorphic to G2}&=\mathbf{P}\left(\Cok_*^{(a)}\begin{pmatrix}M_n/a & \bm{\xi}/a\\
(\bm{\xi}/a)^T & z/a
\end{pmatrix}
\simeq(G^{(2)},\langle\cdot,\cdot\rangle_*^{(a)})\right)\\
&=\mathbf{P}\left(\Cok_*^{(a)}\begin{pmatrix}M_n/a & \bm{\xi}'/a\\
(\bm{\xi}'/a)^T & z'/a
\end{pmatrix}
\simeq(G^{(2)},\langle\cdot,\cdot\rangle_*^{(a)})\right)\\
&+O_{a,\epsilon}(\exp(-\Omega_{a,\epsilon}(n^{1/2})))\\
&=\prod_{i=1}^l\mathbf{P}\left(\Cok_*^{(p_i^{e_{p_i}})}\begin{pmatrix}M_n/p_i^{e_{p_i}} & \bm{\xi}'/p_i^{e_{p_i}}\\
(\bm{\xi}'/p_i^{e_{p_i}})^T & z'/p_i^{e_{p_i}}
\end{pmatrix}
\simeq(G_{p_i}^{(2)},\langle\cdot,\cdot\rangle_*^{(p_i^{e_{p_i}})})\right)\\
&+O_{a,\epsilon}(\exp(-\Omega_{a,\epsilon}(n^{1/2})))\\
&=\prod_{i=1}^l\mathbf{P}\left(\left(G_{p_i}^{(2)},\langle\cdot,\cdot\rangle_{G_{p_i}^{(2)}}\right)\bigg|\left(G_{p_i}^{(1)},\langle\cdot,\cdot\rangle_{G_{p_i}^{(1)}}\right)\right)+O_{a,\epsilon}(\exp(-\Omega_{a,\epsilon}(n^{1/2}))).\\
&=\text{RHS}\eqref{eq: probability isomorphic to G2}.
\end{split}
\end{align}
Here, the last line is due to the fact that $a$ depends on $G^{(1)},G^{(2)},P$. 
\end{proof}

\begin{proof}[Proof of \Cref{thm: quotient simulate_sym}]
Suppose $a$ has prime factorization $a=p_1^{e_{p_1}}\cdots p_l^{e_{p_l}}$. Consider the exposure process 
$$\Cok_*^{(a)}(M_{n/20}/a),\ldots,\Cok_*^{(a)}(M_n/a)$$ 
obtained by iteratively revealing $M_t$ for $n/20\le t\le n$ by adding a new row and column each time and considering the resulting cokernel with quasi-pairing $\Cok_*^{(a)}(M_t/a)$. Also, we construct another exposure process 
$$\Cok_*^{(a)}(M_{n/20}'/a),\ldots,\Cok_*^{(a)}(M_n'/a),$$
where $M_{n/20}':=M_{n/20}$, and the matrices $M_{n/20+1}',\ldots,M_n'$ have entries in $\Z$ and are obtained by adding a uniformly distributed row and column in $\{0,1,\ldots,a-1\}$ every time.

As we will see, the random variable $\Cok_*^{(a)}(M_n'/a)$ plays an intermediate role between $\mu_\infty^{\sym,(a)}$ and the law of $\Cok_*^{(a)}(M_n/a)$: that is, we first bound the $L^1$-distance between the law of $\Cok_*^{(a)}(M_n/a)$ and $\Cok_*^{(a)}(M_n'/a)$, and then bound the $L^1$-distance between $\mu_\infty^{\sym,(a)}$ and the law of $\Cok_*^{(a)}(M_n'/a)$.

\noindent{\bf Step 1. The $L^1$-distance between the law of $\Cok_*^{(a)}(M_n/a)$ and $\Cok_*^{(a)}(M_n'/a)$. }
We claim that for all $t\in[n/20,n-1]$, we have
\begin{multline}\label{eq: increment of distance adding new row and column_sym}
D_{L_1}\left(\Cok_*^{(a)}(M_{t+1}/a),\Cok_*^{(a)}(M_{t+1}'/a)\right)-D_{L_1}\left(\Cok_*^{(a)}(M_t/a),\Cok_*^{(a)}(M_t'/a)\right)\\
=O_{a,\epsilon}(\exp(-\Omega_{a,\epsilon}(t^{1/2}))).
\end{multline}
Indeed, by \Cref{thm: not sparse} and \Cref{thm: not sparse revised p=2}, the matrix $M_t\in\Sym_t(\Z)$ satisfies the properties $\mathcal{E}_{t,p_1}^{\sym},\ldots,\mathcal{E}_{t,p_l}^{\sym}$ and $\mathcal{E}_t^{\text{\textdagger}}$ with probability at least
$$1-\sum_{i=1}^lO_{p_i}(\exp(-\epsilon\Omega_{p_i}(t)))-O_\epsilon(\exp(-\Omega_\epsilon(t^{1/2})))=1-O_{a,\epsilon}(\exp(-\Omega_{a,\epsilon}(t^{1/2}))).$$ 
Guaranteed by \Cref{prop: transition probability relies on cokernel_sym}, as long as $M_t$ satisfies these properties, we can apply \Cref{prop: universal transition_sym} so that \eqref{eq: increment of distance adding new row and column_sym} immediately follows. Since $M_{n/20}'/a=M_{n/20}/a$, we deduce that
$$D_{L_1}\left(\Cok_*^{(a)}(M_n/a),\Cok_*^{(a)}(M_n'/a)\right)=\sum_{t=n/20}^n O_{a,\epsilon}(\exp(-\Omega_{a,\epsilon}(t^{1/2})))=O_{a,\epsilon}(\exp(-\Omega_{a,\epsilon}(n^{1/2}))).$$
Furthermore, if $a$ is odd, we no longer need the assumption $\mathcal{E}_n^{\text{\textdagger}}$, and therefore the stronger bound
$$D_{L_1}\left(\Cok_*^{(a)}(M_n/a),\Cok_*^{(a)}(M_n'/a)\right)=O_{a,\epsilon}(\exp(-\epsilon\Omega_a(n)))$$
follows the same way.

\noindent{\bf Step 2. The $L^1$-distance between $\mu_\infty^{\sym,(a)}$ and the law of $\Cok_*^{(a)}(M_n'/a)$.} 
By \Cref{thm: not sparse}, the matrix $M_{n/20}'\in\Sym_{n/20}(\Z)$ satisfies the properties $\mathcal{E}_{n/20,p_1}^{\sym},\ldots,\mathcal{E}_{n/20,p_l}^{\sym}$ with probability at least
$$1-\sum_{i=1}^lO_{p_i}(\exp(-\epsilon\Omega_{p_i}(n/20)))=1-O_{a,\epsilon}(\exp(-\epsilon\Omega_a(n))).$$
As long as $M_{n/20}'$ satisfies these properties, we have
$\corank(M_{n/20}'/p_i)\le (n/20)^{2/3}\le n^{2/3}$ for all $1\le i\le l$. Moreover, the exposure process naturally factorizes over the primes $p_1,\ldots,p_l$ when we add new uniform rows and columns, since each newly revealed entry is uniform in $\mathbb Z/a\mathbb Z$, hence (by the Chinese remainder theorem) its reductions modulo $p_i^{e_{p_i}}$ are independent and uniform for $1\le i\le l$. Thus, conditioned on $M_{n/20}\in\Sym_{n/20}(\Z)$, the random variables
$$\Cok_*^{(p_i^{e_{p_i}})}(M_n'/p_i^{e_{p_i}}),\quad 1\le i\le l$$
are independent. Also, we have that for all $t\in[n/20,n-1]$ and $1\le i\le l$,
\begin{align}
\begin{split}
\mathbf{P}(\corank(M_{t+1}'/p_i)-\corank(M_t'/p_i)=-1\mid\corank(M_t'/p_i)\ge 2)&\ge 1-p_i^{-2},\\
\mathbf{P}(\corank(M_{t+1}'/p_i)-\corank(M_t'/p_i)=-1\mid\corank(M_t'/p_i)=1)&=1-p_i^{-1}.\\
\end{split}
\end{align}
Indeed, working modulo $p_i$, we choose an
invertible principal minor of maximal size, whose existence is guaranteed by
\Cref{lem: full rank principal minor}. Using paired row-column operations (i.e., congruence
transformations), we can eliminate the entries in the corresponding rows and columns outside this
minor, reducing to a block form with an invertible block and a remaining complement block of
size equal to the corank, from which the desired one-step transition estimate follows.

Therefore, when $\corank(M_{n/20}'/p_i)\le n^{2/3}$, we can apply a standard Hoeffding
large-deviation bound to deduce that with probability at least $1-O_{p_i}(\exp(-\Omega_
{p_i}(n)))$, there exists $\tau_i\in[n/20+1,n/2]$ such that $\corank(M_{\tau_i}'/p_i)=0$. Notice that for any fixed $\tau_i\in[n/20+1,n/2]$, when we condition on $\corank(M_{\tau_i}'/p_i)=0$, we have
$$\Cok_*^{(p_i^{e_{p_i}})}(M_n'/p_i^{e_{p_i}})\stackrel{d}{=}\Cok_*^{(p_i^{e_{p_i}})}(H_{n-\tau_i}/p_i^{e_{p_i}}),$$
where $H_{n-\tau_i}/p_i^{e_{p_i}}\in\Sym_{n-\tau_i}(\Z/p_i^{e_{p_i}}\Z)$ is uniformly distributed. Therefore, recalling the exponential convergence given in \Cref{prop: uniform model exponential convergence_sym}, we have
\begin{equation}\label{eq: Mn' and H_infty distance at p_i}
D_{L^1}\left(\mathcal{L}\left(\Cok_*^{(p_i^{e_{p_i}})}(M_n'/p_i^{e_{p_i}})\Bigg|\mathcal{E}_{n/20,p_i}^{\sym}\right),\mu_\infty^{\sym,(p_i^{e_{p_i}})}\right)=O_{p_i^{e_{p_i}}}(\exp(-\Omega_{p_i^{e_{p_i}}}(n))).
\end{equation}
Summing up the $O_{a,\epsilon}(\exp(-\epsilon\Omega_a(n)))$ probability lost when assuming $\mathcal{E}_{n/20,p_1}^{\sym},\ldots,\mathcal{E}_{n/20,p_l}^{\sym}$ for $M_{n/20}$, and the $L^1$-distance estimate \eqref{eq: Mn' and H_infty distance at p_i} for all $1\le i\le l$, we deduce that
$$D_{L_1}(\mathcal{L}(\Cok_*^{(a)}(M_n'/a)),\mu_\infty^{\sym,(a)})=O_{a,\epsilon}(\exp(-\epsilon\Omega_a(n))).$$

These two steps together provide the proof.
\end{proof}

As a byproduct, we also obtain the convergence of the joint distribution of corners, which is given by the following proposition.

\begin{prop}\label{prop: joint distribution_sym}
Let $\epsilon>0$ and $M_n$ be the same as in \Cref{thm: exponential convergence_sym}. Moreover, for all $1\le t\le n$, let $M_t$ be the $t\times t$ upper-left corner of $M_n$. Let $j\ge 1$ be a fixed integer, 
$$\left(G^{(1)},\langle\cdot,\cdot\rangle_{G^{(1)}}\right),\ldots,\left(G^{(j)},\langle\cdot,\cdot\rangle_{G^{(j)}}\right)$$ 
be finite paired abelian groups, and $P=\{p_1,\ldots,p_l\}$ be a finite set of primes that includes all those that divide $\#G^{(1)}\cdots\#G^{(j)}$. Then we have
\begin{multline}
\mathbf{P}\left(M_{n-j+i}\text{ is nonsingular},(\Cok(M_{n-j+i})_P,\langle\cdot,\cdot\rangle_P)\simeq\left(G^{(i)},\langle\cdot,\cdot\rangle_{G^{(i)}}\right)\quad\forall 1\le i\le j\right)\\
=\prod_{i_1=1}^l\left(\frac{\prod_{k\ge 1}(1-p_{i_1}^{1-2k})}{\Bigg|G_{p_{i_1}}^{(1)}\Bigg|\Bigg|\Aut\Bigg(G_{p_{i_1}}^{(1)},\langle\cdot,\cdot\rangle_{G_{p_{i_1}}^{(1)}}\Bigg)\Bigg|}\cdot\prod_{i_2=1}^{j-1}\mathbf{P}\left(\left(G^{(i_2+1)}_{p_{i_1}},\langle\cdot,\cdot\rangle_{G^{(i_2+1)}_{p_{i_1}}}\right)\Bigg|\left(G^{(i_2)}_{p_{i_1}},\langle\cdot,\cdot\rangle_{G^{(i_2)}_{p_{i_1}}}\right)\right)\right)\\
+O_{G^{(1)},\ldots,G^{(j)},P,\epsilon}(\exp(-\Omega_{G^{(1)},\ldots,G^{(j)},P,\epsilon}(n^{1/2}))).
\end{multline}
Moreover, when $P$ does not contain the prime $2$, the error term on the last line can be further reduced to $O_{G^{(1)},\ldots,G^{(j)},P,\epsilon}(\exp(-\epsilon\Omega_{G^{(1)},\ldots,G^{(j)},P}(n)
))$.
\end{prop}

\begin{proof}
We will only prove the case $2\in P$, since the case $2\notin P$ can be deduced routinely. It suffices to prove the stronger statement that for all $1\le j_0\le j$, 
\begin{multline}\label{eq: joint distribution for j0_sym}
\mathbf{P}\left(M_{n-j+i}\text{ is nonsingular},(\Cok(M_{n-j+i})_P,\langle\cdot,\cdot\rangle_P)\simeq\left(G^{(i)},\langle\cdot,\cdot\rangle_{G^{(i)}}\right)\quad\forall 1\le i\le j_0\right)\\
=\prod_{i_1=1}^l\left(\frac{\prod_{k\ge 1}(1-p_{i_1}^{1-2k})}{\Bigg|G_{p_{i_1}}^{(1)}\Bigg|\Bigg|\Aut\Bigg(G_{p_{i_1}}^{(1)},\langle\cdot,\cdot\rangle_{G_{p_{i_1}}^{(1)}}\Bigg)\Bigg|}\cdot\prod_{i_2=1}^{j_0-1}\mathbf{P}\left(\left(G^{(i_2+1)}_{p_{i_1}},\langle\cdot,\cdot\rangle_{G^{(i_2+1)}_{p_{i_1}}}\right)\Bigg|\left(G^{(i_2)}_{p_{i_1}},\langle\cdot,\cdot\rangle_{G^{(i_2)}_{p_{i_1}}}\right)\right)\right)\\
+O_{G^{(1)},\ldots,G^{(j_0)},P,\epsilon}(\exp(-\Omega_{G^{(1)},\ldots,G^{(j_0)},P,\epsilon}(n^{1/2}))).
\end{multline}
The proof proceeds by induction over $j_0$. The case $j_0=1$ can be directly deduced from \Cref{thm: exponential convergence_sym}. Now, suppose \eqref{eq: joint distribution for j0_sym} already holds for some $1\le j_0\le j-1$. By \Cref{thm: not sparse} and \Cref{thm: not sparse revised p=2}, with probability no less than $1-O_{P,\epsilon}(\exp(-\Omega_{P,\epsilon}(n^{1/2})))$, the matrix $M_{n-j+j_0}$ satisfies the properties
$\mathcal{E}_{n-j+j_0,p_1}^{\sym},\ldots,\mathcal{E}_{n-j+j_0,p_l}^{\sym},\mathcal{E}_{n-j+j_0}^{\text{\textdagger}}.$ Thus, applying \Cref{cor: integral universal transition_sym}, 
we have
\begin{multline}
\mathbf{P}\bigg(M_{n-j+i}\text{ is nonsingular},(\Cok(M_{n-j+i})_P,\langle\cdot,\cdot\rangle_P)\simeq\left(G^{(i)},\langle\cdot,\cdot\rangle_{G^{(i)}}\right)\quad\forall 1\le i\le j_0+1,\\
\text{and }M_{n-j+j_0}\text{ satisfies the properties } \mathcal{E}_{n-j+j_0,p_1}^{\sym},\ldots,\mathcal{E}_{n-j+j_0,p_l}^{\sym},\mathcal{E}_{n-j+j_0}^{\text{\textdagger}}\bigg)\\
=\prod_{i_1=1}^l\left(\frac{\prod_{k\ge 1}(1-p_{i_1}^{1-2k})}{\Bigg|G_{p_{i_1}}^{(1)}\Bigg|\Bigg|\Aut\Bigg(G_{p_{i_1}}^{(1)},\langle\cdot,\cdot\rangle_{G_{p_{i_1}}^{(1)}}\Bigg)\Bigg|}\cdot\prod_{i_2=1}^{j_0}\mathbf{P}\left(\left(G^{(i_2+1)}_{p_{i_1}},\langle\cdot,\cdot\rangle_{G^{(i_2+1)}_{p_{i_1}}}\right)\Bigg|\left(G^{(i_2)}_{p_{i_1}},\langle\cdot,\cdot\rangle_{G^{(i_2)}_{p_{i_1}}}\right)\right)\right)\\
+O_{G^{(1)},\ldots,G^{(j_0+1)},P,\epsilon}(\exp(-\Omega_{G^{(1)},\ldots,G^{(j_0+1)},P,\epsilon}(n^{1/2}))).
\end{multline}
Therefore, the statement \eqref{eq: joint distribution for j0_sym} also holds for $j_0+1$. This completes the proof.
\end{proof}

%% file: Proof_of_the_alternating_case.tex
\section{Proof of the alternating case}\label{sec: Proof of the alternating case}

\begin{thm}\label{thm: quotient simulate_alt}
Let $\epsilon>0$ be a real number. For all $n\ge 1$, let $A_n\in\Alt_n(\Z)$ be the same as in \Cref{thm: exponential convergence_alt}, and $a\ge 2$ be an integer. Then we have
$$D_{L^1}\left(\mathcal{L}(\Cok(A_{2n}/a)),\mu_{2\infty}^{\alt,(a)}\right)=O_{a,\epsilon}(\exp(-\epsilon\Omega_a(n))),$$
$$D_{L^1}\left(\mathcal{L}(\Cok(A_{2n+1}/a)),\mu_{2\infty+1}^{\alt,(a)}\right)=O_{a,\epsilon}(\exp(-\epsilon\Omega_a(n))).$$
Here, $\mu_{2\infty}^{\alt,(a)},\mu_{2\infty+1}^{\alt,(a)}$ are the same as in \Cref{prop: uniform model exponential convergence_alt}. 
\end{thm}

Before we return to the proof of \Cref{thm: quotient simulate_alt}, we see how it implies \Cref{thm: exponential convergence_alt}.

\begin{proof}[Proof of \Cref{thm: exponential convergence_alt}, assuming \Cref{thm: quotient simulate_alt}]
Following \Cref{prop: simplest matrix_alt}, when $G\notin\mathcal{S}_P$, it is clear that
$$\mathbf{P}(\Cok(A_{2n})_P\simeq G)=\mathbf{P}(\Cok_{\tors}(A_{2n+1})_P\simeq G)=0.$$
Therefore, for the rest of this proof, we assume $G\in\mathcal{S}_P$. Let $a:=p_1^{e_{p_1}}\cdots p_l^{e_{p_l}}$ be the same as in \Cref{cor: globally same pairing_alt}. In this case, for all $1\le i\le l$ and $n\ge 1$,
$$\mu_{2n}^{\alt,(p_i^{e_{p_i}})}\left(G_{p_i}\right)=\mathbf{P}\left(\Cok(A_{2n}^{(p_i')})\simeq G_{p_i}\right),$$
$$\mu_{2n+1}^{\alt,(p_i^{e_{p_i}})}\left(G_{p_i}\oplus(\Z/p_i^{e_{p_i}}\Z)\right)=\mathbf{P}\left(\Cok_{\tors}(A_{2n+1}^{(p_i')})\simeq G_{p_i}\right).$$
Here, $A_{2n}^{(p_i')}\in\Alt_{2n}(\Z_{p_i}),A_{2n+1}^{(p_i')}\in\Alt_{2n+1}(\Z_{p_i})$ are Haar-distributed. For the even size case, applying \cite[Theorem 3.9]{bhargava2015modeling}, we have
\begin{align}\label{eq: uniform model even size_alt}
\begin{split}
\mu_{2\infty}^{\alt,(a)}\left(G\right)&=\prod_{i=1}^l\mu_{2\infty}^{\alt,(p_i^{e_{p_i}})}\left(G_{p_i}\right)\\
&=\prod_{i=1}^l\lim_{n\rightarrow\infty}\mathbf{P}\left(\Cok(A_{2n}^{(p_i')})\simeq G_{p_i}\right)\\
&=\prod_{i=1}^l\frac{|G_{p_i}|}{|\Sp(G_{p_i})|}\prod_{k\ge 1}(1-p_i^{1-2k})\\
&=\frac{|G|}{|\Sp(G)|}\prod_{i=1}^l\prod_{k\ge 1}(1-p_i^{1-2k}).
\end{split}
\end{align}
For the odd size case, following \cite[Section 4]{nguyen2025local}, which gives a detailed exposition of \cite[Theorem 3.11]{bhargava2015modeling}, we have
\begin{align}\label{eq: uniform model odd size_alt}
\begin{split}
\mu_{2\infty+1}^{\alt,(a)}\left(G\oplus(\Z/a\Z)\right)&=\prod_{i=1}^l\mu_{2\infty+1}^{\alt,(p_i^{e_{p_i}})}\left(G_{p_i}\oplus(\Z/p_i^{e_{p_i}}\Z)\right)\\
&=\prod_{i=1}^l\lim_{n\rightarrow\infty}\mathbf{P}\left(\Cok_{\tors}(A_{2n+1}^{(p_i')})\simeq G_{p_i}\right)\\
&=\prod_{i=1}^l\frac{1}{|\Sp(G_{p_i})|}\prod_{k\ge 1}(1-p_i^{-1-2k})\\
&=\frac{1}{|\Sp(G)|}\prod_{i=1}^l\prod_{k\ge 1}(1-p_i^{-1-2k}).
\end{split}
\end{align}

Therefore, due to our assumption in \Cref{thm: quotient simulate_alt}, we have
$$\left|\mathbf{P}\left(\Cok(A_{2n})_P\simeq G\right)-\frac{|G|}{|\Sp(G)|}\prod_{i=1}^l\prod_{k\ge 1}(1-p_i^{1-2k})\right|=O_{a,\epsilon}(\exp(-\epsilon\Omega_a(n))),$$
$$\left|\mathbf{P}\left(\Cok(A_{2n+1})_P\simeq G\right)-\frac{1}{|\Sp(G)|}\prod_{i=1}^l\prod_{k\ge 1}(1-p_i^{-1-2k})\right|=O_{a,\epsilon}(\exp(-\epsilon\Omega_a(n))).$$
Notice that the integer $a$ defined in \Cref{cor: globally same pairing_alt} depends on $G,P$. Therefore, we can rewrite $O_{a,\epsilon}(\exp(-\epsilon\Omega_a(n)))$ as $O_{G,P,\epsilon}(\exp(-\epsilon\Omega_{G,P}(n)))$. This completes the proof.
\end{proof}

Now we return to our proof of \Cref{thm: quotient simulate_alt}. Let us start from the following proposition.

\begin{prop}\label{prop: universal transition_alt}
Let $\epsilon>0$ be a real number, and $a\ge 2$ be an integer with prime factorization $a=p_1^{e_{p_1}}\cdots p_l^{e_{p_l}}$. Suppose $A_n\in\Alt_n(\Z)$ is a fixed matrix that satisfies $\mathcal{E}_{n,p_i}^{\alt}$ for $1\le i\le l$. Let $\xi_1,\ldots,\xi_n$ be independent $\epsilon$-balanced random integers, and $\xi_1',\ldots,\xi_n'$ be independent and uniformly distributed in $\{0,1,\ldots,a-1\}$. Let $\bm{\xi}=(\xi_1,\ldots,\xi_n),\bm{\xi}'=(\xi_1',\ldots,\xi_n')\in\Z^n$. Then we have
$$
D_{L^2}\left(\Cok
\begin{pmatrix}A_n/a & \bm{\xi}/a\\
-(\bm{\xi}/a)^T & 0
\end{pmatrix}
,\Cok\begin{pmatrix}A_n/a & \bm{\xi}'/a\\
-(\bm{\xi}'/a)^T & 0
\end{pmatrix}\right)=O_{a,\epsilon}(\exp(-\epsilon\Omega_{a}(n))).\\
$$
\end{prop}

\begin{proof}
For all $1\le i\le l$, by \Cref{lem: full rank principal minor}, there exists $I_{p_i}\subset[n]$ with $\# I_{p_i}=n-\rank(A_n/{p_i})$, such that $(A_n/p_i)_{I_{p_i}^c\times I_{p_i}^c}$ has full rank. By \ref{item: large rank_alt}, we have $\#I_{p_i}\le n^{2/3}$. 

Now, we treat the $p
_i$-part of the cokernel with its quasi-pairing, where we reduce all matrix entries modulo $p_i^{e_{p_i}}$. By paired row-column operations, we can use the invertible matrix $(A_n/p_i^{e_{p_i}})_{I_{p_i}^c\times I_{p_i}^c}$ to eliminate the components of $\bm{\xi}/p_i^{e_{p_i}}$ and $(\bm{\xi}/p_i^{e_{p_i}})^T$ indexed from $I_{p_i}^c$, while the alternating structure is preserved. In this way, the new added column $\bm{\xi}/p_i^{e_{p_i}}$ will be transformed into the form
$$\bm{w}_i:=
(\bm{\xi}/p_i^{e_{p_i}})_{I_{p_i}}-(A_n/p_i^{e_{p_i}})_{I_{p_i}\times I_{p_i}^c}(A_n/p_i^{e_{p_i}})_{I_{p_i}^c\times I_{p_i}^c}^{-1}(\bm{\xi}/p_i^{e_{p_i}})_{I_{p_i}^c}\in(\Z/p_i^{e_{p_i}}\Z)^{\#I_{p_i}}.\\
$$
Here, we ignore and skip the components with index $I_{p_i}^c$ because they have been eliminated to zero. When $i$ runs over all integers from $1$ to $l$, we obtain a combination of vectors $(\bm{w}_i)_{1\le i\le l}$,
and this combination is a random variable in $\left((\Z/p_i^{e_{p_i}}\Z)^{\#I_{p_i}}\right)_{1\le i\le l}$. If we replace $\xi_1,\ldots,\xi_n$ by the input $\xi_1',\ldots,\xi_n'$, we can obtain a combination of vectors $(\bm{w}_i')_{1\le i\le l}\in\left((\Z/p_i^{e_{p_i}}\Z)^{\#I_{p_i}}\right)_{1\le i\le l}$ in the similar way. It is clear that $(\bm{w}_i')_{1\le i\le l}$ is uniformly distributed. We only need to prove that $$D_{L^2}\left((\bm{w}_i)_{1\le i\le l},(\bm{w}_i')_{1\le i\le l}\right)=
O_{a,\epsilon}(\exp(-\epsilon\Omega_a(n))).
$$
To achieve this, applying Parseval's identity, we have
$$D_{L^2}^2\left((\bm{w}_i)_{1\le i\le l},(\bm{w}_i')_{1\le i\le l}\right)=\prod_{i=1}^lp_i^{-e_{p_i}\#I_{p_i}}\sum\left|\E_{\bm{\xi},z}\exp\left(\sum_{i=1}^l\frac{2\pi \sqrt{-1}}{p_i^{e_{p_i}}}\bm{\alpha}_i\cdot\bm{w}_i\right)\right|^2.$$
Here, the sum on the right hand side ranges through the nonzero elements in the set
$$
\mathcal{V}:=\{(\bm{\alpha}_i)_{1\le i\le l}: \bm{\alpha}_i:=(\alpha_{i,j_i})_{j_i\in I_{p_i}}\in(\Z/p_i^{e_{p_i}}\Z)^{\#I_{p_i}}\}.
$$
Therefore, it suffices to show that for all nonzero $(\bm{\alpha}_i)_{1\le i\le l}\in\mathcal{V}$, we have
\begin{multline}\label{eq: small at nonzero vector_alt}
\left|\E_{\bm{\xi},z}\exp\left(\sum_{i=1}^l\frac{2\pi \sqrt{-1}}{p_i^{e_{p_i}}}\bm{\alpha}_i\cdot\left((\bm{\xi}/p_i^{e_{p_i}})_{I_{p_i}}-(A_n/p_i^{e_{p_i}})_{I_{p_i}\times I_{p_i}^c}(A_n/p_i^{e_{p_i}})_{I_{p_i}^c\times I_{p_i}^c}^{-1}(\bm{\xi}/p_i^{e_{p_i}})_{I_{p_i}^c}\right)\right)\right|\\=
O_{a,\epsilon}(\exp(-\epsilon\Omega_a(n))).
\end{multline}
Without loss of generality, we can let $\bm{\alpha}_1\ne 0$. In this case, let $k$ be the smallest valuation of all the entries of $\bm{\alpha}_1$. Then $0\le k<e_{p_1}$, and there exists 
$\bm{\beta}_1\in(\Z/p_1^{e_{p_1}}\Z)^{\#I_{p_1}}$, such that
$$\bm{\alpha}_1=p_1^k\bm{\beta_1},\bm{\beta}_1/p_1\ne 0.$$
Now, we consider the $n$-dimensional vector 
$$\left(-(\bm{\beta}_1/p_1),(A_n/p_1)_{I_{p_1}^c\times I_{p_1}^c}^{-1}(A_n/p_1)_{I_{p_1}^c\times I_{p_1}}(\bm{\beta}_1/p_1)\right)\in\F_p^n.$$ 
Here, we rearranged the entries: $-(\bm{\beta}_1/p_1)$ refers to the index set $I_{p_1}$, and $$(A_n/p_1)_{I_{p_1}^c\times I_{p_1}^c}^{-1}(A_n/p_1)_{I_{p_1}^c\times I_{p_1}}(\bm{\beta}_1/p_1)$$ 
refers to the index set $I_{p_1}^c$.
This vector is orthogonal to the rows of $A_n/p_1$ with index in $I_{p_1}^c$, which has cardinality $\ge n-n^{2/3}$. By \ref{item: ortho nonsparse_alt}, we deduce that the vector $$(A_n/p_1^{e_{p_1}})_{I_{p_1}^c\times I_{p_1}^c}^{-1}(A_n/p_1^{e_{p_1}})_{I_{p_1}^c\times I_{p_1}}\bm{\alpha}_1=p^k(A_n/p_1^{e_{p_1}})_{I_{p_1}^c\times I_{p_1}^c}^{-1}(A_n/p_1^{e_{p_1}})_{I_{p_1}^c\times I_{p_1}}\bm{\beta}_1$$ 
has at least $\frac{n}{100}-n^{2/3}$ nonzero entries. Applying \Cref{lem: expectation of not concentrated}, we have
$$\text{LHS}\eqref{eq: small at nonzero vector_alt}\le\exp(-\frac{\epsilon (n/100-n^{2/3})}{a^2})=O_{a,\epsilon}(\exp(-\epsilon\Omega_a(n))).$$
This completes the proof.
\end{proof}

\begin{cor}\label{cor: integral universal transition_alt}
Let $\epsilon>0$ be a real number, and $n\ge 1$. Let $P=\{p_1,\ldots,p_l\}$ be a finite set of primes, and $G^{(1)},G^{(2)}\in\mathcal{S}_p$. Let $\xi_1,\xi_2,\ldots$ be independent $\epsilon$-balanced random integers.
\begin{enumerate}
\item (From even to odd) Let $A_{2n}\in\Alt_{2n}(\Z)$ be a fixed matrix that satisfies $\mathcal{E}_{2n,p_i}^{\alt}$ for $1\le i\le l$, and $\Cok(A_{2n})_P\simeq G^{(1)}$. Then, we have
\begin{multline}\label{eq: even to odd isomorphic to G2_alt}
\mathbf{P}\left(\Cok_{\tors}\left(
\begin{array}{c|c}
A_{2n} & 
\begin{matrix}
\xi_1 \\
\vdots \\
\xi_{2n}
\end{matrix}
\\ \hline
\begin{matrix}
-\xi_1 & \cdots & -\xi_{2n}
\end{matrix}
& 0
\end{array}
\right)_P\simeq G^{(2)},\corank\left(
\begin{array}{c|c}
A_{2n} & 
\begin{matrix}
\xi_1 \\
\vdots \\
\xi_{2n}
\end{matrix}
\\ \hline
\begin{matrix}
-\xi_1 & \cdots & -\xi_{2n}
\end{matrix}
& 0
\end{array}
\right)=1\right)\\
=\prod_{i=1}^l\mathbf{P}\left(G_{p_i}^{(2)},\text{ odd }\bigg| G_{p_i}^{(1)},\text{ even}\right)+O_{G^{(1)},G^{(2)},P,\epsilon}(\exp(-\epsilon\Omega_{G^{(1)},G^{(2)},P}(n))).
\end{multline}
\item (From odd to even) Let $A_{2n+1}\in\Alt_{2n+1}(\Z)$ be a fixed matrix that satisfies $\mathcal{E}_{2n+1,p_i}^{\alt}$ for $1\le i\le l$, and $\Cok_{\tors}(A_{2n+1})_P\simeq G^{(1)},\corank(A_{2n+1})=1$. Then, we have
\begin{multline}\label{eq: odd to even isomorphic to G2_alt}
\mathbf{P}\left(\Cok\left(
\begin{array}{c|c}
A_{2n+1} & 
\begin{matrix}
\xi_1 \\
\vdots \\
\xi_{2n+1}
\end{matrix}
\\ \hline
\begin{matrix}
-\xi_1 & \cdots & -\xi_{2n+1}
\end{matrix}
& 0
\end{array}
\right)_P\simeq G^{(2)}\right)\\
=\prod_{i=1}^l\mathbf{P}\left(G_{p_i}^{(2)},\text{ even }\bigg| G_{p_i}^{(1)},\text{ odd}\right)+O_{G^{(1)},G^{(2)},P,\epsilon}(\exp(-\epsilon\Omega_{G^{(1)},G^{(2)},P}(n))).
\end{multline}
\end{enumerate}
\end{cor}

\begin{proof}
We will only prove the ``from even to odd'' case, since the ``from odd to even'' can be deduced routinely. Let $a:=p_1^{e_{p_1}}\ldots p_l^{e_{p_l}}$, where $e_{p_i}=\max\{\Dep_{p_i}(G^{(1)}),\Dep_{p_i}(G^{(2)})\}+1$ for all $1\le i\le l$. Let $\xi_1',\xi_2',\ldots$ be independent and uniformly distributed in $\{0,1,\ldots,a-1\}$. Denote by $\bm{\xi}:=(\xi_1,\ldots,\xi_{2n})$, and $\bm{\xi}':=(\xi_1',\ldots,\xi_{2n}')$. Applying \Cref{prop: universal transition_alt}, we have 
\begin{align}
\begin{split}
\text{LHS}\eqref{eq: even to odd isomorphic to G2_alt}&=\mathbf{P}\left(\Cok\begin{pmatrix}A_{2n}/a & \bm{\xi}'/a\\
-(\bm{\xi}'/a)^T & 0
\end{pmatrix}
\simeq G^{(2)}\oplus(\Z/a\Z)\right)\\
&=\mathbf{P}\left(\Cok\begin{pmatrix}A_{2n}/a & \bm{\xi}/a\\
-(\bm{\xi}/a)^T & 0
\end{pmatrix}
\simeq G^{(2)}\oplus(\Z/a\Z)\right)+O_{a,\epsilon}(\exp(-\epsilon\Omega_a(n)))\\
&=\prod_{i=1}^l\mathbf{P}\left(\Cok\begin{pmatrix}A_{2n}/p_i^{e_{p_i}} & \bm{\xi}'/p_i^{e_{p_i}}\\
-(\bm{\xi}'/p_i^{e_{p_i}})^T & 0
\end{pmatrix}
\simeq G^{(2)}\oplus(\Z/a\Z)\right)+O_{a,\epsilon}(\exp(-\epsilon\Omega_a(n)))\\
&=\prod_{i=1}^l\mathbf{P}\left(G_{p_i}^{(2)},\text{ odd }\bigg| G_{p_i}^{(1)},\text{ even}\right)+O_{a,\epsilon}(\exp(-\epsilon\Omega_a(n))).\\
&=\prod_{i=1}^l\mathbf{P}\left(G_{p_i}^{(2)},\text{ odd }\bigg| G_{p_i}^{(1)},\text{ even}\right)+O_{G^{(1)},G^{(2)},P,\epsilon}(\exp(-\epsilon\Omega_{G^{(1)},G^{(2)},P}(n))).
\end{split}
\end{align}
Here, the last line is due to the fact that $a$ depends on $G^{(1)},G^{(2)},P$.
\end{proof}

\begin{proof}[Proof of \Cref{thm: quotient simulate_alt}]
We will only prove the even size case, since the odd size case can be deduced routinely. Suppose $a$ has prime factorization $a=p_1^{e_{p_1}}\cdots p_l^{e_{p_l}}$. Consider the exposure process 
$$\Cok(A_{n/20}/a),\ldots,\Cok(A_{2n}/a)$$ 
obtained by iteratively revealing $A_t$ for $n/20\le t\le 2n$ by adding a new row and column each time and considering the resulting cokernel with quasi-pairing $\Cok(A_t/a)$. Also, we construct another exposure process 
$$\Cok(A_{n/20}'/a),\ldots,\Cok(A_{2n}'/a),$$
where $A_{n/20}':=A_{n/20}$, and the matrices $A_{n/20+1}',\ldots,A_{2n}'$ have entries in $\Z$ and are obtained by adding a uniformly distributed row and column in $\{0,1,\ldots,a-1\}$ every time.

As we will see, the random variable $\Cok(A_{2n}'/a)$ plays an intermediate role between $\mu_{2\infty}^{\alt,(a)}$ and the law of $\Cok(A_{2n}/a)$: that is, we first bound the $L^1$-distance between the law of $\Cok(A_{2n}/a)$ and $\Cok(A_{2n}'/a)$, and then bound the $L^1$-distance between $\mu_{2\infty}^{\alt,(a)}$ and the law of $\Cok(A_{2n}'/a)$.

\noindent{\bf Step 1. The $L^1$-distance between the law of $\Cok(A_{2n}/a)$ and $\Cok(A_{2n}'/a)$. }
We claim that for all $t\in[n/20,2n-1]$, we have
\begin{equation}\label{eq: increment of distance adding new row and column_alt}
D_{L_1}\left(\Cok(A_{t+1}/a),\Cok(A_{t+1}'/a)\right)-D_{L_1}\left(\Cok(A_t/a),\Cok(A_t'/a)\right)\\
=O_{a,\epsilon}(\exp(-\epsilon\Omega_{a}(t))).
\end{equation}
Indeed, by \Cref{thm: not sparse}, the matrix $A_t\in\Alt_t(\Z)$ satisfies the properties $\mathcal{E}_{t,p_1}^{\alt},\ldots,\mathcal{E}_{t,p_l}^{\alt}$ with probability at least
$$1-\sum_{i=1}^lO_{p_i}(\exp(-\epsilon\Omega_{p_i}(t)))=1-O_{a,\epsilon}(\exp(-\epsilon\Omega_{a}(t))).$$ 
As long as $A_t$ satisfies these properties, we can apply \Cref{prop: universal transition_alt} so that \eqref{eq: increment of distance adding new row and column_alt} immediately follows. Since $A_{n/20}'/a=A_{n/20}/a$, we deduce that
$$D_{L_1}\left(\Cok(A_{2n}/a),\Cok(A_{2n}'/a)\right)=\sum_{t=n/20}^{2n} O_{a,\epsilon}(\exp(-\epsilon\Omega_a(t)))=O_{a,\epsilon}(\exp(-\epsilon\Omega_a(n))).$$

\noindent{\bf Step 2. The $L^1$-distance between $\mu_{2\infty}^{\alt,(a)}$ and the law of $\Cok(A_{2n}'/a)$.} 
By \Cref{thm: not sparse}, the matrix $A_{n/20}'\in\Alt_{n/20}(\Z)$ satisfies the properties $\mathcal{E}_{n/20,p_1}^{\alt},\ldots,\mathcal{E}_{n/20,p_l}^{\alt}$ with probability at least
$$1-\sum_{i=1}^lO_{p_i}(\exp(-\epsilon\Omega_{p_i}(n/20)))=1-O_{a,\epsilon}(\exp(-\epsilon\Omega_a(n))).$$
As long as $A_{n/20}'$ satisfies these properties, we have
$\corank(A_{n/20}'/p_i)\le (n/20)^{2/3}\le n^{2/3}$ for all $1\le i\le l$. Moreover, the exposure process naturally factorizes over the primes $p_1,\ldots,p_l$ when we add new uniform rows and columns, since each newly revealed entry is uniform in $\mathbb Z/a\mathbb Z$, hence (by the Chinese remainder theorem) its reductions modulo $p_i^{e_{p_i}}$ are independent and uniform for $1\le i\le l$. Thus, conditioned on $A_{n/20}\in\Alt_{n/20}(\Z)$, the random variables
$$\Cok(A_{2n}'/p_i^{e_{p_i}}),\quad 1\le i\le l$$
are independent. Also, we have that for all $t\in[n/20,n-1]$ and $1\le i\le l$,
\begin{align}
\begin{split}
\mathbf{P}(\corank(A_{t+1}'/p_i)-\corank(A_t'/p_i)=-1\mid\corank(A_t'/p_i)\ge 2)&\ge 1-p_i^{-2},\\
\mathbf{P}(\corank(A_{t+1}'/p_i)-\corank(A_t'/p_i)=-1\mid\corank(A_t'/p_i)=1)&=1-p_i^{-1}.\\
\end{split}
\end{align}
Indeed, working modulo $p_i$, we choose an
invertible principal minor of maximal size, whose existence is guaranteed by
\Cref{lem: full rank principal minor}. Using paired row-column operations (i.e., congruence
transformations), we can eliminate the entries in the corresponding rows and columns outside this
minor, reducing to a block form with an invertible block and a remaining complement block of
size equal to the corank, from which the desired one-step transition estimate follows.

Therefore, when $\corank(A_{n/20}'/p_i)\le n^{2/3}$, we can apply a standard Hoeffding
large-deviation bound to deduce that with probability at least $1-O_{p_i}(\exp(-\Omega_
{p_i}(n)))$, there exists an even integer $\tau_i\in[n/20+1,n/2]$ such that $\corank(A_{\tau_i}'/p_i)=0$. Notice that for any fixed even integer $\tau_i\in[n/20+1,n/2]$, when we condition on $\corank(A_{\tau_i}'/p_i)=0$, we have
$$\Cok(A_{2n}'/p_i^{e_{p_i}})\stackrel{d}{=}\Cok(H_{2n-\tau_i}/p_i^{e_{p_i}}),$$
where $H_{2n-\tau_i}/p_i^{e_{p_i}}\in\Alt_{2n-\tau_i}(\Z/p_i^{e_{p_i}}\Z)$ is uniformly distributed. Therefore, recalling the exponential convergence given in \Cref{prop: uniform model exponential convergence_alt}, we have
\begin{equation}\label{eq: An' and H_infty distance at p_i}
D_{L^1}\left(\mathcal{L}\left(\Cok(A_{2n}'/p_i^{e_{p_i}})\Bigg|\mathcal{E}_{n/20,p_i}^{\alt}\right),\mu_{2\infty}^{\alt,(p_i^{e_{p_i}})}\right)=O_{p_i^{e_{p_i}}}(\exp(-\Omega_{p_i^{e_{p_i}}}(n))).
\end{equation}
Summing up the $O_{a,\epsilon}(\exp(-\epsilon\Omega_a(n)))$ probability lost when assuming $\mathcal{E}_{n/20,p_1}^{\alt},\ldots,\mathcal{E}_{n/20,p_l}^{\alt}$ for $A_{n/20}$, and the $L^1$-distance estimate \eqref{eq: An' and H_infty distance at p_i} for all $1\le i\le l$, we deduce that
$$D_{L_1}(\mathcal{L}(\Cok(A_{2n}'/a)),\mu_{2\infty}^{\alt,(a)})=O_{a,\epsilon}(\exp(-\epsilon\Omega_a(n))).$$

These two steps together provide the proof.
\end{proof}

As a byproduct, we also obtain the convergence of the joint distribution of corners, which is given by the following proposition.

\begin{prop}\label{prop: joint distribution_alt}
Let $\epsilon>0$ and $A_{2n}$ be the same as in \Cref{thm: exponential convergence_alt}. Moreover, for all $1\le t\le 2n$, let $A_t$ be the $t\times t$ upper-left corner of $A_{2n}$. Let $j\ge 1$ be a fixed integer, $P=\{p_1,\ldots,p_l\}$ be a finite set of primes, and
$$G^{(0)},\ldots,G^{(2j)}\in\mathcal{S}_P.$$ 
Then we have 
\begin{align}
\begin{split}
&\mathbf{P}\Bigg(\Cok(A_{2n-2j})_P\simeq G^{(0)},\Cok_{\tors}(A_{2n-2j+1})_P\simeq G^{(1)},\corank(A_{2n-2j+1})=1,\cdots,\\
&\Cok(A_{2n-2})_P\simeq G^{(2j-2)},\Cok_{\tors}(A_{2n-1})_P\simeq G^{(2j-1)},\corank(A_{2n-1})=1,\Cok(A_{2n})_P\simeq G^{(2j)}\Bigg)\\
&=\prod_{i_1=1}^l\Bigg(\frac{\prod_{k\ge 1}(1-p_{i_1}^{-1-2k})}{\left|\Sp\left(G_{p_{i_1}}^{(0)}\right)\right|}\cdot\prod_{i_2=0}^{j-1}\mathbf{P}\left(G^{(2i_2+2)}_{p_{i_1}},\text{ even }\bigg|G^{(2i_2+1)}_{p_{i_1}},\text{ odd}\right)\\
&\cdot\mathbf{P}\left(G^{(2i_2+1)}_{p_{i_1}},\text{ odd }\bigg|G^{(2i_2)}_{p_{i_1}},\text{ even}\right)\Bigg)+O_{G^{(0)},\ldots,G^{(2j)},P,\epsilon}(\exp(-\epsilon\Omega_{G^{(0)},\ldots,G^{(2j)},P}(n))).
\end{split}
\end{align}
\end{prop}

\begin{proof}
It suffices to prove the stronger statement that for all $1\le j_0\le j$, 
\begin{align}\label{eq: joint distribution for j0_alt} 
\begin{split}
&\mathbf{P}\Bigg(\Cok(A_{2n-2j})_P\simeq G^{(0)},\Cok_{\tors}(A_{2n-2j+1})_P\simeq G^{(1)},\corank(A_{2n-2j+1})=1,\cdots,\\
&\Cok_{\tors}(A_{2n-2j+2j_0-1})_P\simeq G^{(2j_0-1)},\corank(A_{2n-2j+2j_0-1})=1,\Cok(A_{2n-2j+2j_0})_P\simeq G^{(2j_0)}\Bigg)\\
&=\prod_{i_1=1}^l\Bigg(\frac{\prod_{k\ge 1}(1-p_{i_1}^{-1-2k})}{\left|\Sp\left(G_{p_{i_1}}^{(0)}\right)\right|}\cdot\prod_{i_2=0}^{j_0-1}\left(G^{(2i_2+2)}_{p_{i_1}},\text{ even }\bigg|G^{(2i_2+1)}_{p_{i_1}},\text{ odd}\right)\\
&\cdot\mathbf{P}\left(G^{(2i_2+1)}_{p_{i_1}},\text{ odd }\bigg|G^{(2i_2)}_{p_{i_1}},\text{ even}\right)\Bigg)+O_{G^{(0)},\ldots,G^{(2j_0)},P,\epsilon}(\exp(-\epsilon\Omega_{G^{(0)},\ldots,G^{(2j_0)},P}(n))).
\end{split}
\end{align}
The proof proceeds by induction over $j_0$. The case $j_0=1$ can be directly deduced from \Cref{thm: exponential convergence_alt}. Now, suppose \eqref{eq: joint distribution for j0_alt} already holds for some $1\le j_0\le j-1$. By \Cref{thm: not sparse}, with probability no less than $1-O_{P,\epsilon}(\exp(-\epsilon\Omega_{P}(n)))$, the matrix $A_{2n-2j+2j_0}$ satisfies the properties
$\mathcal{E}_{2n-2j+2j_0,p_1}^{\alt},\ldots,\mathcal{E}_{2n-2j+2j_0,p_l}^{\alt}$, and the matrix $A_{2n-2j+2j_0+1}$ satisfies the properties
$\mathcal{E}_{2n-2j+2j_0+1,p_1}^{\alt},\ldots,\mathcal{E}_{2n-2j+2j_0+1,p_l}^{\alt}$. Thus, applying \Cref{cor: integral universal transition_alt} once for the transition from even to odd and once from odd to even, 
we have
\begin{align}
\begin{split}
&\mathbf{P}\Bigg(\Cok(A_{2n-2j})_P\simeq G^{(0)},\Cok_{\tors}(A_{2n-2j+1})_P\simeq G^{(1)},\corank(A_{2n-2j+1})=1,\cdots,\\
&\Cok_{\tors}(A_{2n-2j+2j_0+1})_P\simeq  G^{(2j_0+1)},\corank(A_{2n-2j+2j_0+1})=1,\Cok(A_{2n-2j+2j_0+2})_P\simeq  G^{(2j_0+2)},\\
&\text{$A_{2n-2j+2j_0}$ satisfies the properties
$\mathcal{E}_{2n-2j+2j_0,p_1}^{\alt},\ldots,\mathcal{E}_{2n-2j+2j_0,p_l}^{\alt}$},\\
&\text{and $A_{2n-2j+2j_0+1}$ satisfies the properties
$\mathcal{E}_{2n-2j+2j_0+1,p_1}^{\alt},\ldots,\mathcal{E}_{2n-2j+2j_0+1,p_l}^{\alt}$} \Bigg)\\
&=\prod_{i_1=1}^l\Bigg(\frac{\prod_{k\ge 1}(1-p_{i_1}^{-1-2k})}{\left|\Sp\left(G_{p_{i_1}}^{(0)}\right)\right|}\cdot\prod_{i_2=0}^{j_0}\mathbf{P}\left(G^{(2i_2+2)}_{p_{i_1}},\text{ even }\bigg|G^{(2i_2+1)}_{p_{i_1}},\text{ odd}\right)\\
&\cdot\mathbf{P}\left(G^{(2i_2+1)}_{p_{i_1}},\text{ odd }\bigg|G^{(2i_2)}_{p_{i_1}},\text{ even}\right)\Bigg)+O_{G^{(0)},\ldots,G^{(2j_0+2)},P,\epsilon}(\exp(-\epsilon\Omega_{G^{(0)},\ldots,G^{(2j_0+2)},P}(n))).
\end{split}
\end{align}
Therefore, the statement \eqref{eq: joint distribution for j0_alt} also holds for $j_0+1$. This completes the proof. 
\end{proof}

\begin{rmk}
\Cref{prop: joint distribution_alt} shows that the joint law of the corner data admits universal asymptotics.
In particular, fixing a single prime $p$ (i.e.\ taking $P=\{p\}$), as $n\to\infty$ we obtain the
same limiting distribution as in the Haar model over $\mathbb Z_p$.
This Haar setting can be viewed as the non-archimedean analogue of the aGUE corners process studied in \cite{shen2024non}:
the transition kernel appears in \cite[Theorem 1.3]{shen2024non}, the limiting law is given in~\cite[(1.7)]{shen2024non},
and the full joint distribution of corners is described in \cite[Theorem 1.2]{shen2024non}.
Taken together, these results identify the corner process as a Hall--Littlewood process.

It is also worth mentioning that, in the present paper, we use the forward exposure process obtained by adjoining one new row and column at each step,
whereas \cite{shen2024non} studies the same joint distribution from the reverse direction:
one first samples the cokernel of the largest matrix and then applies the backward transition obtained by deleting a row and column at each step.
\end{rmk}

%% file: Index_of_notations.tex
\appendix
\crefalias{section}{appendix}
\begin{appendix}

\section{Index of notations}\label{sec: index of notations}

\renewcommand{\arraystretch}{1.15}
\setlength{\tabcolsep}{6pt}

\begin{longtable}{@{}>{\raggedright\arraybackslash}p{0.28\textwidth}
                  >{\raggedright\arraybackslash}p{0.68\textwidth}@{}}

\label{tab:notation-index}\\
\toprule
\textbf{Notation} & \textbf{Meaning} \\
\midrule
\endfirsthead

\toprule
\textbf{Notation} & \textbf{Meaning} \\
\midrule
\endhead

\midrule
\multicolumn{2}{r}{\emph{(Continued on next page)}}\\
\endfoot

\bottomrule
\endlastfoot

\multicolumn{2}{@{}l}{\textbf{Basic sets and counting}}\\
$[n]$ & The set $\{1,2,\dots,n\}$.\\
$I^c$ & Complement of an index set $I\subseteq [n]$ in $[n]$.\\
$|\,\cdot\,|,\ \#(\cdot)$ & Cardinality / order of a finite set or group.\\[2pt]

\multicolumn{2}{@{}l}{\textbf{Probability and convergence}}\\
$\mathbb{P}(\cdot)$ & Probability of an event.\\
$\mathbb{E}(\cdot)$ & Expectation.\\
$\mathcal{L}(\cdot)$ & Law (distribution) of a random variable.\\
$\mathbb{P}(\cdot\mid\cdot)$, $\mathbb{E}(\cdot\mid\cdot)$, $\mathcal{L}(\cdot\mid\cdot)$ & Conditional probability / expectation / law.\\
$D_{L^q}(\mu,\nu)$ & $L^q$-distance between measures $\mu,\nu$ on a countable set.\\
$D_{L^q}(\mathcal{X},\mathcal{Y})$ & Abbreviation for $D_{L^q}(\mathcal{L}(\mathcal{X}),\mathcal{L}(\mathcal{Y}))$.\\
$\mathcal{X} \stackrel{d}{=} \mathcal{Y}$ & $\mathcal{X}$ and $\mathcal{Y}$ have the same distribution.\\
$\mathcal{X}_n \stackrel{d}{\to} \mathcal{X}$ & Weak convergence in distribution.\\[2pt]

\multicolumn{2}{@{}l}{\textbf{Asymptotic notation}}\\
$\exp(\cdot)$ & Exponential function.\\
$f = O_S(g)$ & $|f|\le K|g|$ for some $K>0$ depending on parameter set $S$.\\
$f = \Omega_S(g)$ & $f\ge c|g|$ for some $c>0$ depending on $S$.\\[2pt]

\multicolumn{2}{@{}l}{\textbf{$p$-adic / finite field notation}}\\
$\mathbb{Q}_p$ & $p$-adic field.\\
$\mathrm{val}(\cdot)$ & $p$-adic valuation $\mathbb{Q}_p^\times\to\mathbb{Z}$.\\
$\mathbb{Z}_p$ & $p$-adic integers.\\
$\mathbb{F}_p$ & Finite field $\mathbb{Z}/p\mathbb{Z}$.\\[2pt]

\multicolumn{2}{@{}l}{\textbf{Matrices and linear algebra}}\\
$\mathrm{Mat}_n(R)$ & $n\times n$ matrices over a ring $R$.\\
$\GL_n(R)$ & $n\times n$ invertible matrices over a ring $R$.\\
$\mathrm{Sym}_n(\mathbb{Z})$ & Symmetric integer matrices $M$ with $M^T=M$.\\
$\mathrm{Alt}_n(\mathbb{Z})$ & Alternating integer matrices $A$ with $A^T=-A$.\\
$\mathrm{supp}(\bm{\xi})$ & Support $\{i\in[n]: \xi_i\neq 0\}$ of a vector $\bm{\xi}=(\xi_1,\dots,\xi_n)$.\\
$\mathrm{wt}(\bm{\xi})$ & Weight $\#\,\mathrm{supp}(\bm{\xi})$.\\
$\xi_J$ & Subvector of $\xi$ restricted to coordinates indexed by $J\subseteq[n]$.\\
$B_{I\times J}$ & Submatrix of $B$ with row indices $I$ and column indices $J$.\\
$B/a$ & Reduction of an integer matrix $B$ modulo $a$ (entries in $\mathbb{Z}/a\mathbb{Z}$).\\
$B/a'$ & Further reduction when $a'\mid a$.\\
$\xi/a$ & Reduction of a vector modulo $a$.\\
$\xi\cdot \eta$ & Dot product of two vectors of the same length.\\
$B_n \cong B_n'$ & Congruence: $\exists U\in \mathrm{GL}_n(\mathbb{Z}_p)$ s.t.\ $U B_n U^T=B_n'$.\\[2pt]

\multicolumn{2}{@{}l}{\textbf{Cokernels and pairings}}\\
$\mathrm{Cok}(Q)$ & Cokernel $\mathbb{Z}^n/Q\mathbb{Z}^n$ for $Q\in \mathrm{Mat}_n(\mathbb{Z})$.\\
$\mathrm{Cok}(M)_p$ & $p$-primary / localized cokernel $\mathbb{Z}_p^n/M\mathbb{Z}_p^n$ (when working over $\mathbb{Z}_p$).\\
$\mathrm{Cok}(M)_P$ & $P$-primary part for a finite set of primes $P$.\\
$\langle\cdot,\cdot\rangle$ & Canonical perfect symmetric pairing on $\mathrm{Cok}(M)$ for nonsingular symmetric $M$.\\
$\mathrm{Aut}(G,\langle\cdot,\cdot\rangle_G)$ & Automorphisms of $G$ preserving the given pairing.\\
$\mathrm{Cok}^{(a)}_\ast(H)$ & Cokernel with quasi-pairing induced by a symmetric matrix $H$ over $\mathbb{Z}/a\mathbb{Z}$.\\
$(G,\langle\cdot,\cdot\rangle^{(a)}_\ast)$ & The equivalence class under $\mathrm{Cok}^{(a)}_\ast$ associated to a paired abelian group $(G,\langle\cdot,\cdot\rangle_G)$.\\
$\mathcal{S}_p$ & The set of finite abelian $p$-groups in $\mathcal{S}$, the set of squares of abelian groups.\\
$\mathcal S_P$ &
For a finite set of primes $P$, the set of finite abelian $P$-primary groups, equivalently
$\mathcal S_P=\prod_{p\in P}\mathcal S_p$.\\
$\mathrm{Sp}(G)$ & The automorphism group of $G$ preserving a perfect alternating pairing on $G$ (when $G\in \mathcal{S}_p$).\\
\parbox[t]{\linewidth}{\raggedright
$\mathbf{P}\!\left(\left(G^{(2)},\langle\cdot,\cdot\rangle_{G^{(2)}}\right)\right.$\\
\hspace*{1.2em}$\left.\mid\left(G^{(1)},\langle\cdot,\cdot\rangle_{G^{(1)}}\right)\right)$
}
&
\parbox[t]{\linewidth}{\raggedright
(Symmetric case) One-step transition probability between paired $p$-groups when
adding a new Haar-distributed row and column
}\\
$\mathbf{P}\left(G^{(2)},\text{ odd}\Big|G^{(1)},\text{ even}\right)$ & (Alternating case) One-step transition probability from even to odd size, between $p$-groups in $\mathcal{S}_p$ when
adding a new Haar-distributed row and column\\
$\mathbf{P}\left(G^{(2)},\text{ even}\Big|G^{(1)},\text{ odd}\right)$ & (Alternating case) One-step transition probability from odd to even size, between $p$-groups in $\mathcal{S}_p$ when
adding a new Haar-distributed row and column\\
$\mu_{\infty}^{\sym,(a)}$ & The limiting probability measure of $\mu_n^{\sym,(a)}:=\mathcal{L}(\mathrm{Cok}(a)_\ast(H_n))$ for $H_n\in \mathrm{Sym}_n(\mathbb{Z}/a\mathbb{Z})$ uniform.\\
$\mu_{2\infty}^{\alt,(a)}$ & The limiting probability measure of $\mu_{2n}^{\alt,(a)}:=\mathcal{L}(\mathrm{Cok}(H_{2n}))$ for $H_{2n}\in \mathrm{Alt}_{2n}(\mathbb{Z}/a\mathbb{Z})$ uniform.\\
$\mu_{2\infty+1}^{\alt,(a)}$ & The limiting probability measure of $\mu_{2n+1}^{\alt,(a)}:=\mathcal{L}(\mathrm{Cok}(H_{2n+1}))$ for $H_{2n+1}\in \mathrm{Alt}_{2n+1}(\mathbb{Z}/a\mathbb{Z})$ uniform.\\
[2pt]

\multicolumn{2}{@{}l}{\textbf{Finite abelian groups and partitions}}\\
$\Y$ & Set of partitions $\lambda=(\lambda_1,\lambda_2,\dots)$, eventually $0$, with $\lambda_1\ge\lambda_2\ge\cdots$.\\
$|\lambda|$ & Size $\sum_{i\ge 1}\lambda_i$.\\
$\mathrm{Len}(\lambda)$ & Length $\#\{i:\lambda_i>0\}$.\\
$G_p$ & Sylow $p$-subgroup of a finite abelian group $G$.\\
$G_P$ & $P$-primary subgroup $\bigoplus_{p\in P} G_p$.\\
$\mathrm{Dep}_p(G)$ & $p$-depth: smallest $e_p$ s.t.\ $p^{e_p} G_p=0$.\\
$\mathrm{Len}_p(G)$ & $p$-length: number of cyclic factors in the type of $G_p$.\\[2pt]

\multicolumn{2}{@{}l}{\textbf{Coding-theory quantities}}\\
$H_p(x)$ & $p$-ary entropy function (defined for $0<x<1$).\\
$\mathrm{Vol}_p(n,t)$ & Size of the Hamming ball in $\mathbb{F}_p^n$ of radius $t$.\\

\end{longtable}

\end{appendix}